\documentclass{conm-p-l}

\usepackage{amssymb}

\usepackage{amsmath,array,amsthm,mathtools,slashed,stmaryrd,tikz-cd, bbm, hyperref}
\tikzset{shorten <>/.style={shorten >=#1,shorten <=#1}}
\usepackage{tkz-euclide}
\usetikzlibrary{decorations.markings,arrows}

\theoremstyle{definition}

\theoremstyle{remark}

\numberwithin{equation}{section}

\theoremstyle{plain}
\newtheorem{thm}{Theorem}[section]

\newtheorem{cor}[thm]{Corollary}
\newtheorem{lem}[thm]{Lemma}
\newtheorem{prop}[thm]{Proposition}
\newtheorem{conjecture}[thm]{Conjecture}
 
\theoremstyle{remark}
\newtheorem{rmk}[thm]{Remark}
\newtheorem{ex}[thm]{Example} 
\newtheorem{defn}[thm]{Definition}
\newtheorem{notation}[thm]{Notation}
 
\newtheorem*{thm*}{Theorem}

\newsavebox{\pullback}
\sbox\pullback{%
\begin{tikzpicture}%
\draw (0,0) -- (1ex,0ex);%
\draw (1ex,0ex) -- (1ex,1ex);%
\end{tikzpicture}}

\RequirePackage{mathrsfs}

\newcommand{\Fun}{{\sf {Fun}}}

\newcommand{\calw}{\mathcal{W}}

\newcommand{\Or}{{\sf{or}}}
\newcommand{\one}{\mathbbm{1}}

\newcommand{\CS}{{\rm CS}}
\newcommand{\SW}{{\rm SW}}
\newcommand{\pre}{{\rm pre}}

\newcommand{\Pic}{{\sf Pic}}

\newcommand{\beq}{\begin{eqnarray}}
\newcommand{\eeq}{\end{eqnarray}}
\newcommand{\bp}{\begin{proof}[Proof]}
\newcommand{\ep}{\end{proof}}

\newcommand{\rmO}{{\rm{O}}}
\newcommand{\rmH}{{\rm H}}

\newcommand{\K}{{\mathbb K}}
\newcommand{\Spin}{{\rm Spin}}
\newcommand{\String}{{\rm String}}

\newcommand{\SO}{{\rm SO}}

\newcommand{\cG}{\mathcal{G}}

\newcommand{\cP}{\mathcal{P}}

\newcommand{\Gerbe}{{\sf Gerbe}}

\newcommand{\pt}{{\rm pt}}

\newcommand{\N}{{\mathbb{N}}}
\newcommand{\id}{{{\rm id}}}

\newcommand{\Z}{{\mathbb{Z}}}
\newcommand{\Aut}{{\sf Aut}}
\newcommand{\Bun}{{\sf Bun}}

\newcommand{\sq}{/\!\!/}

\newcommand{\nohat}[2]{{\sf{Bun}}_{#1}^\text{pre}(#2)}
\newcommand{\nohatG}[1]{{\sf{Bun}}_{\cG}^\text{pre}(#1)}
\newcommand{\LieGrpd}{\sf{LieGrpd}}

\newcommand{\twocommute}{\ensuremath{\rotatebox[origin=c]{30}{$\Rightarrow$}}}
\newcommand{\twocommuter}{\ensuremath{\rotatebox[origin=c]{-150}{$\Rightarrow$}}}

\newcommand\EquivTo{\xrightarrow{
   \,\smash{\raisebox{-0.65ex}{\ensuremath{\scriptstyle\sim}}}\,}}

\newcommand{\mf}\mathfrak
\newcommand{\mc}\mathcal
\newcommand{\mb}\mathbb

\usepackage[textwidth=25mm, textsize=tiny]{todonotes}
\usepackage[normalem]{ulem}
\newcommand{\out}{\bgroup\markoverwith
	{\textcolor{red}{\rule[.6ex]{3pt}{0.6pt}}}\ULon}

\begin{document}

\title{Flat principal 2-group bundles and flat string structures}


\author[D. Berwick-Evans]{Daniel Berwick-Evans}
\address{Department of Mathematics, University of Illinois Urbana-Champaign, Champaign, Illinois 61820}
\curraddr{}
\email{}
\thanks{The first author was supported in part by NSF Grant DMS 2205835}

\author[E. Cliff]{Emily Cliff}
\address{Department of Mathematics, Universit\'e de Sherbrooke, Sherbrooke, Quebec J1K 2R1, Canada}
\curraddr{}
\email{}
\thanks{The second author was supported in part by NSERC under Discovery Grant RGPIN-2022-04104 and FRQNT Research Support for New Academics under Grant 327706}

\author[L. Murray]{Laura Murray}
\address{Department of Mathematics \& Computer Science, Providence College, Providence, Rhode Island 02918}
\curraddr{}
\email{}
\thanks{The third author was supported in part by NSF Grant DMS 2316646}

\author[A. Nakade]{Apurva Nakade}
\address{Department of Mathematics, Northwestern University, Evanston, Illinois 60208}
\curraddr{}
\email{}
\thanks{}

\author[E. Phillips]{Emma Phillips}
\address{Department of Mathematics \& Statistics, University of New Hampshire, Durham, New Hampshire 03824}
\curraddr{}
\email{}
\thanks{}


\date{}

\begin{abstract}
For a weak 2-group $\mathcal{G}$, we construct a bicategory of flat $\mathcal{G}$-bundles over differentiable stacks as a localization of a functor bicategory. This description is amenable to explicit geometric constructions. For example, we show that flat $\mathcal{G}$-bundles can be described in terms of ordinary $G$-bundles together with a trivialization of a certain 2-gerbe. This specializes to a characterization of flat string structures on vector bundles over differentiable stacks. 
\end{abstract}

\maketitle

\section{Introduction}

The \emph{string group} is the smooth categorical central extension
\beq\label{eq:stringgroupdef}
1\to \pt\sq U(1)\to \String_n\to \Spin_n\to 1,\qquad n\in \N
\eeq
classified by the image of the fractional Pontryagin class under the isomorphism $\frac{p_1}{2}\in \rmH^4(B\Spin_n;\Z)\simeq \rmH^3_{\rm SM}(\Spin_n;U(1))$, where $\rmH_{\rm SM}$ denotes Segal--Mitchison cohomology~\cite{SP11}. 
The smooth extension~\eqref{eq:stringgroupdef} is in the bicategory of smooth 2-groups, defined as group objects in smooth stacks. Every smooth 2-group has an underlying discrete 2-group~\cite[page 28]{SP11}. In particular, the discrete version of ~\eqref{eq:stringgroupdef} yields the discrete string group $\String_n^\delta$ (see~\S\ref{sec:disstring}). 

\begin{defn}[{Definition~\ref{defn:stringstruc}}]\label{defn:string}
A \emph{flat string structure} on a flat vector bundle~$V$ is a lift 
\beq\begin{tikzpicture}[baseline=(basepoint)];
\node (B) at (4,0) {$\pt\sq \String_n^\delta$};
\node (D) at (4,-1) {$\pt\sq \rmO_n^\delta$};
\node (C) at (0,-1) {$X$};
\draw[->,dashed] (C) to (B);
\draw[->] (B) to (D);
\draw[->] (C) to node [below] {$V$} (D);
\path (0,-.75) coordinate (basepoint);
\end{tikzpicture}\label{eq:stringlift1}
\eeq
where $\rmO_n^\delta$ is the orthogonal group $\rmO_n$ with the discrete topology, and the arrow labeled $V$ is the map classifying the frame bundle of the flat vector bundle~$V$. 
\end{defn}

String structures were born out of physicists' work on the 2-dimensional supersymmetric sigma model~\cite{WittenDirac,Killingback}, and have been codified in various mathematical ways, e.g., \cite{ST04,Bunke}. Among these approaches, our Definition~\ref{defn:string} is a flat version of Schommer-Pries's notion in~\cite{SP11} using the language of 2-group principal bundles. Another approach to string structures is Waldorf's definition involving trivializations of the Chern--Simons 2-gerbe~\cite{Waldorfstring}. The following connects these a priori distinct approaches in the flat case; we note that flat string structures are the objects of a bicategory.

\begin{thm}[Theorem~\ref{thm:flatstringdata}]\label{thm2}The bicategory of flat string structures on a vector bundle $V\to X$ is equivalent to the bicategory of trivializations of $\CS_V\to X$ where $\CS_V$ is the flat Chern--Simons 2-gerbe of $V$.
\end{thm}

Theorem~\ref{thm2} can be situated within a general framework of flat 2-group bundles. 
To explain this, note that the continuous homomorphism 
$\Spin^\delta_n\to \Spin_n$ induces a map 
\beq\label{eq:discretecohomo}
\rmH^3_{\rm SM}(\Spin_n;U(1))\to \rmH^3(\Spin_n^\delta;U(1))
\eeq
from Segal--Mitchison cohomology to group cohomology of the discrete group $\Spin_n^\delta$ (see Lemma~\ref{lem:SMtoGRP}). The discrete extension underlying~\eqref{eq:stringgroupdef} is classified by an ordinary group cohomology class in $\rmH^3(\Spin_n^\delta;U(1))$, given by the image of $\frac{p_1}{2}$ under~\eqref{eq:discretecohomo}. 
More generally, a discrete group~$G$, discrete abelian group~$A$, and 3-cocycle $\alpha\in Z^3(G;A)$ for trivial $G$-action on $A$ determine a discrete 2-group $\cG$ sitting in the  categorical central extension 
\beq\label{eq:extension}
1\to \pt\sq A\to \cG \to G\to 1
\eeq
classified by $[\alpha]\in \rmH^3(G;A)$. This gives a generalization of~\eqref{eq:stringlift1} for the moduli of lifts of a flat $G$-bundle $P$ to flat $\cG$-bundles,
\beq\begin{tikzpicture}[baseline=(basepoint)];
\node (B) at (4,0) {$\pt\sq \cG$};
\node (D) at (4,-1) {$\pt\sq G.$};
\node (C) at (0,-1) {$X$};
\draw[->,dashed] (C) to (B);
\draw[->] (B) to (D);
\draw[->] (C) to node [below] {$P$} (D);
\path (0,-.75) coordinate (basepoint);
\end{tikzpicture}\label{eq:2grpbundlelifting}
\eeq
Letting $P$ vary, the collection of dashed arrows is the bicategory $\Bun_\cG(X)$ of $\cG$-bundles on a differentiable stack~$X$, and there is a forgetful functor $\Bun_\cG(X)\to \Bun_G(X)$ that extracts the  ordinary $G$-bundle $P\to X$. 

\begin{rmk}
    Since the bundles considered in~\eqref{eq:2grpbundlelifting} use the discrete topologies on $G$ and $A$, they are necessarily flat. However, to underscore the parallel with Definition~\ref{defn:string}, we often call a lift~\eqref{eq:2grpbundlelifting} a \emph{flat} $\cG$-bundle.  
\end{rmk}

We prove the following generalization of Theorem~\ref{thm2}.

\begin{thm}[Theorem~\ref{thm:2gerbetriv}]\label{thm1}
A flat $\cG$-bundle on $X$ is equivalent to a $G$-bundle $P\to X$ with a trivialization of the 2-gerbe $\lambda_{P,\alpha}$ determined by $P$ and $\alpha$. 
\end{thm}

\begin{rmk} \label{rmk:strictweak}
For 2-groups with trivial associator (e.g., \emph{strict} 2-groups), the 2-gerbe $\lambda_{P_\alpha}$ is canonically trivial; a change of trivialization is then a choice of (1-)gerbe. This agrees with previous descriptions of strict 2-group bundles, e.g., \cite{BaezSchreiber},~\cite[\S7.1]{NikolausWaldorf}. Theorem~\ref{thm1} therefore reveals a subtle distinction: the theory of strict 2-group bundles reduces to the geometry of gerbes, whereas (weak) 2-group bundles involve the geometry of 2-gerbes and their trivializations. 
\end{rmk}

\subsection{String structures from homotopy theory} 
For a topological space~$X$, isomorphism classes of $n$-dimensional (metrized) vector bundles $V\to X$ are in bijection with homotopy classes of maps $X\to B\rmO_n$ to the classifying space $B\rmO_n$. This leads to homotopical notions of orientations, spin structures, and string structures corresponding to lifts,
\beq\begin{tikzpicture}[baseline=(basepoint)];
\node (AA) at (4,2) {$B\String_n$};
\node (A) at (4,1) {$B\Spin_n$};
\node (AAA) at (7,1) {${\rm K}(\Z,4)$};
\node (B) at (4,0) {$B\SO_n$};
\node (BB) at (7,0) {${\rm K}(\Z/2,2)$};
\node (D) at (4,-1) {$B\rmO_n$};
\node (DD) at (7,-1) {${\rm K}(\Z/2,1)$.};
\node (C) at (0,-1) {$X$};
\draw[->,dashed] (C) to (B);
\draw[->,dashed] (C) to (A);
\draw[->,dashed] (C) to (AA);
\draw[->] (B) to (D);
\draw[->] (AA) to (A);
\draw[->] (A) to (B);
\draw[->] (D) to node [above] {$w_1$} (DD);
\draw[->] (B) to node [above] {$w_2$} (BB);
\draw[->] (A) to node [above] {$\frac{p_1}{2}$} (AAA);
\draw[->] (C) to node [below] {$V$} (D);
\path (0,.5) coordinate (basepoint);
\end{tikzpicture}\label{eq:homotopicallifting}
\eeq
In the above diagram, 
\beq\label{eq:Whitehead}
\dots \to B\String_n\to B\Spin_n\to B\SO_n\to B\rmO_n
\eeq
is the Whitehead tower of $B\rmO_n$, i.e., each successive space in the tower~\eqref{eq:Whitehead} is obtained by killing the lowest-degree homotopy group of the previous space. In the cases of interest, these spaces admit equivalent descriptions as bundles classified by the maps $w_1,w_2,$ and $p_1/2$ given by the Stiefel--Whitney classes and the fractional Pontryagin class. 

In~\eqref{eq:homotopicallifting}, the notation $BG$ for various groups $G$ is not entirely precise, as we now explain. By construction, the spaces in the tower~\eqref{eq:Whitehead} are only defined up to homotopy equivalence. The notation is justified by the fact that the spaces $B\SO_n$ and $B\Spin_n$ have the same homotopy type as the classifying spaces of the Lie groups $\SO_n$ and $\Spin_n$, respectively. This leads to differential-geometric descriptions of lifts in~\eqref{eq:homotopicallifting} in terms of $\SO_n$- and $\Spin_n$-principal bundles on~$X$. As an added bonus, this geometry leads to the ``correct" definition of $G$-equivariant orientations and $G$-equivariant spin structures via $G$-actions on the corresponding $\SO_n$- and $\Spin_n$-principal bundles. 

On the other hand, the homotopy type $B\String_n$ in~\eqref{eq:Whitehead} does not admit a description in terms of the classifying space of a Lie group~\cite[page~3]{SP11}. It does admit various models as the classifying space of an infinite-dimensional group~\cite{Stolz,ST04,BL04,BCSS,Andre,NSW,Waldorf,Bunk}, but these descriptions are not well-suited to differential geometric applications. Instead, the most promising model for $B\String_n$ is as the classifying space of the (weak) smooth 2-group $\String_n$~\cite[\S1]{SP11}. This identifies $B\String_n$ with a classifying space for higher categorical bundles, namely 2-group principal bundles for the smooth 2-group~$\String_n$. One of the primary goals of this paper is to gain a better understanding of the geometry of such $\String_n$-bundles and string structures via Theorem~\ref{thm2} and Conjecture~\ref{conjecture} below. 

\subsection{String structures and topological modular forms}

One reason for pursuing a deeper understanding of string geometry comes from the Stolz--Teichner program~\cite{ST04,ST11}. It is expected that a certain (yet to be defined mathematically) class of 2-dimensional field theories furnishes a geometric model for the cohomology theory of topological modular forms (TMF). Physical reasoning predicts that string structures encode anomaly cancellation data for these field theories~\cite{WittenDirac, Killingback}. The most compelling evidence for Stolz and Teichner's proposal is that string structures are orientation data for TMF: vector bundles with string structure have Thom isomorphisms in TMF. However, a differential geometric description of TMF has evaded discovery for over 30 years. This is partly because there are so many seemingly distinct mathematical approaches to 2-dimensional quantum field theory (e.g.,~\cite{Segal_Elliptic,SegalCFT,HU2004325,ST04,ST11,baas_dundas_rognes_2004,DouglasHenriques,Costello1,Costello2}), and partly because it has been difficult to construct a map from geometric objects to~TMF. 

An important development in this vein is Lurie's universal property that characterizes TMF~\cite[\S5.5]{Lurie}. The key idea is that of a \emph{2-equivariant cohomology theory}, which is a higher-categorical equivariant refinement involving 2-groups acting on topological stacks~\cite[\S5.1-5.3]{Lurie}. Lurie sketches the construction of a 2-equivariant refinement of TMF, and further describes how TMF satisfies a universal property among 2-equivariant cohomology theories. This suggests a possible road to constructing a comparison map between TMF and field theories, utilizing Lurie's universal property. Furthermore, 2-equivariance is naturally accommodated in only a few current paradigms for mathematical quantum field theory. This cuts down the zoo of choices, and makes is harder to guess the ``wrong" definition of quantum field theory that incorporates 2-equivariance. 

The geometric origins of 2-equivariance are undoubtedly wrapped up in the string group and string structures, understood in terms of the geometry of 2-group principal bundles. 
Combined with the geometry of anomaly cancellation in physics,
these ideas come together in a (as-yet undeveloped) theory of \emph{equivariant} string structures, e.g., for vector bundles defined on quotient stacks $M\sq G$. Definition~\ref{defn:string} accommodates such examples in the case of flat vector bundles, and Conjecture~\ref{conjecture} below indicates the anticipated generalization to the non-flat case.

\subsection{2-group principal bundles}
The majority of the previous literature on 2-group principal bundles concerns \emph{strict} 2-groups, implying that the associator is trivial e.g., see~\cite{BaezSchreiber,BaezStevenson,Wockel,SchreiberWaldorf,NikolausWaldorf}. Crucially, however, the string group~\eqref{eq:stringgroupdef} is not strict. 
Hence, a robust theory of string structures requires a better developed theory of 2-group principal bundles for weak 2-groups. In the previous approaches to principal $\cG$-bundles for a weak 2-group $\cG$ (e.g.,~\cite{Bartels,LurieTopos,SP11,NikolausSchreiberStevenson1,NikolausSchreiberStevenson2}), a solid global theory was developed, but the bicategory of $\cG$-bundles on a manifold $M$ is somewhat difficult to access directly.

A main goal of this paper is to develop an explicit presentation of the bicategory $\Bun_\cG(X)$ of flat principal $\cG$-bundles over a differentiable stack $X$ for $\cG$ a weak 2-group. This presentation specializes to the standard \v{C}ech description of flat $A$-gerbes when $\cG=\pt\sq A$, see \S\ref{ex: gerbes}. Similarly, when $\cG=G$ is an ordinary group we obtain the description of flat $G$-bundles in terms of a cover and gluing data, see Example~\ref{ex:ordinaryGbund}. Furthermore, flat $\cG$-bundles in the sense of this paper give examples of principal $\cG$-bundles in the sense of \cite{SP11}, see Lemma~\ref{lem: gluing 2-bundles}. The explicit description of this bicategory $\Bun_\cG(X)$ is what allows for the proofs of Theorems~\ref{thm2} and~\ref{thm1} via concrete computations. 

Ultimately, one would hope for a generalization of our main results for arbitrary categorical central extensions in the bicategory of smooth 2-groups, and in particular for the string group~\eqref{eq:stringgroupdef}.

\begin{conjecture}\label{conjecture}
For a smooth categorical central extension 
$$
1\to\pt\sq A\to \cG\to G\to 1
$$
of a Lie group~$G$ classified by $[\alpha]\in \rmH^3_{\rm SM}(G;A)$, a $\cG$-bundle on a smooth stack $X$ is equivalent to a $G$-bundle $P\to X$ with a trivialization of the 2-gerbe $\lambda_{P,\alpha}\to X$ determined by $P$ and the Segal--Mitchison cocycle $\alpha$. 

In particular, string structures on a vector bundle $V\to X$ are equivalent to trivializing sections $\sigma\colon X\to \CS_V$ for $\CS_V$ the Chern-Simons 2-gerbe of $V$.
\end{conjecture}

The current paper is the first in a program of the authors to prove Conjecture~\ref{conjecture}. We anticipate this will reveal the correct notion of \emph{string connection} on a $\String_n$-principal bundle, an important step in the Stolz--Teichner program to construct a geometric model for the string orientation of TMF~\cite{ST04}. 

The flatness assumptions of the current paper simplify some of the stacky geometry expected in Conjecture~\ref{conjecture}, while also retaining the anticipated categorical complexities of the complete picture. 
In particular, studying $\cG$-bundles for a weak 2-group $\cG$ reveals that the nontrivial associator~$\alpha$ changes the flavor of the geometry rather dramatically when compared to bundles for strict 2-groups; see Remark~\ref{rmk:strictweak}. Similar behavior must persist for principal bundles for arbitrary smooth 2-groups.

\subsection{Outline and conventions}

We begin in \S\ref{sec: review of 2-groups} with a review of 2-groups. In particular, we define a 2-group $\cG(G,A,\alpha)$ associated to a discrete group $G$, a discrete abelian group $A$ with an action of $G$, and a 3-cocycle $\alpha$. In the remainder of the paper, we focus on the case where the $G$-action on $A$ is trivial. 
In~\S\ref{sec: bundles on groupoids}, we introduce a bicategory of principal $\cG$-bundles on a Lie groupoid $X = \lbrace X_1 \rightrightarrows X_0\rbrace$. Initially we formalize this notion via the bicategory $\Fun(X,\pt\sq \cG)$, but in the case of a 2-group $\cG(G,A,\alpha)$, we unwind the data to give an equivalent bicategory, which we will denote $\Bun_\cG^\pre(X)$. These bicategories have as objects principal $\cG$-bundles that trivialize over the objects $X_0$ of $X$. Hence, the assignment $X \mapsto \Bun_\cG^\pre(X)$ defines a prestack over the bicategory $\LieGrpd$ of Lie groupoids and smooth functors; the superscript ``pre" emphasizes the fact that this prestack is (almost never) a stack. 

In \S\ref{sec: localization}, we extend the assignment $X \mapsto \Bun_\cG^\pre(X)$ to a stack. To this end, we first define a bicategory $\Bun_\cG^\pre$ which is fibered in 2-groupoids over $\LieGrpd$ and localize with respect to the class of 1-morphisms covering essential equivalences in $\LieGrpd$. This results in a bicategory $\Bun_\cG$ over the localized bicategory $\LieGrpd[{\sf W}^{-1}]$ of differential stacks, where ${\sf W}$ is the class of essential equivalences of Lie groupoids. By construction, the homotopy fiber of $\Bun_\cG$ at $X$ is a bicategory $\Bun_\cG(X)$ of principal $\cG$-bundles on $X$. As an application of this construction, we study the natural functor $\pi: \Bun_\cG(X) \to \Bun_G(X)$; this enables us to characterize a principal $\cG$-bundle $\cP$ in terms of its underlying $G$-bundle $P= \pi(\cP)$ and the additional data of a trivialization of a 2-gerbe $\lambda_{P, \alpha}$ (Theorem \ref{thm1}/\ref{thm:2gerbetriv}).

In~\S\ref{sec:flatstringserction} we turn our attention to structures on flat vector bundles over a geometric stack. We show that orientations, flat spin structures, and flat string structures can equivalently be phrased in terms of lifts as in~\eqref{eq:stringlift1}, or as trivializations of line bundles, gerbes, or 2-gerbes in the respective cases. This proves Theorem~\ref{thm2}. 

Throughout, a functor between bicategories will mean a pseudofunctor, i.e., a particular type of lax 2-functor. 

A group $G$ can be viewed as a groupoid with a single object $\pt$ and morphisms~$G$, which we denote by $\pt\sq G$, where multiplication in $G$ corresponds to composition in the category. 
Given a group $G$ acting on a set $S$ on the right, the \emph{action groupoid} has objects the elements of $S$, morphisms $S\times G$ (where the source map is the projection and the target map is the action map), and composition is determined by multiplication in $G$. We denote this groupoid by $S\sq G$.

\subsection{Acknowledgements}
This paper began in conversations at the 2019 Mathematics Research Community, ``Geometric Representation Theory and Equivariant Elliptic Cohomology." We thank the AMS and NSF for supporting this program, as well as Eric Berry and Joseph Rennie who contributed to this project in its early stages. We also thank Matt Ando, Nora Ganter, Connor Grady, Eugene Lerman, Christopher Schommer-Pries, Nat Stapleton, and Stephan Stolz for stimulating conversations on these ideas.
We also thank the anonymous referees for their helpful comments.
This material is based upon work supported by the National Science Foundation under Grant Number DMS 1641020. The first author was partially supported by the National Science Foundation under Grant Number DMS 2205835. The second author was partially supported by the Natural Sciences and Engineering Research Council of Canada under Discovery Grant RGPIN-2022-04104 and the Fonds de Recherche du Qu\'ebec (Nature et Technologie) Research Support for New Academics under Grant 327706. The third author was partially supported by the National Science Foundation under Grant Number DMS 2316646.

\section{Review of 2-groups} \label{sec: review of 2-groups}

In this section, we introduce some foundational terminology and examples of 2-groups. The example of crucial interest in later sections of the paper is the 2-group $\cG(G, A, \alpha)$ associated to a discrete group $G$, an abelian group $A$ (with trivial $G$-action), and a 3-cocycle $\alpha$ as in Remark \ref{rmk:normalized2}; however, in this section we consider the general case of a nontrivial $G$-action on $A$. 
None of the definitions or results in this section are original; standard references include~\cite{BL04,SP11}

\subsection{Basic definitions}

\begin{defn} \label{defn:2group}
A (discrete) \emph{2-group} is a monoidal groupoid $(\cG,\otimes,\one)$ where every object is (weakly) $\otimes$-invertible, meaning for every object $x$ there exists an object $x^{-1}$ and isomorphisms $x\otimes x^{-1}\simeq \one\simeq x^{-1}\otimes x.$ A \emph{1-homomorphism} between 2-groups is a lax monoidal functor.
A \emph{2-homomorphism} between 1-homomorphisms is a lax monoidal transformation.
\end{defn}

\begin{prop}
The collection of 2-groups, 1-homomorphisms, and 2-homo-morphisms has the structure of a bicategory. 
\end{prop}

We often use the abbreviated notation $\cG$ for a 2-group, and adopt the following terminology that generalizes similar language for ordinary groups. Let $\pt$ denote the terminal category. A 2-group has a \emph{unit} functor $\one \colon \pt\to \cG$ and a \emph{multiplication} functor $\mu\colon \cG\times \cG\to \cG$. It also has an \emph{associator} natural transformation $\alpha$ and \emph{unitor} natural transformations~$\epsilon_l,\epsilon_r$ fitting into the diagrams
\beq
\begin{tikzpicture}[baseline=(basepoint)];
\node (A) at (0,0) {$\cG\times\cG\times\cG$};
\node (B) at (3,0) {$\cG\times\cG$};
\node (D) at (3,-1.5) {$\cG$};
\node (C) at (0,-1.5) {$\cG\times \cG$,};
\node (H) at (1.5,-.75) {$\alpha\ \twocommute$};
\draw[->] (A) to (B);
\draw[->] (B) to node [right] {$\mu\times \id_\cG $} (D);
\draw[->] (A) to node [left] {$\id_\cG\times \mu$} (C);
\draw[->] (C) to node [below] {$\mu$} (D);
\path (0,-.75) coordinate (basepoint);
\end{tikzpicture}\nonumber\qquad\begin{tikzpicture}[baseline=(basepoint)];
\node (A) at (0,0) {$\cG$};
\node (B) at (2,0) {$\cG\times\cG$};
\node (C) at (-2.5,0) {$\cG\times \cG$};
\node (D) at (0,-2) {$\cG$};
\node (H) at (.7,-.5) {$\epsilon_l\ \Rightarrow$};
\node (H) at (-.7,-.5) {$\Leftarrow\ \epsilon_r$};
\draw[->] (A) to node [above] {$\one\times \id_\cG$} (B);
\draw[->] (B) to node [right] {$\mu$} (D);
\draw[->] (A) to node [above] {$\id_\cG\times \one$} (C);
\draw[->] (C) to node [below] {$\mu$} (D);
\draw[->] (A) to node [left] {$\id_\cG$} (D);
\path (0,-.75) coordinate (basepoint);
\end{tikzpicture}
\eeq
and satisfying standard properties (e.g., the pentagon axiom).  A 1-homomorphism is the data of a functor $\varphi \colon \cG\to \cG'$ and natural transformation $\gamma$, 
\beq
\begin{tikzpicture}[baseline=(basepoint)];
\node (A) at (0,0) {$\cG\times\cG$};
\node (B) at (4,0) {$\cG'\times\cG'$};
\node (D) at (4,-1.5) {$\cG'$};
\node (C) at (0,-1.5) {$\cG$,};
\node (H) at (2,-.75) {$\gamma \ \twocommute$};
\draw[->] (A) to node [above] {$\varphi\times\varphi$} (B);
\draw[->] (B) to node [right] {$\mu$} (D);
\draw[->] (A) to node [left] {$\mu'$} (C);
\draw[->] (C) to node [below] {$\varphi$} (D);
\path (0,-.75) coordinate (basepoint);
\end{tikzpicture}\nonumber
\eeq
making $\varphi$ into a monoidal functor via the standard properties, including a commuting cube intertwining  the associators for $\cG$ and $\cG'$ using $\gamma$. A 2-homomorphism is a natural transformation $(\varphi,\gamma) \Rightarrow (\varphi',\gamma')$ satisfying compatibility properties. 

There are a few variations of Definition~\ref{defn:2group} in the literature, which we record below. 

\begin{defn}[{\cite[Definition 29]{BL04}}]
    A 2-group with trivial associator and unitors is called a \emph{strict 2-group}. 
\end{defn}

\begin{defn}[{\cite[Definition~7]{BL04}}]\label{defn:coherent}
A \emph{coherent 2-group} $\cG$ is a 2-group equipped with the additional data of an object $x^{-1}$ for each $x\in \cG$ together with specified isomorphisms $x\otimes x^{-1}\simeq \one \simeq x^{-1}\otimes x$ satisfying the axioms of an adjoint equivalence. 
\end{defn}

\subsection{Examples of 2-groups}

\begin{ex} An ordinary group $G$ determines a 2-group whose underlying groupoid has only identity morphisms with monoidal structure given by the group multiplication on~$G$. A homomorphism of groups is equivalent data to a monoidal functor between the corresponding monoidal categories. Hence, the category of groups and homomorphisms admits a faithful embedding into the bicategory of 2-groups. 
\end{ex}

\begin{ex}\label{ex:ptmodA}
For an abelian group $A$, the groupoid $\pt\sq A$ has the structure of a 2-group, where the monoidal structure on morphisms is determined by multiplication in $A$.
\end{ex}

\begin{ex}
For any monoidal category $({\sf C},\otimes, \one)$, the \emph{Picard 2-group} of ${\sf C}$, denoted $\Pic({\sf C})$, is the monoidal subgroupoid of ${\sf C}$ with objects the $\otimes$-invertible objects of ${\sf C}$ and morphisms the isomorphisms. In the case that $({\sf C},\otimes, \one)=({\sf Vect}_\K,\otimes,\K)$ is the (symmetric) monoidal category of vector spaces over $\K$ with tensor product, then $\Pic({\sf Vect}_\K)$ is the category of 1-dimensional $\K$-vector spaces. 
\end{ex}

\begin{ex} 
Given a category ${\sf C}$, the functor category ${\sf Fun}({\sf C},{\sf C})$ is a (strict) monoidal category with monoidal structure coming from composition of functors. The subcategory 
$$
\Aut({\sf C}):=\Pic({\sf Fun}({\sf C},{\sf C}))
$$
is the 2-group of automorphisms of ${\sf C}$. Composition of functors is strictly associative and unital; hence the associator and unitors for this 2-group are trivial. 
\end{ex}

An early example of 2-groups came from the concept of crossed modules, first formulated by Whitehead in \cite{Whitehead1941, Whitehead1949}.

\begin{ex}
A \emph{crossed module} is the data of homomorphisms $t\colon H\to G$ and $a\colon G\to \Aut(H)$ between (ordinary) groups such that the diagrams commute
\beq
\begin{tikzpicture}[baseline=(basepoint)];
\node (A) at (0,0) {$H\times H$};
\node (B) at (3,0) {$G\times H$};
\node (C) at (1.5,-1.5) {$H$};
\draw[->] (A) to node [above] {$t\times \id_H$} (B);
\draw[->] (A) to node [left] {${\rm Ad}$} (C);
\draw[->] (B) to node [right] {$a$} (C);
\path (0,-.75) coordinate (basepoint);
\end{tikzpicture}\qquad \begin{tikzpicture}[baseline=(basepoint)];
\node (A) at (0,0) {$G\times H$};
\node (B) at (3,0) {$H$};
\node (D) at (3,-1.5) {$G$,};
\node (C) at (0,-1.5) {$G\times G$};
\draw[->] (A) to node [above] {$a$} (B);
\draw[->] (B) to node [right] {$t$} (D);
\draw[->] (A) to node [left] {$\id_G\times t$} (C);
\draw[->] (C) to node [below] {${\rm Ad}$} (D);
\path (0,-.75) coordinate (basepoint);
\end{tikzpicture}
\eeq
where ${\rm Ad}$ is the adjoint action and (in a mild abuse of notation) the arrows labeled $a$ are determined by the homomorphism~$a$. This determines a 2-group with underlying groupoid the action groupoid $G\sq H$ for the (left) $H$-action on $G$ via $t$. The monoidal structure on this groupoid comes from the multiplication in $G$ and the semidirect product~$H\rtimes G$ determined by~$a$. The associator for this monoidal structure is trivial.
\end{ex}

The bicategories of strict 2-groups and crossed modules are equivalent; see~\cite{BL04}.

\begin{ex}
    Given any 2-group $\cG$ and (1-)category $\sf C$, the category of functors from $\sf C$ to $\cG$ is naturally a 2-group. 
\end{ex}

\subsection{Classification of 2-groups as categorical extensions}

For a 2-group~$\cG$, let $G:=\pi_0(\cG)$ denote the isomorphism classes of objects in~$\cG$, $A:=\Aut(\one)$ denote the automorphism group of the unit object~$\one \colon \pt\to  \cG$, and $\iota\colon \pt\sq A\to \cG$ the inclusion of the automorphisms of the unit object. We refer to $G$ as the ordinary group underlying $\cG$. Note that we have the sequence of 1-homomorphisms
\beq\label{eq:catextension}
1\to \pt\sq A\xrightarrow{\iota} \cG\xrightarrow{\pi_0} G\to 1,
\eeq
where $G$ is an ordinary group regarded as a 2-group with group structure inherited from the monoidal structure on~$\cG$. The set $A$ has a pair of multiplications given by composition or multiplication in $\cG$; by the  Eckmann--Hilton argument these agree and are commutative. Furthermore, there is a $G$-action on $A$ from conjugating an automorphism $a\in \Aut(\one)$ by $\id_g$ for an object $g$ of $\cG$,
$$
g\mapsto \id_{g^{-1}}\otimes (a \otimes \id_g).
$$

Finally, the associator of $\cG$ determines a 3-cocycle $\alpha\colon G\times G\times G\to A$, where the pentagon axiom implies that $\alpha$ satisfies the cocycle condition. In analogy to the structures in ordinary extensions of a group $G$ by an abelian group $A$, we call~\eqref{eq:catextension} a \emph{categorical extension} of $G$ by $\pt\sq A$. When the $G$-action on $A$ is trivial,~\eqref{eq:catextension} is a \emph{categorical central extension}.

\begin{prop}\label{prop:2grpconstr}
A 2-group $\cG$ is determined up to 1-isomorphism by
\begin{enumerate}
\item a group $G=\pi_0(\cG)$, the isomorphim classes of objects in $\cG$; 
    \item an abelian group $A=\Aut(1_G)$;
    \item a $G$-action on $A$, $G=\pi_0(\cG) \to \Aut(A)$;
    \item a class $[\alpha]\in \rmH^3(G;A)$.
\end{enumerate}
\end{prop} 

\bp
This follows from \cite[Theorem 43]{BL04} (based on \cite{sinh}); for central extensions, it is an immediate consequence of \cite[Theorem 99]{SP11}. For the reader's convenience we review the construction on objects. 

Above we already described how to extract data (1)--(4) from a 2-group. Next we describe how to construct a 2-group from the data (1)--(4). Start with the groupoid $\{G\ltimes A\rightrightarrows G\}$ where the source and target are the projections. Define a monoidal structure via multiplication in $G$ and $G\ltimes A$. The 3-cocyle $\alpha$ gives the associator. The left and right unitors are determined by $\alpha$ and the requirement that the triangle identities hold. This yields a 2-group~$\cG(G, A, \alpha)$.
\ep

\begin{rmk}\label{rmk:normalized}
For trivial $G$-action on $A$, a 3-cocycle $\alpha\in Z^3(G;A)$ is \emph{normalized} if
$$
\alpha(1_G,g_1,g_2)=\alpha(g_1,1_G,g_2)=\alpha(g_1,g_2,1_G)=1_A,
$$
where $1_G\in G$ and $1_A\in A$ are the identity elements. Every 3-cocycle is cohomologous to a normalized one.
Hereafter we shall assume that all 3-cocycles are normalized. 
\end{rmk}

\begin{rmk}\label{rmk:normalized2}
Let $\cG = \cG(G,A,\alpha)$ be a 2-group determined by data $(G, A, \alpha)$ as in Proposition~\ref{prop:2grpconstr}, where $A$ is a trivial $G$-module and $\alpha$ is a normalized 3-cocycle. In this case the left and right unitors are trivial. Furthermore, this 2-group can be endowed with a coherent structure (see Definition~\ref{defn:coherent}), where the monoidal inverse of an object $g \in G$ is given by its multiplicative inverse $g^{-1}$, and the unit and counit structure morphisms are
\beq\label{eq:unitcounitnormalized}
i_g = 1_A \colon 1_G \to g\otimes g^{-1}, \quad e_g=\alpha(g^{-1},g,g^{-1})\colon g^{-1}\otimes g\to 1_G
\eeq
for all objects $g\in G$, i.e., the unit $i$ is the identity and counit $e$ is determined by~$\alpha$. (Here we use that $g^{-1} \otimes g$ and $g \otimes g^{-1}$ are both equal to the object $1_G$, so $i_g$ and $e_g$ are both elements of $\Aut(1_G) = A$.) Coherent 2-groups of this form will be the main examples of interest in this paper.
\end{rmk}

\section{Flat principal 2-group bundles on a Lie groupoid}\label{sec: bundles on groupoids}

In this section, we define for each 2-group $\cG$ and Lie groupoid $X$ a bicategory $\Bun_\cG^\pre(X)$ of flat principal $\cG$-bundles on $X$ admitting trivializations over the objects~$X_0$ of the Lie groupoid $X$. 

\subsection{Background: Lie groupoids and essential equivalences} \label{subsec: background}

The material in this subsection is classical, e.g., see the references~\cite{Lerman, Mackenzie, Moerdijk2002}. Readers already familiar with Lie groupoids may want to skip to Section~\ref{subsec: ordinary G-bundles}.  

\begin{defn}[{e.g., {\cite[Defn 2.11]{Lerman}}}]
A  \emph{Lie groupoid} is a groupoid internal to the category of manifolds whose source and target maps are submersions. We use the notation $X= \{X_1 \rightrightarrows X_0\}$, where $X_0$ is the manifold of objects and $X_1$ is the manifold of morphisms. We use the notation $s,t\colon X_1\to X_0$ for the source and target maps, and $c\colon X_1\times_{X_0}X_1\to X_1$ for the composition. For $n \ge 2$, we denote by $X_n = X_1 \times_{X_0} \ldots \times_{X_0} X_1$, the manifold of sequences of $n$ composable morphisms in $X$. 

A \emph{(smooth) functor} between Lie groupoids $X=\{X_1\rightrightarrows X_0\}$ and $Y=\{Y_1 \rightrightarrows Y_0\}$ is a functor between their underlying categories that is a smooth map on objects $f_0\colon X_0\to Y_0$ and a smooth map on morphisms $f_1\colon X_1\to Y_1$. A \emph{(smooth) natural transformation} $\eta \colon f\Rightarrow g$ between functors is a natural transformation between underlying functors with the property that the map $\eta\colon X_0\to Y_1$ is smooth. 
\end{defn}

\begin{ex}\label{eq:quotientgroupoid}
Given a compact Lie group $G$ acting smoothly on a manifold $M$ on the right, the action groupoid $M\sq G$ has objects $M$, morphisms $M\times G$, source map the projection, target map the action map, and composition determined by multiplication in $G$. Regarding $M$ as a Lie groupoid with only identity morphisms, there is a smooth functor $M\to M\sq G$.
\end{ex}

Given a smooth functor $f\colon X\to Y$ between Lie groupoids, consider the diagrams of smooth manifolds 
\beq\label{eq:essentialequiv}
\begin{tikzpicture}[baseline=(basepoint)];
\node (A) at (0,0) {$Y_1\times_{Y_0}X_0$};
\node (B) at (3,0) {$X_0$};
\node (C) at (0,-1.5) {$Y_1$};
\node (D) at (3,-1.5) {$Y_0$};
\draw[->] (A) to node [above] {$p_2$} (B);
\draw[->] (B) to node [right] {$f_0$} (D);
\draw[->] (A) to node [left] {$p_1$} (C);
\draw[->] (C) to node [below] {$s$} (D);
\path (0,-.75) coordinate (basepoint);
\end{tikzpicture}\qquad 
\begin{tikzpicture}[baseline=(basepoint)];
\node (A) at (0,0) {$X_1$};
\node (B) at (3,0) {$Y_1$};
\node (C) at (0,-1.5) {$X_0\times X_0$};
\node (D) at (3,-1.5) {$Y_0\times Y_0.$};
\draw[->] (A) to node [above] {$f_1$} (B);
\draw[->] (B) to node [right] {$s\times t$} (D);
\draw[->] (A) to node [left] {$s\times t$} (C);
\draw[->] (C) to node [below] {$f_0\times f_0$} (D);
\path (0,-.75) coordinate (basepoint);
\end{tikzpicture}
\eeq

\begin{defn}[{e.g., \cite[Definition 3.5]{Lerman}}] \label{defn:localequiv}
A smooth functor $f\colon X\to Y$ is an \emph{essential equivalence} if $Y_1\times_{Y_0}X_0\to Y_0$ on the left in~\eqref{eq:essentialequiv} is a surjective submersion, and the diagram on the right in~\eqref{eq:essentialequiv} is a pullback. 
\end{defn}

\begin{rmk} 
The two conditions in Definition~\ref{defn:localequiv} are smooth versions of the conditions for a functor to be essentially surjective and fully faithful. 
\end{rmk}

\begin{ex}[\cite{Pronk}, Section 2.1]\label{eq:dumbessential}
Suppose that $f$ is a smooth functor $f\colon X\to Y$ for which there exists a smooth functor $g\colon Y\to X$ with $f\circ g\simeq \id$ and $g\circ f\simeq \id$. Then $f$ is an essential equivalence. 
\end{ex}

\begin{ex}[{\cite[Example~1.1.9]{Mackenzie}, \cite[Example 3.29]{Lerman}}] Let $P\to X$ be a principal $G$-bundle. Then there is an essential equivalence $P\sq G \to X=P/G$. It is invertible in the sense of Example~\ref{eq:dumbessential} if and only if $P \cong X \times G$ is the trivial principal $G$-bundle. 
\end{ex}

\begin{ex}[{\cite[Example~1.1.8]{Mackenzie}}] \label{ex: ess equiv to manifold}
Given a smooth manifold $M$ and a surjective submersion $Y\to M$, set $Y^{[k]}$ to be the $k$-fold fibered product of $Y$ over $M$, e.g., $Y^{[2]}=Y\times_M Y$. The \emph{\v{C}ech groupoid} associated with $Y\to M$ is 
\beq
\check{C}(Y):=\{Y^{[2]}\rightrightarrows Y\}
\eeq
with source and target maps given by projection to the factors, and composition given by the map $Y^{[3]}\to Y^{[2]}$ that projects out the middle factor in the fibered product. As a special case of this example, given an open cover $\{U_i\}$ of $M$, the map $\coprod_i U_i\to M$ is a surjective submersion, and we use the common notation for intersections, e.g., $\coprod U_{ij}=\coprod U_i\bigcap U_j=Y^{[2]}$ and $\coprod U_{ijk}=\coprod U_i\bigcap U_j\bigcap U_k=Y^{[3]}$. In this case, we use the notation $\check{C}(U_i)$ for the \v{C}ech groupoid associated to the cover. Returning to the general case of a surjective submersion $Y\to M$, there is an essential equivalence
\beq\label{eq:Cechequivalence}
\check{C}(Y)\to M
\eeq
from the \v Cech groupoid to the manifold $M$ viewed as a Lie groupoid $M$.

Conversely, suppose that $f \colon X \to M$ is an essential equivalence from a Lie groupoid $X$ to a manifold $M$. Since $Y_0 = Y_1 = M$ in the notation of the diagrams \eqref{eq:essentialequiv}, we deduce from the first diagram that the morphism $f_0 \colon X_0 \to M$ is a surjective submersion, and from the second that $X_1$ is homeomorphic to $X_0 \times_{M} X_0$. In other words, the only essential equivalences of a Lie groupoid into a manifold are isomorphic to those coming from \v Cech groupoids.  
\end{ex}

\begin{ex}[{\cite[Proposition~2.26]{Lerman}}]\label{ex:pullbackgrpod}
Generalizing the previous example, given a Lie groupoid $X$ and a surjective submersion $f\colon Y\to X_0$, define the \emph{pullback groupoid} $f^*X$ with objects $(f^*X)_0=Y$ and morphisms the pullback
\beq\nonumber
\begin{tikzpicture}[baseline=(basepoint)];
\node (A) at (0,0) {$(f^*X)_1$};
\node (B) at (3,0) {$X_1$};
\node (C) at (0,-1.5) {$Y\times Y$};
\node (D) at (3,-1.5) {$X_0\times X_0$};
\draw[->] (A) to  (B);
\draw[->] (B) to node [right] {$s\times t $} (D);
\draw[->] (A) to (C);
\draw[->] (C) to node [below] {$f\times f$} (D);
\path (0,-.75) coordinate (basepoint);
\end{tikzpicture}\qquad 
\eeq
i.e., for any pair $x,y\in Y$ we take the morphisms between $f(x)$ and $f(y)$ in $X_1$. There is an evident smooth functor $f^*Y\to X$; it is an essential equivalence. 
\end{ex}

\subsection{Motivation: Ordinary $G$-bundles in the language of Lie groupoids}\label{subsec: ordinary G-bundles}
\begin{lem}\label{lem:Gbundlefromgrpod}
    Given a Lie group $G$ and a surjective submersion $u\colon Y \to M$, a smooth functor
    \beq\label{eq:Gbundasfunctor1}
\check{C}(Y)\to \pt\sq G. 
\eeq
determines a principal $G$-bundle $P\to M$ with a trivialization of the pullback bundle $u^*P\to Y$. More precisely, a smooth functor~\eqref{eq:Gbundasfunctor1} determines a (strictly) commuting square of Lie groupoids
\beq\label{eq:totalspacefunctor}
    &&\begin{tikzpicture}[baseline=(basepoint)];
\node (A) at (0,0) {$(Y^{[2]}\times G \rightrightarrows  Y\times G)$};
\node (B) at (5,0) {$P$};
\node (C) at (0,-1.5) {$\check{C}(Y)$};
\node (D) at (5,-1.5) {$M$};
\draw[->] (A) to node [above] {$\sim$} (B);
\draw[->] (B) to (D);
\draw[->] (A) to (C);
\draw[->] (C) to node [below] {$\sim$} (D);
\path (0,-.75) coordinate (basepoint);
\end{tikzpicture}
\eeq
where $P\to M$ is a principal $G$-bundle, the vertical arrows are the natural projections, and the horizontal arrows are essential equivalences. 
Furthermore, a natural transformation of functors~\eqref{eq:Gbundasfunctor1} determines an isomorphism of principal bundles~$P\xrightarrow{\sim}P'$ over $M$. 
\end{lem}

\bp
Cocycle data for a $G$-bundle is a surjective submersion
$u\colon Y\to M$ and a transition function $g\colon Y^{[2]}\to G$ satisfying the cocycle condition
\beq\label{eq:strictcocycle}
g(y_1,y_2) \cdot  g(y_2, y_3) = g(y_1, y_3) & \text{ for all } (y_1, y_2, y_3) \in Y^{[3]};
\eeq 
e.g., see \cite{NikolausSchweigert}.
This is precisely the data and property of a smooth functor~\eqref{eq:Gbundasfunctor1} between Lie groupoids.
For cocycle data $g,g'$ defining a pair of $G$-bundles relative to the same open cover $Y$, an isomorphism between them is the data of a map $h\colon Y\to G$ satisfying
\beq\label{eq:strict1iso}
h(y_1) \cdot g(y_1, y_2) \cdot h(y_2)^{-1} = g'(y_1, y_2) \text{ for all } (y_1, y_2) \in Y^{[2]}.
\eeq
This is the same data and property as a smooth natural transformation of functors.
\ep

Let $\Bun_G(M)$ denote the category of $G$-bundles on $M$. An object of this category can be specified by a choice of surjective submersion $u \colon Y \to M$ and a smooth functor~\eqref{eq:Gbundasfunctor1}, which is equivalent to a zig-zag
\beq\label{eq:zigzag1}
M\xleftarrow{\sim} \check{C}(Y)\to \pt\sq G
\eeq
where the left arrow is an essential equivalence. 
We would like to understand the generalization of this category to the case where the base is allowed to be an arbitrary Lie groupoid $X=\{X_1\rightrightarrows X_0\}$. It is easy to generalize the presentation of objects in terms of covers and smooth functors \eqref{eq:zigzag1}; see Definition \ref{defn:prinbund0} below. However, it becomes more awkward to describe the morphisms between such objects, as these can involve pairs of distinct essential equivalences, allowing for e.g., refinements of covers.

We describe a resolution to this issue that readily generalizes to 2-group bundles. First consider the category $\Bun_G^\pre(X):=\Fun(X,\pt\sq G)$ of $G$-bundles on the Lie groupoid~$X$ that trivialize on objects. Since $G$-bundles pull back, letting $X$ vary gives a fibered category $\Bun_G^\pre$ with a forgetful functor to Lie groupoids. Localizing this category with respect to essential equivalences of Lie groupoids (i.e., stackifying) gives a category $\Bun_G$ whose objects are zig-zags~\eqref{eq:zigzag1}; this is the desired category of $G$-bundles over smooth stacks~\cite{Pronk}. We explain this carefully in the setting of principal bundles for a 2-group $\cG$ in \S\ref{sec: localization}, which reduces to this setting in the case $\cG = G$.

\subsection{Defining flat principal 2-group bundles over Lie groupoids}\label{subsec: definition of principal 2-group bundles}
The following definition generalizes ordinary $G$-bundles as described in Lemma~\ref{lem:Gbundlefromgrpod}. 
The smoothness condition below is the requirement that the structure maps (e.g., maps on objects and morphisms) are smooth relative to the Lie groupoid $X$ and the discrete topology on the 2-group~$\cG$, see Remark~\ref{rmk:pi0}. These structure maps are spelled out in detail in Lemmas~\ref{lem:prinbund2functor} and~\ref{lem:1-morphism no hats}.

\begin{defn}\label{defn:prinbund0}
For a Lie groupoid $X=\{X_1\rightrightarrows X_0\}$, a \emph{(flat) $\cG$-bundle on~$X$} is a zig-zag 
    \beq\label{eq:Gbundasfunctor0}
 X\xleftarrow{\sim} Y \xrightarrow{\cP} \pt\sq \cG
\eeq
for an essential equivalence $Y\xrightarrow{\sim} X$, and a smooth functor from $Y$ to the 1-object bicategory $\pt\sq \cG$ associated to the monoidal category $\cG$. 
 \end{defn}

The goal of the next two sections is to construct a bicategory whose objects are $\cG$-bundles~\eqref{eq:Gbundasfunctor0}. This bicategory is constructed by localizing (i.e., stackifying) the following intermediate definition of $\cG$-bundle; we construct this localization in~\S\ref{sec: localization}. 
\begin{defn}\label{defn:prinbund}
For a Lie groupoid $X=\{X_1\rightrightarrows X_0\}$, a \emph{(flat) $\cG$-bundle on~$X$ that trivializes on $X_0$} is a smooth functor
    \beq\label{eq:Gbundasfunctor}
\cP\colon X \to \pt\sq \cG
\eeq
 from $X$ to the 1-object bicategory $\pt\sq \cG$ associated to the monoidal category $\cG$. These are the objects of the bicategory of smooth functors, natural transformations, and modifications~\cite[Theorem 4.4.11]{JohnsonYau}
\beq
\Bun_\cG^\pre(X):=\Fun(X,\pt\sq \cG).
\eeq
\end{defn}
 
The superscript ``pre" refers to the fact that the assignment $X\mapsto \Bun_\cG^\pre(X)$ is a prestack that is (almost never) a stack. 

\begin{rmk}
For readers familiar with smooth stacks, note the following repackaging of Definition~\ref{defn:prinbund}. When considering the smooth stack $[X]$ underlying a Lie groupoid $X$, a $\cG$-bundle in the sense of Definition~\ref{defn:prinbund} is a $\cG$-bundle that trivializes when pulled back along the smooth functor $X_0\to X$, i.e., trivializes on the chosen atlas for $[X]$; compare~\cite[Definition 70]{SP11} or~\cite[Proposition 4.12]{BehrendXu}. 
\end{rmk}

\begin{rmk}\label{rmk:pi0}
Since the objects and morphisms of $\cG$ have the discrete topology, a smooth functor~\eqref{eq:Gbundasfunctor} factors through connected components on objects and morphisms of $X$, and hence is equivalent to an ordinary functor (with no smoothness condition), 
$$
\{\pi_0X_1 \rightrightarrows \pi_0X_0\}\to \pt\sq \cG
$$ 
where the source is the groupoid of connected components, regarded as a 2-groupoid, and the target is a 2-groupoid. 
\end{rmk}

\begin{lem}\label{lem:underlyingGbundle}
Let $G:=\pi_0\cG$ be the ordinary group underlying $\cG$.
There is a natural functor
\beq\label{eq:underlyingGbundle}
\Bun_\cG^\pre(X)\to \Bun_G^\pre(X)
\eeq
from flat $\cG$-bundles over $X$ that trivialize over $X_0$ to (flat) $G=\pi_0(\cG)$-bundles on~$X$ that trivialize on $X_0$. 
\end{lem}
\bp
The source and target of~\eqref{eq:underlyingGbundle} are defined as functor categories, and the induced functor between them is
\beq\label{eq:underlying}
u_*\colon \Fun(X,\pt\sq \cG)\to \Fun(X,\pt\sq G),\qquad u\colon  \pt\sq \cG\xrightarrow{\pi_0} \pt\sq G
\eeq
for $u$ the map of bicategories induced by the monoidal functor $\cG\to \pi_0\cG=G$. 
\ep

\subsection{$\mathcal{G}$-bundles as categorical 1-cocycles}\label{subsec: categorical 1-cocycles}

In~\cite[\S2.5.1]{Bartels}, $\cG$-bundles are defined as 1-cocycles valued in the monoidal category $\cG$. Definition~\ref{defn:prinbund} recovers this, as demonstrated by the following lemma. 

\begin{lem}\label{lem:prinbund2functor}
Given a Lie groupoid $X=\lbrace X_1\rightrightarrows X_0\rbrace$ and a 2-group $\cG$, a flat $\cG$-bundle over $X$ that trivializes on $X_0$ is the data of
\begin{enumerate}
    \item a smooth (and hence locally constant) functor $\widehat{\rho}\colon X_1\to \cG$, where $X_1$ is viewed as a Lie groupoid with only identity morphisms, and $\cG$ is viewed as a Lie groupoid with the discrete topology on objects and morphisms;
    \item smooth (and  hence locally constant) natural transformations $\widehat{\gamma}$ and $\widehat{\epsilon}$
    \beq\label{eq:weakcocycle}
    &&\begin{tikzpicture}[baseline=(basepoint)];
\node (A) at (0,0) {$X_2$};
\node (B) at (3,0) {$X_1$};
\node (D) at (3,-1.5) {$\cG$};
\node (C) at (0,-1.5) {$\cG\times \cG$,};
\node (H) at (1.5,-.75) {$\widehat{\gamma}\ \twocommuter$};
\draw[->] (A) to node [above] {$c$} (B);
\draw[->] (B) to node [right] {$\widehat{\rho}$} (D);
\draw[->] (A) to node [left] {$p_1^*\widehat{\rho} \times p_2^*\widehat{\rho}$} (C);
\draw[->] (C) to node [below] {$\mu$} (D);
\path (0,-.75) coordinate (basepoint);
\end{tikzpicture}\qquad\begin{tikzpicture}[baseline=(basepoint)];
\node (A) at (0,0) {$X_0$};
\node (B) at (3,0) {$X_1$};
\node (D) at (3,-1.5) {$\cG$,};
\node (C) at (0,-1.5) {$*$};
\node (H) at (1.5,-.75) {$\widehat{\epsilon}\ \twocommuter$};
\draw[->] (A) to node [above] {$u$} (B);
\draw[->] (B) to node [right] {$\widehat{\rho}$} (D);
\draw[->] (A) to (C);
\draw[->] (C) to node [below] {$e$} (D);
\path (0,-.75) coordinate (basepoint);
\end{tikzpicture}
\eeq
\end{enumerate}
These data must satisfy unity constraints: the composition of 2-morphisms in
\beq\nonumber
\resizebox{\columnwidth}{!}{
\begin{tikzpicture}[baseline=(basepoint)];
\node (AA) at (-3,0) {$X_1$};
\node (A) at (0,0) {$X_2$};
\node (B) at (3,0) {$X_1$};
\node (D) at (3,-1.5) {$\cG$,};
\node (C) at (0,-1.5) {$\cG\times \cG$};
\node (CC) at (-3,-1.5) {$\cG$};
\node (H) at (1.5,-.75) {$\widehat{\gamma}\ \twocommuter$};
\node (H) at (-1.5,-.75) {$\widehat{\epsilon}\times \id_{\widehat{\rho}} \twocommuter$};
\node (H) at (0,-2) {$\varepsilon_l \Uparrow$};
\draw[->] (CC) to node [below] {$e\times \id_\cG$} (C);
\draw[->] (AA) to node [above] {$l_\id$} (A);
\draw[->] (AA) to node [left] {$\widehat{\rho}$} (CC);
\draw[->,bend left=25] (AA) to node [above] {$\id$} (B);
\draw[->,bend right=25] (CC) to node [below] {$\id$} (D);
\draw[->] (A) to node [above] {$c$} (B);
\draw[->] (B) to node [right] {$\widehat{\rho}$} (D);
\draw[->] (A) to (C);
\draw[->] (C) to node [below] {$\mu$} (D);
\path (0,-.75) coordinate (basepoint);
\end{tikzpicture}\quad
\begin{tikzpicture}[baseline=(basepoint)];
\node (AA) at (-3,0) {$X_1$};
\node (A) at (0,0) {$X_2$};
\node (B) at (3,0) {$X_1$};
\node (D) at (3,-1.5) {$\cG$,};
\node (C) at (0,-1.5) {$\cG\times \cG$};
\node (CC) at (-3,-1.5) {$\cG$};
\node (H) at (1.5,-.75) {$\widehat{\gamma}\ \twocommuter$};
\node (H) at (-1.5,-.75) {$\id_{\widehat{\rho}}\times \widehat{\epsilon} \twocommuter$};
\node (H) at (0,-2) {$\varepsilon_r \Uparrow$};
\draw[->] (CC) to node [below] {$\id_\cG\times e$} (C);
\draw[->] (AA) to node [above] {$r_\id$} (A);
\draw[->] (AA) to node [left] {$\widehat{\rho}$} (CC);
\draw[->,bend left=25] (AA) to node [above] {$\id$} (B);
\draw[->,bend right=25] (CC) to node [below] {$\id$} (D);
\draw[->] (A) to node [above] {$c$} (B);
\draw[->] (B) to node [right] {$\widehat{\rho}$} (D);
\draw[->] (A) to (C);
\draw[->] (C) to node [below] {$\mu$} (D);
\path (0,-.75) coordinate (basepoint);
\end{tikzpicture}}\eeq
equals the identity 2-morphism on $\widehat{\rho}$. Above, $l_\id,r_\id\colon X_1\to X_2$ are $l_\id(f)=\id \circ f$ and $r_\id(f)= f\circ \id$, respectively. The data (1) and (2) are also required to satisfy an associativity constraint in the form of a 2-commuting cube,
\beq\label{diag:prettycube}
\begin{tikzcd}[row sep=scriptsize, column sep=scriptsize]
& X_3  \arrow[dl] \arrow[rr] \arrow[dd] & & X_2 \arrow[dl] \arrow[dd] \\ 
X_2 \arrow[rr, crossing over] \arrow[dd] & & X_1 \\ 
& \cG\times\cG\times\cG \arrow[dl] \arrow[rr] & & \cG\times \cG \arrow[dl] \\ 
\cG\times \cG \arrow[rr] & & \cG \arrow[from=uu, crossing over]\\ 
\end{tikzcd}\label{eq:highercompatcond}
\eeq
where the top face commutes on the nose (by associativity of composition in $X$), the bottom face has 2-commuting data given by the associator $\alpha$ for $\cG$, and the sides carry pullbacks of 2-commuting data $\gamma$.
\end{lem}
\bp
This follows from unpacking the definition of a functor, see~\cite[\S4.1]{JohnsonYau}. 
Smoothness of $\widehat{\rho}$, $\widehat{\gamma}$ and $\widehat{\epsilon}$ is equivalent to local constancy, as the objects and morphisms of $\cG$ have the discrete topology. 
\ep

The 1- and 2-morphisms between $\cG$-bundles from~\cite[\S2.5.1]{Bartels} are stated in terms of coboundaries for $\cG$-valued 1-cocycles. This generalizes the usual cocycle description~\eqref{eq:strict1iso} of an isomorphism between ordinary $G$-bundles as follows. 

\begin{lem}\label{lem:1-morphism no hats}
A 1-morphism between flat $\cG$-bundles $(\widehat{\rho_1},\widehat{\gamma_1},\widehat{\epsilon_1})$ and $(\widehat{\rho_2},\widehat{\gamma_2},\widehat{\epsilon_2})$ over $X$ that trivialize on $X_0$ is
\begin{enumerate}
    \item a smooth (and hence locally constant) functor $\widehat{h}\colon X_0\to \cG$ 
    \item a smooth (and hence locally constant) natural transformation $\widehat{\eta}$, 
    \beq
    \begin{tikzpicture}[baseline=(basepoint)];
\node (A) at (0,0) {$X_1$};
\node (B) at (4,0) {$\cG\times \cG$};
\node (D) at (4,-1.5) {$\cG$};
\node (C) at (0,-1.5) {$\cG\times \cG$};
\node (H) at (2,-.75) {$\widehat{\eta}\ \twocommute$};
\draw[->] (A) to node [above] {$t^*\widehat{h}\times \widehat{\rho_1}$} (B);
\draw[->] (B) to node [right] {$\mu$} (D);
\draw[->] (A) to node [left] {$\widehat{\rho_2}\times s^*\widehat{h} $} (C);
\draw[->] (C) to node [below] {$\mu$} (D);
\path (0,-.75) coordinate (basepoint);
\end{tikzpicture}\label{eq:1iso}
\eeq
\end{enumerate}
These data are required to satisfy a compatibility condition; see~\cite[\S4.2]{JohnsonYau}. 

A 2-morphism is a smooth natural isomorphism between functors $\widehat{\omega}\colon\widehat{h_1}\Rightarrow \widehat{h_2}$ satisfying the condition that the square commutes:
\beq
\begin{tikzcd}
    \widehat{\rho_2}\cdot s^* \widehat{h_1} \arrow[r, "\widehat{\eta_1}"]\arrow[d, "s^*\widehat{\omega}"'] & t^*\widehat{h_1}\cdot \widehat{\rho_1} \arrow[d, "t^*\widehat{\omega}"]\\
    \widehat{\rho_2}\cdot s^*\widehat{h_2} \arrow[r, "\widehat{\eta_2}"'] & t^*\widehat{h_2}\cdot \widehat{\rho_1}
\end{tikzcd}
\eeq 
\end{lem}

\bp
The description of 1-morphisms and 2-morphisms follows from unpacking the definitions of natural transformations and modifications, see~\cite[\S4.2]{JohnsonYau} and \cite[\S4.4]{JohnsonYau} respectively. 
\ep

\subsection{2-group bundles as \v{C}ech cocycles}\label{subsec: remove hats}
In the special case where the 2-group $\cG$ is a central extension of the form $\cG=\cG(G, A, \alpha)$ as in Remark~\ref{rmk:normalized2}, then one can unpack the data of a $\cG$-bundle in terms of group cohomology data.

\begin{prop} \label{prop: remove hats} Let $\cG=\cG(G,A,\alpha)$ be as in Remark~\ref{rmk:normalized2}. A flat $\cG$-bundle over a Lie groupoid $X=\lbrace X_1\rightrightarrows X_0\rbrace$ that trivializes on $X_0$ is given by the same data as a pair $(\rho,\gamma)$ where
\begin{enumerate}
    \item $\rho\colon X_1\to G$ is a locally constant map satisfying the (ordinary) cocycle condition
\beq \label{eq: condition for rho}
c^*\rho = p_1^*\rho\cdot p_2^*\rho \colon X_2\to G
\eeq
for $p_1,p_2,c\colon X_2\to X_1$ the projections and the composition;
    \item $\gamma\colon X_2\to A$ is a locally constant map satisfying the conditions
 \beq\label{eq:highercompatcond2}
&&\rho^*\alpha=d \gamma\colon X_3\to A,\\
&&\gamma(\phi, \id_{s(\phi)})=\gamma(\id_{t(\psi)}, \psi), \quad \forall (\phi,\psi)\in X_2 \label{eq:epsilonconditions}
 \eeq
    where $d$ is the \v{C}ech differential on $A$-valued cochains, and \[ s(f)=x=t(g)\in X_0. \]  
\end{enumerate}
\end{prop}

\bp
By construction,  the objects of $\cG$ are the elements of the discrete group~$G$. Hence a locally constant functor~$\widehat{\rho}\colon X_1\to \cG$ is the same as a locally constant map $X_1\to G$. Since $\cG$ is skeletal (meaning objects in $\cG$ are isomorphic if and only if they are equal), the existence of a natural transformation~\eqref{eq:weakcocycle} requires $\rho$ to satisfy the cocycle condition on $X_2$. Since the morphisms of $\cG$ are $G\times A$, natural transformations~\eqref{eq:weakcocycle} are maps $\gamma\colon X_2\to A$ and $\epsilon \colon X_0\to A$. The condition~\eqref{eq:highercompatcond} is precisely~\eqref{eq:highercompatcond2}.
Since $\cG$ is defined using a normalized 3-cocycle, the conditions on $\epsilon \colon X_0\to A$ are
\beq\label{eq:conditionfromepsilon}
&&(s^*\epsilon)\cdot (l_\id^*\gamma)=1\colon X_1\to A,\quad (t^*\epsilon)\cdot (r_\id^*\gamma)=1\colon X_1\to A
\eeq
where $1\colon X_1\to A$ denotes the constant function to the unit of $A$. This implies that $\epsilon$ is determined by $\gamma$. Furthermore,~\eqref{eq:conditionfromepsilon} imposes the condition on $\gamma$, 
$$\gamma(\phi, \id_x)^{-1} =\epsilon(x)= \gamma(\id_x, \psi)^{-1}, \qquad (\phi,\psi)\in X_2
$$
where $s(\phi) = x=t(\psi)$. This verifies~\eqref{eq:epsilonconditions}.
\ep

\begin{prop}\label{prop: remove hats for 1-morphisms}
 Let $\cG = \cG(G, A, \alpha)$ be as in Remark~\ref{rmk:normalized2} and let $X = \lbrace X_1 \rightrightarrows X_0\rbrace$ be a Lie groupoid. Suppose that $(\widehat{h}, \widehat{\eta}) \colon (\widehat{\rho}_1, \widehat{\gamma}_1) \to (\widehat{\rho}_2, \widehat{\gamma}_2)$ is a 1-morphism between flat principal $\cG$-bundles that trivialize on $X_0$ which are specified by data $(\rho_i, \gamma_i), i = 1,2$ as in Proposition \ref{prop: remove hats}. Then $(\widehat{h}, \widehat{\eta})$ is given by the same data as a pair $(h, \eta)$, where 
 \begin{enumerate}
     \item $h \colon X_0 \to G$ is a locally constant map such that 
     \beq\label{eq: h condition}
         \rho_2 \cdot s^* h =  t^* h \cdot \rho_1 \colon X_1 \to G;
     \eeq 
     \item $\eta \colon X_1 \to A$ is a locally constant map such that the following diagram commutes, for any pair of composable morphisms $(\phi, \psi)\in X_2$:
     \begin{center}
     \beq
    \begin{tikzcd} \label{eq:eta condition}
        \rho_2(\phi\circ \psi)s^*h(\psi)\arrow[r, "\eta(\phi\circ \psi)"]\arrow{d}{\gamma_2(\phi,\psi)\times id} & t^*h(\phi)\rho_1(\phi\circ \psi)\arrow{d}{id\times\gamma_1(\phi,\psi)}\\
        (\rho_2(\phi)\rho_2(\psi))s^*h(\psi)\arrow[d, "\alpha"] &t^*h(\phi) (\rho_1(\phi)\rho_1(\psi)) \\
        \rho_2(\phi)(\rho_2(\psi)s^*h(\psi)) \arrow[d, "id\times \eta(\psi)"]& (t^*h(\phi)\rho_1(\phi))\rho_1(\psi)\arrow[u, "\alpha"]\\
        \rho_2(\phi)(t^*h(\psi)\rho_1(\psi))  & (\rho_2(\phi)s^*h(\phi))\rho_1(\psi). \arrow[l, "\alpha"]\arrow[u, "\eta(\phi)\times id"]
    \end{tikzcd}
    \eeq 
\end{center}
 \end{enumerate}
\end{prop}

\begin{proof}
Because the objects of $\cG$ are given by the discrete group $G$, a locally constant functor $\widehat{h} \colon X_0 \to \cG$ is the same as a locally constant map $h \colon X_0 \to G$. 

Now consider the data of the natural transformation $\widehat{\eta} \colon \mu \circ (\widehat{\rho_2} \times s^*h) \Rightarrow \mu \circ (t^*h \times  \widehat{\rho_1})$ of functors $X_1 \to \cG$. It consists of a smooth map from the objects of $X_1$ to the morphisms of $\cG$, that is from $X_1 \to G \times A$, providing a family of isomorphisms
\begin{align*}
    \widehat{\eta}_\phi \colon \widehat{\rho_2}(\phi) \cdot h(s(\phi)) \to h(t(\phi)) \cdot \widehat{\rho_1}(\phi)) \in \cG, \; \phi \in X_1. 
\end{align*}
Because $\cG$ is skeletal with automorphisms given by $A$, we must have $\rho_2(\phi) \cdot h(s(\phi)) = h(t(\phi)) \rho_1(\phi) \in G$, and the data of $\widehat{\eta}_\phi$ is an element $\eta(\phi) \in A$.  

The compatibility condition specifies that $\eta$ behaves well with respect to composition of morphisms $(\phi, \psi) \in X_2$. 
\end{proof}

\begin{rmk}
    Considering $\eta(\phi)$ as the data of an isomorphism $\rho'(\phi) \cdot h(s(\phi)) \to  h(t(\phi)) \rho(\phi)$ in $\cG$ for each $\phi \in X_1$ motivates the formulas for composition of two 1-morphisms. Indeed, given 
    \begin{align*}
        (\rho_1, \gamma_1) \xrightarrow{(h_1, \eta_1)} (\rho_2, \gamma_2) \xrightarrow{(h_2, \eta_2)} (\rho_3, \gamma_3),
    \end{align*}
    we have for every $\phi \colon x \to x'$ in $X$ isomorphisms $\eta_i(\phi) \colon \rho_{i+1}(\phi) \cdot h_i(x) \to  h_i(x') \cdot \rho_i(\phi)$ for $i = 1, 2$. This yields the following chain of isomorphisms:
    \begin{align*}
        \rho_3(\phi)\cdot(h_2(x) \cdot h_1(x)) & \xrightarrow{\alpha^{-1}(\rho_3(\phi), h_2(x), h_1(x))} (\rho_3(\phi)\cdot h_2(x)) \cdot h_1(x) \xrightarrow{\eta_2(\phi)} \\
        (h_2(x') \cdot \rho_2(\phi)) \cdot h_1(x) 
         & \xrightarrow{\alpha(h_2(x'), \rho_x(\phi), h_1(x))} h_2(x') \cdot (\rho_2(\phi) \cdot h_1(x)) \xrightarrow{\eta_1(\phi)} \\
         h_2(x') \cdot (h_1(x') \cdot \rho_1(\phi)) 
         & \xrightarrow{\alpha^{-1}(h_2(x'), h_1(x'), \rho_1(\phi))}  (h_2(x') \cdot h_1(x')) \cdot \rho_1(\phi).
    \end{align*}
In other words, we have a 1-morphism $(h_{21}, \eta_{21}) \colon (\rho_1, \gamma_1) \to (\rho_3, \gamma_3)$ given by 
    \begin{align*}
        h_{21} \colon X_0  & \to G \\
        x & \mapsto h_2(x)h_1(x); \\
        \eta_{21} \colon X_1 & \to A \\
        \phi & \mapsto \eta_2(\phi)\eta_1(\phi) \frac{\alpha(h_2(x'), \rho_2(\phi), h_1(x))}{\alpha(\rho_3(\phi), h_2(x), h_1(x))\alpha(h_2(x'), h_1(x'), \rho_1(\phi))}.
    \end{align*}
    This is the composition $(h_2, \eta_2) \circ (h_1, \eta_1)$. 
\end{rmk}

These satisfy the compatibility conditions of Proposition~\ref{prop: remove hats for 1-morphisms} and hence do indeed give a 1-morphism as desired; this uses multiple iterations of the cocycle condition for $\alpha$, the cocycle condition for $\rho_i$ in equation~\ref{eq: condition for rho}, $i=1,2,3$, the conditions for $h_j$ in equation~\ref{eq: h condition}, and the conditions for $\eta_j$ in equation~\ref{eq:eta condition}, $j=1,2$.

\begin{prop}\label{prop: remove hats for 2-morphisms}
    Let $\cG = \cG(G, A, \alpha)$ be as in Remark~\ref{rmk:normalized2} and let $X = \lbrace X_1 \rightrightarrows X_0\rbrace$ be a Lie groupoid. 
    Suppose that we have two flat principal $\cG$-bundles $(\widehat{\rho}_i, \widehat{\gamma}_i), i = 1, 2$, specified by data $(\rho_i, \gamma_i)$ as in Proposition~\ref{prop: remove hats}, and two 1-morphisms $(\widehat{h}_i, \widehat{\eta}_i) \colon (\widehat{\rho}_1, \widehat{\gamma}_1) \to (\widehat{\rho}_2, \widehat{\gamma}_2)$ of principal $\cG$-bundles, specified by data $(h_i, \eta_i)$ as in Proposition~\ref{prop: remove hats for 1-morphisms}.
    Then  a 2-morphism $\widehat{\omega} \colon (\widehat{h}_1, \widehat{\eta}_1) \to (\widehat{h}_2, \widehat{\eta}_2)$  can exist only if $h_1 = h_2$. 
    Furthermore, the data of $\widehat{\omega}$ is just a morphism 
    \begin{align*}
        \omega \colon X_0 \to A,
    \end{align*}
    satisfying 
    \begin{align*}
        \frac{\omega(x)}{\omega(x')} = \frac{\eta_1(\phi)}{\eta_2(\phi)}
    \end{align*}
    for $\phi \colon x \to x' \in X_1$. 
\end{prop}

\begin{rmk}
    Horizontal and vertical compositions of 2-morphisms are both given by pointwise multiplication in $A$. This is consistent with the interpretation of $\omega(x)$ as an isomorphism $h_1(x) \to h_2(x)$ in $\cG$ for each $x \in X_0$. 
\end{rmk}

\begin{defn}\label{defn: normalized}
A pair $(\rho, \gamma)$ is a \emph{normalized} $\cG$-bundle if 
\[ \gamma(\phi, \id_x) = \gamma(\id_x, \psi)= 1_A \] for all $x \in X_0$ and $\phi, \psi \in X_1$ with $s(\phi) = x=t(\psi)$. 
\end{defn}

The following result is inspired by the ``semi-strictification'' results of \cite[Propositions 25 and 26]{Bartels}; we are able to give a simpler proof here because our 2-group is a skeletal category. 
\begin{lem} \label{lem: normalize gamma}
Let $\cG = \cG(G, A , \alpha)$ be as in Remark~\ref{rmk:normalized2} and  $X = \lbrace X_1 \rightrightarrows X_0\rbrace$ be a Lie groupoid. Let $(\rho, \gamma)$ be as in Proposition \ref{prop: remove hats}. Then $(\rho, \gamma)$ is isomorphic to a normalized $\cG$-bundle.
\end{lem}

\bp
We define a new pair $(\rho, \gamma')$ which will be a normalized $\cG$-bundle over $X$, isomorphic to the original pair $(\rho, \gamma)$.

We keep $\rho$ as before, but we define $\gamma' \colon X_2 \to A$ by $\gamma'(\phi, \psi) = \gamma(\phi,\psi)\epsilon(s(\phi))$. (As before, $\epsilon(s(\phi)) = \gamma(\phi, \id_{s(\phi)})^{-1} = \gamma(\id_{t(\psi)}, \psi)^{-1}$.) It is then clear that $\gamma'(\id_x, \psi) = \gamma(\id_x, \psi) \epsilon(x) = 1_A$ and $\gamma'(\phi, \id_x) = \gamma(\phi, \id_x) \epsilon(x) = 1_A$, so the normalization condition is satisfied. 

We need to show that $d\gamma' = \rho^* \alpha$, so let us take $(\phi, \psi, \xi) \in X_3$. We have 

\begin{align*}
    d\gamma'(\phi,\psi,\xi) & = \frac{\gamma'(\psi,\xi)\gamma'(\phi,\psi\circ \xi)}{\gamma'(\phi \circ \psi,\xi)\gamma'(\phi,\psi)} =\frac{\gamma(\psi,\xi) \epsilon(s(\psi)) \gamma(\phi,\psi\circ \xi) \epsilon(s(\phi))}{\gamma(\phi\circ \psi,\xi) \epsilon(s(\phi\circ \psi)) \gamma(\phi,\psi) \epsilon(s(\phi))} \\
    & = d\gamma(\phi,\psi,\xi) \frac{\epsilon(s(\psi)) \epsilon (s(\phi))}{\epsilon(s(\psi)) \epsilon(s(\phi))} = \alpha(\rho(\phi), \rho(\psi), \rho(\xi)),
\end{align*}
as desired. Therefore $(\rho, \gamma')$ is a normalized $\cG$-bundle. It remains to construct an isomorphism $(\rho, \gamma) \to (\rho, \gamma')$. 

We take $h \colon X_0 \to G$ to be the constant function $h(x) = 1_G$, and we take $\eta \colon X_1 \to A$ to be $\eta(\phi) = \epsilon(s(\phi))^{-1}$. Because $h(x) = 1_G$ and $\alpha$ is normalized, the condition (\ref{eq:eta condition}) for $\eta$ simplifies to 
\begin{align*}
\frac{\eta(\phi)\eta(\psi)}{\eta(\phi \circ \psi)} = \frac{\gamma(\phi,\psi)}{\gamma'(\phi,\psi)}, (\phi, \psi) \in X_2.
\end{align*}
It is easy to see that with our definitions of $\eta$ and $\gamma'$, both sides of the above equation give $\epsilon(s(\phi))^{-1}$, and hence $(1_G, \eta)$ gives the desired isomorphism $(\rho, \gamma) \to (\rho, \gamma')$. 
\ep

To summarize, starting with the bicategory $\Bun_\cG^\pre(X)=\Fun(X, \pt\sq \cG)$ from Definition~\ref{defn:prinbund}, for $\cG$ of a preferred form (Remark~\ref{rmk:normalized2}) we have shown that objects are the same as data $(\rho,\gamma)$, 1-morphisms are the same as data $(h,\eta)$ and 2-morphisms are the same as data $\omega$. Furthermore, there is an equivalent full and faithful subcategory of normalized pairs $(\rho, \gamma)$. Hereafter, we will use the notation $\Bun_\cG^\pre(X)$ to denote this description of the bicategory in terms of explicit normalized \v{C}ech data, whereas $\Fun(X, \pt\sq \cG)$ will refer to Definition~\ref{defn:prinbund}. With this convention the above discussion gives a canonical equivalence $\Bun_\cG^\pre(X)\xhookrightarrow{\sim}\Fun(X, \pt\sq \cG)$. 

\begin{lem}
\label{lem:underlyingGbundle2}
In terms of \v Cech cocycle data, the 2-functor~\eqref{eq:underlyingGbundle}, $\Bun_\cG^\pre(X)\to \Bun_G^\pre(X)$, sends $(\rho,\gamma)$ to $\rho$, and $(h,\eta)$ to $h$. \end{lem}

\bp
As in Lemma~\ref{lem:underlyingGbundle}, the morphism $\cG \to G$ of 2-groups induces a 2-functor $\pt \sq \cG \to \pt \sq G$ of bicategories, which in turn induces a 2-functor 
\begin{align*}
    \Bun_\cG^\pre(X) \xhookrightarrow{\simeq} \Fun(X, \pt \sq \cG) \to \Fun(X, \pt \sq G)=\Bun_G^\pre(X)
\end{align*}
of bicategories with the following values
\begin{align*}
    \textbf{On objects: }& (\rho, \gamma) \mapsto \rho; \\
    \textbf{On 1-morphisms: } & (h, \eta) \mapsto \eta; \\
    \textbf{On 2-morphisms: } & \hspace{1.6em}\omega \mapsto *.\\
\end{align*}
This proves the lemma. 
\ep 

\subsection{Flat gerbes as $\pt\sq A$-bundles}\label{ex: gerbes}

Recall the 2-group $\text{pt}/\!/ A$ from Example~\ref{ex:ptmodA}. 

\begin{defn}[{\cite[Example 73]{SP11} and \cite[\S4-5]{BehrendXu}}]\label{defn:gerbe}
For a Lie groupoid $X$, define the \emph{bicategory of flat $A$-gerbes that trivialize on $X_0$} as
$$
\Gerbe_A^\pre(X):=\Bun^\pre_{\pt\sq A}(X).
$$
\end{defn}

In light of Propositions \ref{prop: remove hats}, \ref{prop: remove hats for 1-morphisms}, \ref{prop: remove hats for 2-morphisms}, and Lemma \ref{lem: normalize gamma}, the bicategory $\Gerbe^{\text{pre}}_A(X)$ can be described as follows:

\begin{description}
\item[Objects] given by $\gamma \colon X_2 \to A$ a locally constant function such that $d\gamma = 1_A$, and $\gamma(\phi,\psi) = 1_A$ for $(\phi,\psi) \in X_2$ with either $\phi$ or $\psi$ equal to $\id_x, x \in X_0$. 
\item[1-morphisms $\gamma_1 \to \gamma_2$] given by $\eta \colon X_1 \to A$ a locally constant function such that for any $(\phi, \psi) \in X_2$
\begin{align*}
    \frac{\eta(\phi)\eta(\psi)}{\eta(\phi \circ \psi)} = \frac{\gamma_1(\phi, \psi)}{\gamma_2(\phi, \psi)}.
\end{align*}
\item[2-morphisms $\eta_1 \to \eta_2$] given by $\omega \colon X_0 \to A$ such that for any $\phi \colon x \to x'$ in $X_1$, 
\begin{align*}
    \frac{\omega(x)}{\omega(x')} = \frac{\eta_1(\phi)}{\eta_2(\phi)}. 
\end{align*}
\end{description}

In particular, if $M$ is a manifold and $Y \to M$ is a surjective submersion, a flat $A$-gerbe on $\check{C}(Y)$ in the sense of Definition~\ref{defn:gerbe} is given by a locally constant function $\gamma \colon \check{C}(Y)_2 \cong Y^{[3]}\to A$. More geometrically, $\gamma$ defines an isomorphism $p_{13}^*P \to p_{12}^*P \otimes_A p_{23}^*P$ between flat $A$-bundles over $Y^{[3]}$ where $P = Y^{[2]} \times A \to Y^{[2]}$ is the trivial $A$-bundle. The cocycle condition for~$\gamma$ implies a compatibility condition among the pullbacks of this isomorphism to $Y^{[4]}$.  
This recovers an $A$-gerbe in the classical definition of~\cite{Murray} with the additional condition that the $A$-bundles on $Y^{[2]}$ are flat.

\section{The bicategory of flat principal 2-group bundles}\label{sec: localization}

For the remainder of the paper, $\cG= \cG(G,A, \alpha)$ is a 2-group as in Remark~\ref{rmk:normalized2}, i.e., a categorical central extension of $G$ by $\pt\sq A$ classified by $[\alpha]\in \rmH^3(G;A)$. 

In the previous section we constructed the bicategory $\Bun_\cG^\pre(X)$ for a Lie groupoid $X$. In this section, we first construct the bicategory $\Bun_\cG^\pre$ whose objects are flat $\cG$-bundles on a (varying) Lie groupoid $X$ that are trivial on $X_0$, and whose 1-morphisms are maps of $\cG$-bundles covering functors between Lie groupoids. We then localize $\Bun_\cG^\pre$ to obtain a bicategory $\Bun_\cG$ whose objects are zig-zags~\eqref{eq:Gbundasfunctor0}. The (homotopy) fiber of $\Bun_\cG$ at a Lie groupoid $X$ therefore recovers the correct bicategory of flat $\cG$-bundles over the smooth stack $[X]$. We show that when $\cG=G$ is an ordinary group and $X$ is a manifold, this construction recovers the usual category of principal $G$-bundles. Similarly when $\cG=\pt\sq A$ is a one-object 2-group, we recover the bicategory of gerbes. We also compare our flat $\cG$-bundles to principal $\cG$-bundles in the sense of \cite[Definition 70]{SP11}.

\subsection{Pulling back $\mathcal{G}$-bundles}

Let $X = \{X_1 \rightrightarrows X_0\}, Y = \{Y_1 \rightrightarrows Y_0\}$ be Lie groupoids, and let $f=(f_1 \colon X_1 \to Y_1, f_0 \colon X_0 \to Y_0)$ be a smooth functor $X \to Y$. For any $n \ge 2$, we let $f_n \colon X_n \to Y_n$ denote the induced map. 

\begin{defn}\label{def: pullback}
    We define a functor \begin{align*}
    f^* \colon \Bun_\cG^\pre(Y) \to \Bun_\cG^\pre(X).
    \end{align*}
as follows:
\begin{align*}
    \textbf{On objects: } & f^*(\rho, \gamma) = (\rho \circ f_1, \gamma \circ f_2);\\
\textbf{On 1-morphisms: } &f^*(h, \eta) = (h \circ f_0, \eta \circ f_1);\\
\textbf{On 2-morphisms: } &f^*\omega = \omega \circ f_0.
\end{align*}
It is straightforward to check that this is a well-defined functor of bicategories. Moreover, for composable smooth functors $f, f'$, we have a strict equality $(f \circ f')^* = (f')^* \circ f^*$.

\end{defn}

\begin{rmk}
    From now on we will leave out the subscripts on the components of $f$ whenever they are clear from context, writing for example $f^*(\rho, \gamma) = (\rho \circ f, \gamma \circ f)$. Instead, the subscripts $f_1$ and $f_2$ will denote two different smooth functors, as in the following proposition. 
\end{rmk}

\begin{prop}\label{prop:pullbackmagic}
    Let $\tau \colon f_1 \Rightarrow f_2$ be a natural transformation of smooth functors. It induces a natural transformation of functors of bicategories $\widetilde{\tau} \colon f_1^* \Rightarrow f_2^*$. 
\end{prop}

\begin{proof}
Let $(\rho, \gamma)$ be an object in $\Bun_\cG^\pre(Y)$. We need to define a 1-morphism
$\widetilde{\tau}_{(\rho, \gamma)} \colon f_1^*(\rho, \gamma) \to f_2^*(\rho, \gamma)$ in $\Bun_\cG^\pre(X)$. For $x \in X_0$, we have $\tau(x) \in Y_1$ giving an isomorphism in the groupoid $Y$:
\begin{align*}
    \tau(x) \colon f_1(x) \to f_2(x).
\end{align*}
We use this to define $\widetilde{\tau}_{(\rho, \gamma)}$ to be the pair $(h_\tau, \eta_\tau)$, where 
\begin{align*}
    h_\tau \colon X_0  \to G, \quad 
    x \xmapsto{h_\tau} \rho(\tau(x));\qquad 
    \eta_\tau \colon X_1,  \to A \quad 
    \phi \xmapsto{\eta_\tau} \frac{\gamma(\tau(y), f_1(\phi))}{\gamma(f_2(\phi), \tau(x))}.
\end{align*}
Furthermore, given a 1-morphism $(h, \eta) \colon (\rho_1, \gamma_1) \to (\rho_2, \gamma_2)$ in $\Bun_\cG^\pre(Y)$, we need a 2-morphism~$\omega$ in $\Bun_\cG^\pre(X)$ making the following diagram commute:
\begin{center}
    \begin{tikzcd}
        f_1^*(\rho_1,\gamma_1) \arrow{d}{f_1^*(h, \eta)} \arrow{r}{\widetilde{\tau}_{(\rho_1, \gamma_1)}} & f_2^*(\rho_1,\gamma_1) \arrow{d}{f_2^*(h, \eta)} \\
         f_1^*(\rho_2,\gamma_2) \arrow{r}{\widetilde{\tau}_{(\rho_2, \gamma_2)}} & f_2^*(\rho_2, \gamma_2). \\
    \end{tikzcd}
\end{center}
We take $\omega \colon X_0 \to G, \omega(x) = \eta(\tau(x))$. It is straightforward to check that these data define a natural transformation as required. 
\end{proof}

With the structure of pullback, the bicategories $\Bun_\cG^\pre(X)$ fit together into a bicategory $\Bun_\cG^\pre$ over Lie groupoids fibered in 2-groupoids as follows. 
\begin{defn}
Define the bicategory $\Bun_\cG^\pre$ as follows. 

\textbf{Objects:} a Lie groupoid $X$, and a pair $(\rho, \gamma) \in \Bun_\cG^\pre(X)$. We often denote this data by $(\rho, \gamma)_X$ to emphasize the base Lie groupoid. 

\textbf{1-morphisms:} a smooth functor $f \colon X \to Y$ of Lie groupoids, together with a 1-morphism $(h, \eta) \colon (\rho_1, \gamma_1) \to f^*(\rho_2, \gamma_2)$ in $\Bun_\cG^\pre(X)$. 

\textbf{2-morphisms:} given 1-morphisms $(f_1, h_1, \eta_1), (f_2, h_2, \eta_2)\colon (\rho_1, \gamma_1)_X \to (\rho_2, \gamma_2)_Y$, a 2-morphism $(f_1, h_1, \eta_1) \to (f_2, h_2, \eta_2)$ is given by a natural transformation $\tau \colon f_1 \Rightarrow f_2$ of smooth functors, together with a 2-morphism $\omega$ in $\Bun_\cG^\pre(X)$ between the 1-morphisms 
\begin{align*}
    \widetilde{\tau}_{(\rho_2, \gamma_2)} \circ (h_1, \eta_1) \text{ and } (h_2, \eta_2) \colon (\rho_1, \gamma_1) \to f_2^*(\rho_2, \gamma_2). 
\end{align*}
\end{defn}

We spell out some of the compatibility conditions implied by the above definition. A 1-morphism $(\rho_1, \gamma_1)_X \to (\rho_2, \gamma_2)_Y$ is the data of a triple $(f, h, \eta)$, where for $\phi \colon x \to x'$ a morphism in $X$ we view $\eta(\phi) \in A$ as giving an isomorphism $\rho_2(\phi)h(x) \to h(x')\rho_1(\phi)$ in $\cG$. This leads to the condition on $h$
\begin{align*}
    \rho_2(f(\phi))h(x) = h(x')\rho_1(\phi) \in G.
\end{align*}
Furthermore, the compatibility condition for $\eta$ in~\eqref{eq:eta condition} becomes 
\begin{align*}
    \frac{\eta(\phi) \eta(\psi)}{\eta(\phi \circ \psi)} = \frac{\gamma_1(\phi, \psi)}{\gamma_2(f(\phi), f(\psi))} \frac{\alpha(\rho_2(f(\phi)), h(x'), \rho_1(\psi))}{\alpha(\rho_2(f(\phi)), \rho_2(f(\psi)), h(x)) \alpha(h(x''), \rho_1(\phi), \rho_1(\psi))},
\end{align*}
for $x \xrightarrow{\psi} x' \xrightarrow{\phi} x''$ in $X$.

Given two composable 1-morphisms $(f_1,h_1, \eta_1) \colon (\rho_1, \gamma_1)_X \to (\rho_2, \gamma_2)_Y$ and $(f_2, h_2, \eta_2) \colon (\rho_2, \gamma_2)_Y \to (\rho_3, \gamma_3)_Z$, the composition is $(f_2 \circ f_1, f_1^*(h_2, \gamma_2) \circ (h_1, \gamma_1))$. Explicitly, its components are the following three maps:
\beq
  &f_2 \circ f_1 \colon X \to Z; \quad \widetilde h \colon X_0 \to G, \quad x  \xmapsto{\tilde{h}} h_2(f_1(x))h_1(x); \nonumber \\[5pt]\nonumber 
   &\widetilde \eta \colon X_1 \to A, \\[5pt]\nonumber 
   &(\phi \colon x \to x') \xmapsto{\tilde{\eta}} \eta_2(f_1(\phi))\eta_1(\phi) \frac{\alpha(h_2(f_1(x')), \rho_2(f_1(\phi)), h_1(x))}{\alpha(\rho_2(f_1(\phi)), h_2(f_1(x)), h_1(x))\alpha(h_2(f_1(x')), h_1(x'), \rho_1(\phi))}.\nonumber
\eeq

The following observation is simple but useful: 
\begin{lem}\label{lem: factor 1-morph}
    Every 1-morphism $(f,h, \eta) \colon (\rho_1, \gamma_1)_X \to (\rho_2, \gamma_2)_Y$ in $\Bun_\cG^\pre$ factors as a composition of a morphism in $\Bun_\cG^\pre(X)$ and one of the form $(f, 1_G, 1_A) = (f, \id_{f^*(\rho_2, \gamma_2)})$:
\begin{align*}
    (\rho, \gamma)_X \xrightarrow{(\id_X, h, \eta)} (\rho_2 \circ f, \gamma_2 \circ f)_X \xrightarrow{(f, 1_G, 1_A)} (\rho_2, \gamma_2)_Y.
\end{align*}
\end{lem}

We now unwind the data and conditions satisfied by 2-morphisms. Given two 1-morphisms $(f_1, h_1, \eta_1), (f_2, h_2, \eta_2)$ between $(\rho_1, \gamma_1)_X \to (\rho_2, \gamma_2)_Y$, a 2-morphism $(f_1, h_1, \eta_1) \to (f_2, h_2, \eta_2)$ is a pair $(\tau, \omega)$ where $\tau \colon f_1 \Rightarrow f_2$ is a natural transformation such that $\rho_2(\tau(x))h_1(x) = h_2(x)$ for all $x \in X_0$; furthermore $\omega \colon X_0 \to A$ provides a family of isomorphisms witnessing this identity, satisfying 
\beq \label{compat for 2-morph}
  &&\resizebox{.9\textwidth}{!}{$\frac{\omega(x)}{\omega(x')} = \frac{\eta_1(\phi)}{\eta_2(\phi)}\frac{\gamma_2(\tau(x'), f_1(\phi))}{\gamma_2(f_2(\phi), \tau(x))} \frac{\alpha(\rho_2(\tau(x')), \rho_2(f_1(\phi)), h_1(x))}{\alpha(\rho_2(\tau(x')), h_1(x'), \rho_1(\phi))\alpha(\rho_2(f_2(\phi)), \rho_2(\tau(x)), h_1(x))}$},
\eeq
for $\phi \colon x \to x'$ in $X$. 

Suppose we also have two 1-morphisms $(f_3, h_3, \eta_3), (f_4, h_4, \eta_4) \colon (\rho_2, \gamma_2)_Y \to (\rho_3, \gamma_3)_Z$ and a 2-morphism $(\tau_2, \omega_2) \colon (f_3, h_3, \eta_3)\to (f_4, h_4, \eta_4) $. The horizontal composition $(\tau_2, \omega_2) \star (\tau, \omega)$ is a 2-morphism $(f_3, h_3, \eta_3) \circ (f_1, h_1, \eta_1) \to (f_4, h_4, \eta_4) \circ (f_2, h_2, \eta_2)$. It is the pair consisting of the natural transformation $\tau_2 \star \tau \colon f_3 \circ f_1 \Rightarrow f_4 \circ f_2$ given by the horizontal composition of $\tau_2$ and $\tau$, and the locally constant map $\omega_2 \star \omega \colon X_ 0 \to A$ given by 
\beq \label{eq: horiz comp 2-morph}
&& \omega_2 \star \omega(x)= \omega_2(f_1(x))\omega(x)\eta_4(\tau(x)) \gamma_3(f_4(\tau(x)), \tau_2(f_1(x))) \\[5pt]
 && \frac{\alpha(\rho_3(f_4(\tau(x))), \rho_3(\tau_2(f_1(x))), h_3(f_1(x))h_1(x))\alpha(h_4(f_2(x)), \rho_2(\tau(x)), h_1(x))}{\alpha(\rho_3(\tau_2(f_1(x))), h_3(f_1(x)), h_1(x))\alpha(\rho_3(f_4(\tau(x))), h_4(f_1(x)), h_1(x))}.\nonumber  
\eeq

Similarly, given a third 1-morphism $(f_3, h_3, \eta_3) \colon (\rho_1, \gamma_1)_X \to (\rho_2, \gamma_2)_Y$ and a 2-morphism $(\tau_2, \omega_2) \colon (f_2, h_2, \eta_2) \to (f_3, h_3, \eta_3)$, we obtain the vertical composition $(\tau_2, \omega_2) \circ (\tau, \omega) \colon (f_1, h_1, \eta_1) \to (f_3, h_3, \eta_3)$. It is the pair consisting of the natural transformation $\tau_2 \circ \tau \colon f_1 \to f_3$ given by the vertical composition of $\tau_2$ and $\tau$, and the locally constant map $\omega_2 \circ \omega \colon X_0 \to A$ given by
\begin{align*}
    (\omega_2 \circ \omega)(x) = \omega_2(x) \omega(x) \alpha(\rho_2(\tau_2(x)), \rho_2(\tau(x)), h_1(x)) \gamma_2(\tau_2(x), \tau(x)).
\end{align*}

\begin{rmk} The above defines a prestack of $\cG$-bundles on Lie groupoids. There are various machines that implement 2-stackification; we will use the calculus of right fractions from~\cite{Pronk}, though see also~\cite{NikolausSchweigert}. 
\end{rmk}

\begin{thm}\label{thm:pronk}
    The bicategory $\Bun_\cG^\pre$ admits a calculus of right fractions with respect to the class of 1-morphisms $\calw$ of the form $(f,h,\eta)$ where $f$ is an essential equivalence. Let $\Bun_\cG:=\Bun_\cG^\pre[\calw^{-1}]$ denote the localized bicategory. 
\end{thm}

We defer the proof of Theorem~\ref{thm:pronk} to Appendix~\ref{sec: localization app}.

\begin{rmk}\label{rmk: refine covers}
By construction, the bicategory $\Bun_\cG$ has the property that for any object $(\rho, \gamma)_X$ and for any essential equivalence $f \colon Y \to X$, the 1-morphism $(f, 1_G, 1_A) \colon f^*(\rho, \gamma) \to (\rho, \gamma) $ is an equivalence. 
\end{rmk}

\begin{rmk}
    The 2-functor $\Bun_\cG^\pre \to \LieGrpd$ sends 1-morphisms in $\calw$ to essential equivalences. By the universal property of localization~\cite[Theorem 21]{Pronk}, we obtain a 2-functor $\Bun_\cG \to \LieGrpd[W^{-1}]$, where $\sf{W}$ is the class of essential equivalences. The bicategory $\LieGrpd[W^{-1}]$ is equivalent to the bicategory of differentiable stacks. Furthermore, morphisms in $\LieGrpd[W^{-1}]$ can be described using \emph{bibundles}~\cite{SP11, Lerman}; we do not use these structures here.
\end{rmk}

\subsection{The bicategory ${\sf Bun}_{\mathcal{G}}(X)$.}\label{subsec: BuncG(X)}

\begin{defn} \label{eq:defncfBund}
For a Lie groupoid $X$, define the bicategory $\Bun_\cG(X)$ of \emph{principal $\cG$-bundles on $X$} as the homotopy fiber of $\Bun_\cG\to \LieGrpd[{\sf W}^{-1}]$ at the object $X\in \LieGrpd[{\sf W}^{-1}]$.
\end{defn}

\begin{cor}
The bicategory of principal $\cG$-bundles forms a 2-stack over the bicategory $\LieGrpd[W^{-1}]$. 
\end{cor}
\bp
After localization, the homotopy fibers of $\Bun_\cG\to \LieGrpd[{\sf W}^{-1}]$ necessarily satisfy descent.
\ep

By definition, an object of the homotopy fiber is a pair $(f \colon Y \to X, (\rho, \gamma)_Y)$ consisting of an invertible 1-morphism $f$ in $\LieGrpd[{\sf W}^{-1}]$ and an object $(\rho, \gamma)_Y$ of $\Bun_\cG$. The invertible morphism $f$ in $\LieGrpd[{\sf W}^{-1}]$ is a zig-zag of essential equivalences, $Y\xleftarrow{g} Z\to X$, and the composition $f \circ g\colon Z\to X$ in $\LieGrpd[W^{-1}]$ is naturally isomorphic to the right leg of this zig-zag. Then following Remark \ref{rmk: refine covers}, the pair $(f \colon Y \to X, (\rho, \gamma)_Y)$ is equivalent to the pair $(f \circ g, g^*(\rho, \gamma))$. Therefore, $\Bun_\cG(X)$ is equivalent to a bicategory whose objects are zig-zags
      \beq\label{eq:zigzag2}
        X \xleftarrow{\;f\;} Y \xrightarrow{(\rho, \gamma)} \pt \sq \cG.
    \eeq
 with $f$ an essential equivalence (c.f. Definition \ref{defn:prinbund0}). We will denote the zig-zag by $(f, (\rho, \gamma)_Y)$ below. 

Similarly, a 1-morphism $(f_1, (\rho_1, \gamma_1)_{Y_1}) \to (f_2, (\rho_2, \gamma_2)_{Y_2}) $ in the homotopy fiber is given up to 2-isomorphism by the following data: $(Z, g_1, g_2, \tau, (h, \eta))$, where $g_i \colon Z \to Y_i, i = 1,2$ are essential equivalences, $\tau$ is a smooth natural transformation $f_1 \circ g_1 \Rightarrow f_2 \circ g_2$, and $(h, \eta)$ is a 1-morphism in $\nohatG{Z}$ between $g_1^*(\rho_1, \gamma_1)$ and $g_2^*(\rho_2, \gamma_2)$. 

Finally, let $(Z_j, g_{1j}, g_{2j}, \tau_j, (h_j, \eta_j)), j = 3,4$ be two 1-morphisms of this form between objects $(f_1, (\rho_1, \gamma_1)_{Y_1})$ and $(f_2, (\rho_2, \gamma_2)_{Y_2})$. Unwinding the definitions, we see that a 2-morphism between them is given by data $(Z, g_3, g_4, \sigma_1, \sigma_2, \omega)$, where $g_j \colon Z \to Z_j$ are essential equivalences, $\sigma_i \colon g_{i3} \circ g_3 \Rightarrow g_{i4} \circ g_4$, $i=1,2$ are smooth natural transformations such that $\sigma_2 \circ \tau_3 = \tau_4 \circ \sigma_1$, and $\omega$ is a 2-morphism in $\nohatG{Z}$ between $\tilde{\sigma}_2 \circ g_3^*(h_3, \eta_3)$ and $g_4^*(h_4, \eta_4) \circ \tilde{\sigma}_1$. 
Two sets of data of the form $(Z, g_3, g_4, \sigma_1, \sigma_2, \omega)$ represent the same 2-morphism if they agree ``upon refinement'', that is, after being pulled back along compatible essential equivalences. 

In particular, all 1-morphisms and 2-morphisms are invertible. 

\begin{notation}
    Hereafter, $\Bun_\cG(X)$ will denote the bicategory equivalent to the one in Definition~\ref{eq:defncfBund} whose objects are zig-zags~\eqref{eq:zigzag2}, and 1-morphisms and 2-morphisms are described above. 
\end{notation}

We complete this subsection with a comparison of the bicategory $\Bun_\cG(X)$ with previously studied special cases.

\begin{ex}\label{ex:ordinaryGbund}
Let $G$ be an ordinary group, viewed as a 2-group, and let $M$ be a smooth manifold regarded as a Lie groupoid. We check that the bicategory $\Bun_G(M)$ just defined is in fact an ordinary category, and agrees with the usual category of principal $G$-bundles in terms of \v{C}ech cocycle data, as discussed in \S\ref{subsec: ordinary G-bundles}. Recall from Example \ref{ex: ess equiv to manifold} that the only essential equivalences to $M$ are of the form $\widetilde{f} \colon \check{C}(Y) \to M$ for $f \colon Y \to M$ a surjective submersion. With this in mind, we see that an object of $\Bun_G(M)$ is given by a pair $(f \colon Y \to M, \rho \colon Y \times_M Y \to G)$ consisting of a surjective submersion and a locally constant 1-cocycle; i.e., a principal $G$-bundle on the manifold $M$. 

To unpack the notion of 1-morphisms, we first note that if $\widetilde g \colon \widetilde Z \to \check{C}(Y)$ is an essential equivalence, then $\widetilde f \circ \widetilde g$ is an essential equivalence from the Lie groupoid $\widetilde Z$ to the manifold $M$. Then as in Example \ref{ex: ess equiv to manifold}, if we let $Z = \widetilde Z _0$ and $g = \widetilde g_0$ we have that $f \circ g \colon Z \to M$ is a surjective submersion, $\widetilde Z$ is the associated \v Cech groupoid $\check{C}(Z)$, and $\widetilde g$ is the essential equivalence $\check{C}(Z) \to \check{C}(Y)$ induced by the smooth map $g$. We also note that the category of morphisms from a Lie groupoid into a manifold is in fact a set; there are no non-identity 2-morphisms. These observations allow us to see that a morphism in $\Bun_G(M)$ from $(f_1, \rho_1)$ to $(f_2, \rho_2)$ consists of a smooth manifold $Z$ with smooth morphisms $g_i \colon Z \to Y_i$ for $i = 1,2$, such that $f_1 \circ g_1 = f_2 \circ g_2$ is a surjective submersion $Z \to Y$, together with a locally constant 0-cochain $h \colon Z \to A$ such that $\rho_2(g_2(z), g_2(z'))h(z') = h(z)\rho_1(g_1(z), g_1(z'))$. (Here $(z, z') \in Z \times_M Z$ is viewed as a morphism $z' \to z$ in $\check{C}(Z)$.) In other words, $h$ is a 1-morphism of principal $G$-bundles between the pullbacks of $\rho_1$ and $\rho_2$ to $Z$. 

Finally, there are no non-identity 2-morphisms. 

\end{ex}

\begin{ex}\label{ex:Hequivprin}
    Let $G$ and $H$ be ordinary groups and $X = M \sq H$ an action groupoid. We claim that $\Bun_G(M\sq H)$ is the category whose objects are flat principal $G$-bundles $p\colon P\to M$ where $P$ is equipped with an $H$-action compatible with the $H$-action on $M$ and commuting with the $G$-action i.e., $H$-equivariant $G$-bundles on $M$. Indeed, the smooth functor $M\to M\sq H$ of Lie groupoids induces a pullback functor 
    \beq\label{ex:HequivprinG}
\Bun_G(M\sq H)\to \Bun_G(M), 
    \eeq
and so Example~\ref{ex:ordinaryGbund} extracts an ordinary $G$-bundle on $M$ from an object in $\Bun_G(M\sq H)$. Conversely, lifting an object $P\in \Bun_G(M)$ along the functor~\eqref{ex:HequivprinG} requires descent data for $P$ along the map $M\to M\sq H$. When $P$ is defined relative to a surjective submersion $f\colon Y\to M$, this descent data is equivalent to a $G$-bundle over the pullback groupoid $f^*(M\sq H)$ (see Example~\ref{ex:pullbackgrpod}) that trivializes on $Y\to M\to M\sq H$. A straightforward computation verifies that this descent data does indeed determine an $H$-action on $P$. 
\end{ex}

\begin{ex}\label{ex: gerbes localized}
Applying the localization construction to Example \ref{ex: gerbes} for a fixed abelian Lie group $A$, we see that gerbes on a Lie groupoid form a bicategory $\Bun_{\pt\sq A}(X) = \Gerbe_A(X)$. Concretely, an object is given by the data of a pair $(f \colon Y \EquivTo X, \gamma \colon Y_2 \to A)$ where $f$ is a smooth essential equivalence and $\gamma$ defines a gerbe on $Y$ trivial over $Y_0$. 

The constructions of \S\ref{sec: localization} recover the fact that the assignment $X\mapsto \Gerbe_A(X)$ is a 2-stack, as was already known (e.g., see~\cite[\S5]{NikolausWaldorf}). 
\end{ex}

Finally, we indicate how our definition of flat $\cG$-bundles compares with \cite[Definition 70]{SP11}. Let $(f \colon Y \to X, (\rho, \gamma)_Y)$ be a flat principal $\cG$-bundle. We define a Lie groupoid $\cP$ with objects $P_0 = Y_0 \times G$ and morphisms $P_1= Y_1 \times G \times A$, where $(\phi \colon y \to y', g, a)$ has source $(y, g)$ and target $(y', \rho(\phi)g)$. Composition is given by 
    $$
    (\phi, \rho(g), a) \circ (\psi, g, b) = (\phi \circ \psi, g, ab\gamma(\phi, \psi)\alpha(\rho(\phi), \rho(\psi), g)),
    $$
    and the identity morphism of an object $(y, g)$ is $(\id_y, g, 1_A)$.

    \begin{lem}\label{lem: gluing 2-bundles}
        This defines a Lie groupoid $\cP = \lbrace P_1 \rightrightarrows P_0\rbrace$ over $X$, equipped with a smooth functor ${\sf act} \colon \cP \times \cG \to \cP$.
    \end{lem}
    \begin{proof}
        The properties of a Lie groupoid follow from the cocycle and normalization conditions on $\alpha, \rho$ and $\gamma$. The smooth functor $\pi \colon \cP \to X$ is given by $(y, g) \mapsto f(y)$ on objects and $(\phi, g, a) \mapsto f(\phi)$ on morphisms. 

        The smooth functor ${\sf act}$ is given on objects by $((y,g),h) \mapsto (y,gh)$ and on morphisms by $((\phi, g, a),(h,b)) \mapsto (\phi, gh, \frac{ab}{\alpha(\rho(\phi), g,h)})$. Functoriality also follows from the properties of $\alpha, \rho,$ and $\gamma$. 
    \end{proof}

To complete the comparison with \cite{SP11}, 
one must extend the functor ${\sf act}$ to an action of the 2-group $\cG$ on $\cP$ with a equivariant isomorphism over the (trivializing) cover $Y\to X$. This is not difficult, but requires introducing higher categorical actions that are beyond the scope of the present paper. 

\subsection{From $\mathcal{G}$-bundles to trivializations of 2-gerbes}\label{sec:trivgerbe}

There is a natural forgetful functor from principal $\cG$-bundles to principal bundles for the underlying ordinary group $G$. The goal of this section is to describe the fibers of this forgetful functor; this will use the language of 2-gerbes reviewed in~\S\ref{sec:2gerbe}.

\begin{lem}\label{lem:forgetfulfunctor}For a 2-group $\cG = \cG(G,A, \alpha)$ and a Lie groupoid $X$, the functor~\eqref{eq:underlyingGbundle}, $\Bun_\cG^\pre(X)\to \Bun_G^\pre(X)$, induces a forgetful functor 
\beq
\pi\colon \Bun_\cG(X)\to \Bun_G(X)\label{eq:forgetfulfunctor}
\eeq 
from the bicategory of $\cG$-bundles on $X$ to the category of $G$-bundles on $X$. \end{lem}

\bp 
The functor~\eqref{eq:underlyingGbundle} sends the class $\calw$ of 1-morphisms $(f, h, \eta)$ with $f$ an essential equivalence to the class $W$ of 1-morphisms $(f, h)$ with $f$ an essential equivalence. By the universal property of localization~\cite[Theorem 21]{Pronk}, we obtain a functor $\Bun_\cG \to \Bun_G$, as follows: 
\begin{align*}
    \textbf{On objects: }& (\rho, \gamma)_X \mapsto \rho_X; \\
    \textbf{On 1-morphisms: } & (f, h, \eta) \mapsto (f, h); \\
    \textbf{On 2-morphisms: } & (\tau, \omega) \mapsto \tau.
\end{align*} For any Lie groupoid $X$, we obtain a functor between homotopy fibers yielding the desired forgetful functor~\eqref{eq:forgetfulfunctor}; its values on objects are as above, while its values on zig-zags giving 1- and 2-morphisms are given by applying the above formulas to the components of the zig-zag.
\ep

The following lemma shows that~\eqref{eq:forgetfulfunctor} satisfies a path-lifting property: morphisms in $\Bun_G(X)$ can be lifted to 1-morphisms in $\Bun_\cG(X)$. 

\begin{lem}\label{lem: lifting 1-morphisms}
    Let $(f \colon Y \to X, (\rho_1, \gamma_1)_Y)$ be an object of $\Bun_\cG(X)$, and suppose that $(f, \rho_2)$ is an object of $\Bun_G(X)$ defined relative to the same essential equivalence $f\colon Y\to X$. If $h \colon \rho_1 \to \rho_2$ is a 1-morphism in $\nohat{G}{Y}$, then there exists locally constant cochains $\gamma_2 \colon Y_2 \to A$ and $\eta \colon Y_1 \to A$ such that $(\rho_2, \gamma_2)_Y$ is a principal $\cG$-bundle, and $(h, \eta)$ is a 1-isomorphism $(\rho_1, \gamma_1) \to (\rho_2, \gamma_2)$ in $\nohatG{Y}$. 
\end{lem}

\begin{proof}
The existence of $h \colon \rho_1 \to \rho_2$ tells us that for $\phi \colon y \to y'$ in the Lie groupoid $Y$, we have an equality $\rho_2(\phi) = \left( h(y') \rho_1(\phi)\right) h(y)^{-1}$ in $G$. 

To lift $\rho_2$ to a pair $(\rho_2, \gamma_2) \in \nohatG{Y}$, we need to provide for every composable pair of morphisms $y \xrightarrow{\psi} y' \xrightarrow{\phi} y''$ in $Y$ an isomorphism 
\begin{align*}
    \gamma_2(\phi, \psi) \colon \rho_2(\phi\circ \psi) \to \rho_2(\phi) \rho_2(\psi) 
\end{align*}
in $\cG$. That is, we need an isomorphism
\begin{align*}
    \left(h(y'')\rho_1(\phi \circ \psi) \right)h(y)^{-1} \to \left[\left(h(y'') \rho_1(\phi)\right)h(y')^{-1}\right]\left[\left(h(y') \rho_1(\psi)\right) h(y)^{-1}\right].
\end{align*}
The associator, the isomorphism $\gamma_1 \colon \rho_1(\phi \circ \psi) \to \rho_1(\phi) \circ \rho_1(\psi)$, and the morphism $e_{h(y')} \colon h(y')^{-1}h(y') \to 1_G$ from \eqref{eq:unitcounitnormalized} provide exactly the ingredients needed to do this, and we take
\begin{align*}
    \gamma_2(\phi, \psi) = \gamma_1(\phi, \psi) \frac{\alpha(h(y'), h(y')^{-1}, h(y')\rho_1(\psi))\alpha(\rho_2(\phi), h(y')\rho_1(\psi), h(y)^{-1})}{\alpha(h(y''), \rho_1(\phi), \rho_1(\psi))\alpha(h(y'')\rho_1(\phi), h(y')^{-1}, h(y')\rho_1(\psi))}.
\end{align*}
The cocycle condition on $\alpha$ and the compatibility condition $d\gamma_1 = \rho_1^*\alpha$ yield the analogous condition $d\gamma_2 = \rho_2^*\alpha$, so that $(\rho_2, \gamma_2)$ is indeed an object of $\nohatG{Y}$.

It remains to extend $h \colon \rho_1 \to \rho_2$ to an isomorphism $(h, \eta) \colon (\rho_1, \gamma_1) \to (\rho_2, \gamma_2)$. For $\phi \colon y \to y'$ in $Y$ we need to specify an isomorphism $\eta(\phi) \colon \rho_2(\phi)h(y) \to h(y')\rho_1(\phi)$. Just as before, we have all the ingredients we need to do this:
\begin{align*}
    \rho_2(\phi)h(y) = \left[(h(y')\rho_1(\phi))h(y)^{-1}\right]h(y) &\xrightarrow{\alpha} (h(y')\rho_1(\phi))(h(y)^{-1}h(y)) \\
    & \xrightarrow{e_{h(y)}} (h(y')\rho_1(\phi))1_G = h(y')\rho_1(\phi).
\end{align*}
That is, we take $\eta(\phi) = \frac{\alpha(h(y')\rho_1(\phi), h(y)^{-1}, h(y))}{\alpha(h(y), h(y)^{-1}, h(y))}$. 
Applying the cocycle condition to the tuple $(\rho_2(\phi), h(y), h(y)^{-1}, h(y))$, we obtain $\eta(\phi) = \alpha(\rho_2(\phi), h(y), h(y)^{-1})$. Again, we use the cocycle condition on $\alpha$, the compatibility condition of $\gamma_1$, and the definition of $\gamma_2$ to show that $(h, \eta)$ is indeed an isomorphism $(\rho_1, \gamma_1) \to (\rho_2, \gamma_2)$ as claimed. 
\end{proof}

We are now ready to describe the fibers of~\eqref{eq:forgetfulfunctor} in terms of 2-gerbes; we refer to~\S\ref{sec:2gerbe} for background on 2-gerbes. Our main examples are the following. 

\begin{ex}
    Let $G$ be a group and $A$ an abelian group. Then a 3-cocycle $\alpha \in Z^3(G;A)$ determines a 2-gerbe on $\pt\sq G$. We will again denote this 2-gerbe by~$\alpha$. 
\end{ex}

\begin{ex}\label{ex:2gerbefromP}
Let $X$ be a Lie groupoid. Given a flat $G$-bundle on $X$ 
\beq\label{eq:cocycledata2gerbe}
X\xleftarrow[]{\sim} Y\xrightarrow[]{\rho}\pt\sq G,
\eeq
define a 2-gerbe by
    \begin{align}\label{eq: 2-gerbe example}
        Y_3\xrightarrow{\rho^{\times 3}} G\times G\times G \xrightarrow{\alpha} A.
    \end{align}
We use the notation $\rho^*\alpha$ for this 2-gerbe that trivializes over the (fixed) cover $Y\to X$, compare Definition  \ref{def: 2-gerbe fixed base}.

Next we describe how $\rho^*\alpha$ changes under isomorphisms of $G$-bundles, and in particular, under changing the cover $Y\to X$: given an isomorphism between principal bundles, we construct an explicit 1-isomorphism between the corresponding 2-gerbes. Let $Z=Y_1\times_X Y_2$ and $h \colon Z_0 \to G$ be the data of the isomorphism $(f_1\colon Y_1\to X,\rho_1)\to (f_2\colon Y_2\to X,\rho_2)$. Let $p_i\colon Z\to Y_i$, $i=1,2$  denote the projections out of the fibered product. 
For $z \xrightarrow{\psi} z' \xrightarrow{\phi} z''$ in $Z$, define 
    \begin{align*}
        \gamma(\phi, \psi) = \frac{\alpha(\rho_2(p_2(\phi)), h(z'), \rho_1(p_1(\psi)))}{\alpha(\rho_2(p_2(\phi)), \rho_2(p_2(\psi)), h(z))\alpha(h(z'') ,\rho_1(p_1(\phi)), \rho_1(p_1(\psi)))}.
    \end{align*} 
    It follows from the cocycle condition and the compatibility condition between $h, \rho_1$ and $\rho_2$ that $\gamma$ gives the claimed 1-morphism $\rho_1^*\alpha\to \rho_2^*\alpha$. 
    Therefore, for any principal $G$-bundle $P$, we obtain a flat 2-gerbe which, up to isomorphism, does not depend on the choice of cocycle data $(f, \rho)$. We denote this flat 2-gerbe by~$\lambda_{P, \alpha}$. 
\end{ex}

\begin{prop}\label{prop: fiber of prestacks}
    Let $\cG = \cG(G, A, \alpha)$ be a 2-group and let $X$ be a Lie groupoid. For $\rho \colon X_1 \to G$ an object of $\Bun_G^\pre(X)$, the homotopy fiber ${\sf fib}(\rho)$ of $\nohatG{X} \to \nohat{G}{X}$ over $\rho$ is canonically equivalent to the bicategory of trivializations of the 2-gerbe $\rho^*\alpha$, constructed in Example~\ref{ex:2gerbefromP}. We emphasize that this is a 2-gerbe that trivializes on a fixed cover $Y\to X$, see Definition \ref{def: 2-gerbe fixed base}.
\end{prop}

\begin{proof}
    An object of ${\sf fib}(\rho)$ is given by an object $(\rho_1, \gamma_1)$ of $\nohatG{X}$ together with an isomorphism $h_1 \colon \rho_1 \to \rho$ in $\nohat{G}{X}$. A 1-morphism $((\rho_1, \gamma_1), h_1) \to ((\rho_2, \gamma_2), h_2)$ is given by a 1-morphism $(h, \eta) \colon (\rho_1, \gamma_1) \to (\rho_2, \gamma_2)$ such that the morphism $h \colon \rho_1 \to \rho_2$ is compatible with the morphisms $h_i \colon \rho_i \to \rho$ in $\nohat{G}{X}$. This compatibility condition amounts to the requirement that $h= h_2^{-1}h_1$ as cochains $X_0 \to G$. Finally, given two 1-morphisms $(h, \eta_1), (h, \eta_2) \colon (\rho_1, \gamma_1) \to (\rho_2, \gamma_2)$, a 2-morphism in  ${\sf fib}(\rho)$ is just a 2-morphism $\omega \colon (h, \eta_1) \to (h, \eta_2)$.

    Given an object $((\rho_1, \gamma_1), h)$, Lemma \ref{lem: lifting 1-morphisms} provides locally constant cochains $\gamma \colon X_2 \to G$ and $\eta \colon X_1 \to A$ such that $(h, \eta)$ gives an isomorphism $(\rho_1, \gamma_1) \to (\rho, \gamma)$. Then $((\rho, \gamma), 1_G)$ is an object of  ${\rm fib}(\rho)$, and $(h, \eta) \colon ((\rho_1, \gamma_1), h) \to ((\rho, \gamma), 1_G)$. So every object in  ${\sf fib}(\rho)$ is canonically isomorphic to one of the form $((\rho, \gamma), 1_G)$, and $\gamma$ satisfies $d\gamma = \rho^*\alpha$. In other words, $\gamma$ provides a trivialization of the 2-gerbe $\rho^*\alpha$. 

    Now let $((\rho, \gamma_1), 1_G), ((\rho, \gamma_2), 1_G)$ be two objects of ${\sf fib}(\rho)$. A 1-morphism between them is a pair $(h, \eta)$ with $h = 1_G^{-1}1_G = 1_G$; that is, it is just a locally constant cochain $\eta \colon X_1 \to A$ satisfying the condition \eqref{eq:eta condition}, which simplifies to $d\eta = \frac{\gamma_1}{\gamma_2}$, since $h = 1$. Thus, $\eta$ is a 2-morphism between the two trivializations $\gamma_1$ and $\gamma_2$. 

    Finally, let $(1_G, \eta_1), (1_G, \eta_2)$ give two 1-morphisms in ${\sf fib}(\rho)$; a 2-morphism between them is a locally constant cochain $\omega \colon X_0 \to A$ such that $d\omega = \frac{\eta_2}{\eta_1}$. This is a 3-morphism between the 2-morphisms $\eta_1, \eta_2$ in the 3-category $2\Gerbe_A^\pre(X)$. 
\end{proof}

We finally consider the homotopy fiber of the functor $\pi \colon \Bun_\cG(X) \to \Bun_G(X)$ over a principal $G$-bundle $P$ on $X$. We refer to it as the \emph{bicategory of $\cG$-bundles with underlying $G$-bundle $P$}, and denote it by $\Bun_\cG(X)_P$.

\begin{thm}  \label{thm:2gerbetriv}
Given the data as above, there is a canonical equivalence between the bicategory $\Bun_\cG(X)_P$ of $\cG$-bundles with underlying flat $G$-bundle $P$ and the bicategory ${\sf Triv}(\lambda_{P, \alpha})$ of trivializations of the flat 2-gerbe $\lambda_{P, \alpha}$. Furthermore,  in the event that these two bicategories are non-empty, there is a non-canonical equivalence between the bicategory $\Bun_\cG(X)_P$ and the bicategory $\Gerbe_A(X)$ of gerbes on $X$, where an equivalence is determined by a choice of principal $\cG$-bundle lifting $P$.
\end{thm}

\begin{proof}
    Let us choose cocycle data $(f\colon Y \to X, \rho)$ in $\Bun_G(X)$ representing the $G$-bundle~$P$. Stackification is a left adjoint, so commutes with homotopy fiber products. Hence $\Bun_\cG(X)_P$ is canonically equivalent to the stackification of the homotopy fiber ${\sf fib}(\rho)$ studied in Proposition \ref{prop: fiber of prestacks}. Combined with Lemma~\ref{lem: trivializations of 2-gerbe torsor for gerbes}, we obtain a sequence of equivalences of bicategories
    \begin{align}
        \Bun_\cG(X)_P \simeq {\sf Triv}(f, \rho^*\alpha) \simeq \Gerbe_A(X),
    \end{align}
    where the first equivalence is canonical and the second depends upon the choice of a trivialization $(g \colon Z \to Y, \gamma)$ of $(f, \rho^*\alpha)$, or equivalently, a lift $(g \circ f, g^*\rho, \gamma)$ of $(f, \rho)$. 
\end{proof}

\begin{rmk}
We can describe the functor from $\Bun_\cG(X)_P$ to ${\sf Triv}(f, \rho^*\alpha)$ more concretely. An object of $\Bun_\cG(X)_P$ is a pair, consisting of an object $(f_1, \rho_1, \gamma_1)$ in $\Bun_\cG(X)$ and an isomorphism between $(f_1, \rho_1)$ and $(f, \rho)$ in $\Bun_G(X)$. Similarly, a 1-morphism is a pair consisting of a 1-morphism between the objects of $\Bun_\cG(X)$, and a 2-morphism providing compatibility data on the level of 1-morphisms in $\Bun_G(X)$. Finally, a 2-morphism is a 2-morphism in $\Bun_\cG(X)$ which satisfies a compatibility condition in $\Bun_G(X)$. 

We claim that every object is isomorphic to one consisting of a pair where the object of $\Bun_\cG(X)$ is of the form $(f \circ g, g^*\rho, \gamma_1)$ for some essential equivalence $g \colon Z \to Y$, and the 1-morphism in $\Bun_G(X)$ corresponds to $\id \colon g^* \rho \to g^*\rho$ in $\Bun_G(Z)$. To prove this claim, we first remark that every object of  $\Bun_\cG(X)_P$ is isomorphic to a pair where the $\cG$-bundle is of the form $(f_1 = f \circ g, \rho_1, \gamma_1)$ for some essential equivalence $g \colon Z \to Y$, and the 1-morphism in $\Bun_G(X)$ is given by a 1-morphism $h \colon \rho_1 \to g^* \rho$ in $\Bun_G(Z)$. Furthermore, it follows from Lemma \ref{lem: lifting 1-morphisms} that up to canonical isomorphism we can also replace the $\cG$-bundle by one of the form $(f_1 = f \circ g, g^* \rho, \gamma)$, and the 1-morphism by the identity as claimed. This object of $\Bun_\cG(X)_P$ corresponds to a trivialization $(g, \gamma)$ of $\rho^*\alpha$. 
\end{rmk}

Given a principal $\cG$-bundle represented by data $(f \colon Y \to X, (\rho, \gamma_1)_Y)$ and an $A$-gerbe $(f \colon Y \to X, \gamma_2)$, we can ``twist'' the principal bundle by the gerbe to obtain a new principal bundle, 
\begin{align}\label{eq: twist}
    (f \colon Y \to X, (\rho, \gamma_1 \cdot \gamma_2)_Y);
\end{align}
indeed, we have 
\begin{align*}
    d(\gamma_1 \gamma_2) = (d\gamma_1)(d\gamma_2) = (\rho^*\alpha)\cdot 1_G = \rho^*\alpha. 
\end{align*}
We remark that this twisting does not change the underlying $G$-bundle, which is still determined by $\rho$. Furthermore, the twisting can be defined even when the original $\cG$-bundle and the gerbe are defined relative to two different surjective submersions $f_1 \colon Y_1 \to X$ and $f_2 \colon Y_2 \to X$: it suffices to pull the data back to $Y_1 \times_X Y_2$ before multiplying the 2-cochains as in \eqref{eq: twist}.

\begin{prop}\label{prop:twist}
    The twisting operation \eqref{eq: twist} extends to a functor
    \begin{align}
        F \colon \Bun_\cG(X) \times \Gerbe_A(X) \to \Bun_\cG(X).
    \end{align}
    The functor $\pi \colon \Bun_\cG \to \Bun_G$ is invariant under $F$
    in the sense that $\pi \circ F \cong \pi$. 
\end{prop}

\begin{proof}
    We will describe a functor
    \begin{align}
        F^\pre \colon \Bun_\cG^\pre(X) \times \Gerbe_A^\pre(X) \to \Bun_\cG^\pre(X)
    \end{align}
    which will induce the desired functor $F$ upon stackification. 

    We know that on objects $F^\pre$ must behave as in the formula \eqref{eq: twist}: 
    \begin{align}
        F^\pre((\rho, \gamma_1), \gamma_2) = (\rho, \gamma_1\cdot \gamma_2).
    \end{align}
    Similarly, given 1-morphisms $(h, \eta_1) \colon (\rho_1, \gamma_1) \to (\rho_2, \gamma_2)$ in $\Bun_\cG^\pre(X)$ and $\eta_2 \colon \gamma_3 \to \gamma_4$ in $\Gerbe_A^\pre(X)$ it is easy to see that setting 
    \begin{align*}
        F^\pre((h, \eta_1), \eta_2) = (h, \eta_1 \cdot \eta_2)
    \end{align*} 
    gives a 1-morphism from $(\rho_1, \gamma_1 \cdot \gamma_3)$ to $(\rho_2, \gamma_2 \cdot \gamma_4)$. The value of $F^\pre$ on 2-morphisms is similarly determined by pointwise multiplication in $A$. It is straightforward to check that this defines a functor, and that for $\pi^\pre \colon \Bun_\cG^\pre(X) \to \Bun_G^\pre(X)$ we have $\pi^\pre \circ F^\pre = \pi^\pre$.

    The desired functor $F$ and invariance $\pi \circ F \cong \pi$ follow upon stackification. 
\end{proof}

In fact, this functor $F$ can be extended to a (weak) action of the symmetric monoidal bicategory $\Gerbe_A(X)$ on the bicategory $\Bun_\cG(X)$ making it into a torsor over $\Bun_G(X)$. Spelling this out in detail requires introducing natural transformations and further higher compatibility data, and is beyond the scope of this paper. However, what is relevant for our purposes is that this higher categorical action is free and transitive on fibers, inducing an equivalence of $\Gerbe_A(X)$ with each non-empty fiber. More precisely, we have the following: 

\begin{thm}
For any object $\cP$ of $\Bun_\cG(X)$ with underlying $G$-bundle $P$, the restriction of $F$ to $\Gerbe_A(X) \cong \{\cP\} \times \Gerbe_A(X)$ provides the equivalence of $\Gerbe_A(X)$ with the bicategory $\Bun_\cG(X)_P$ of Theorem \ref{thm:2gerbetriv}.
\end{thm}

\begin{proof}
    Choose cocycle data $(Y\to X, (\rho, \gamma)_Y)$ classifying $\cP$. Restricting the functor $F^\pre$ from the proof of Proposition \ref{prop:twist}, we obtain a functor
    \begin{align}
        \{(\rho, \gamma)\}\times \Gerbe_A^\pre(Y) \to \Bun_G^\pre(Y),
    \end{align}
    that induces the functor of interest upon stackification. Comparing the formulas given in the proofs of Lemma~\ref{lem: trivializations of 2-gerbe torsor for gerbes} and Proposition~\ref{prop:twist}, one recovers the functor of Theorem~\ref{thm:2gerbetriv}. 
\end{proof}

\section{Flat string structures}\label{sec:flatstringserction}

In this section we apply Theorem~\ref{thm:2gerbetriv} to analyze flat string structures in the sense of Definition~\ref{defn:string}. The material in Sections~\ref{sec: flat vect bun}-\ref{sec: flat spin} is well-known, though the presentation below is tailored to fit our purposes. We begin with a recollection of flat bundles in Section~\ref{sec: flat vect bun}, before translating flat orientations and spin structures into this language in Sections~\ref{sec: flat ori}-\ref{sec: flat spin}. Section~\ref{sec: flat string} then uses the language of 2-groups to give a similar formulation of flat string structures.

\subsection{Background: Flat vector bundles}\label{sec: flat vect bun}

Let $M$ be a manifold. A metrized vector bundle $V\to M$ with a compatible connection $\nabla$ is flat if its curvature $F=\nabla\circ \nabla\in \Omega^2(M;{\rm End}(V))$ vanishes. This can be rephrased in terms of principal bundles: a bundle is flat if and only if its $\rmO_n$-frame bundle admits a reduction of structure group to an $\rmO_n^\delta$-bundle, where $\rmO_n^\delta$ is the orthogonal group $\rmO_n$ endowed with the discrete topology. By standard principal bundle theory, such a reduction of structure group is equivalent to admitting a \v{C}ech cocycle whose transition data uses only \emph{constant} $\rmO_n$-valued functions. This generalizes to flat bundles on Lie groupoids as follows; compare Example~\ref{ex:standardOn} below. 

\begin{defn}\label{eq:defnflatvb}
    The groupoid of \emph{vector bundles} on a Lie groupoid $X$ is the groupoid of maps $X\to \pt\sq \rmO_n$ in $\LieGrpd[W^{-1}]$. A \emph{flat 
    structure} on a vector bundle is a lift in $\LieGrpd[W^{-1}]$, 
\beq
\begin{tikzpicture}[baseline=(basepoint)];
\node (B) at (4,0) {$\pt\sq \rmO_n^\delta$};
\node (D) at (4,-1) {$\pt\sq \rmO_n.$};
\node (C) at (0,-1) {$X$};
\draw[->,dashed] (C) to (B);
\draw[->] (B) to (D);
\draw[->] (C) to node [below] {$V$} (D);
\path (0,-.75) coordinate (basepoint);
\end{tikzpicture}\label{eq:flatbundle}
\eeq
The groupoid of \emph{flat vector bundles on $X$} is the groupoid of maps $X\to \pt\sq \rmO_n^\delta$  in $\LieGrpd[W^{-1}]$, which we identify with the (1-)groupoid $\Bun_{\rmO_n^\delta}(X)$ from the previous section.
\end{defn}

\begin{ex}\label{ex:standardOn}
When $X=M$ is a manifold, the dashed arrow in~\eqref{eq:flatbundle} is a zig-zag 
\beq\nonumber
M\xleftarrow[]{\sim}\check{C}(Y)\xrightarrow{V} \pt\sq \rmO_n^\delta
\eeq
for an essential equivalence determined by a surjective submersion $Y\to M$. Hence, \eqref{eq:flatbundle} recovers a standard definition of flat bundles in terms of locally constant transition functions; see Definition~\ref{defn:prinbund0}. 
\end{ex}

\begin{ex}\label{ex:flatmfld}
A \emph{flat manifold} is a Riemannian manifold $M$ whose curvature vanishes, i.e., such that the tangent bundle $TM$ is equipped with a flat structure. More generally, for a Lie group $G$ acting on a manifold $M$, a \emph{$G$-equivariant flat structure on $M$} is a lift
\beq
\begin{tikzpicture}[baseline=(basepoint)];
\node (B) at (4,0) {$\pt\sq \rmO_n^\delta$};
\node (D) at (4,-1) {$\pt\sq \rmO_n$};
\node (C) at (0,-1) {$M\sq G$};
\draw[->,dashed] (C) to (B);
\draw[->] (B) to (D);
\draw[->] (C) to node [below] {$TM$} (D);
\path (0,-.75) coordinate (basepoint);
\end{tikzpicture}\nonumber
\eeq
where $TM$ classifies the $G$-equivariant frame bundle of $M$, compare Example~\ref{ex:Hequivprin}. 
\end{ex}

\subsection{Orientations of flat bundles}\label{sec: flat ori}

The inclusion $\SO_n^\delta\hookrightarrow \rmO_n^\delta$ determines a functor
\beq\label{eq:orientationlift}
\Bun_{\SO_n^\delta}(X)\to \Bun_{\rmO_n^\delta}(X)
\eeq
for any Lie groupoid $X$. 

\begin{defn}
Given a flat vector bundle $V\in \Bun_{\rmO_n^\delta}(X)$ the category of \emph{orientations} on~$V$ is the homotopy fiber of~\eqref{eq:orientationlift}. 
\end{defn}

Unwinding the definition, an orientation on $V$ is 2-commutative triangle in $\LieGrpd[W^{-1}]$
\beq
\begin{tikzpicture}[baseline=(basepoint)];
\node (B) at (4,0) {$\pt\sq \SO_n^\delta$};
\node (D) at (4,-1) {$\pt\sq \rmO_n^\delta$};
\node (C) at (0,-1) {$X$};
\draw[->,dashed] (C) to (B);
\draw[->] (B) to (D);
\draw[->] (C) to node [below] {$V$} (D);
\path (0,-.75) coordinate (basepoint);
\end{tikzpicture}\label{eq:orlift2}
\eeq
where we emphasize that 2-commutativity is additional data beyond the dashed arrow. Similarly, an isomorphism between oriented vector bundles is a 2-morphism between dashed arrows where the 2-commutative data is required to satisfy the evident compatibility condition. 

\begin{ex}\label{ex:mfldorienation}
When $X=M$ is a manifold, the dashed arrow in~\eqref{eq:orlift2} specifies a \v Cech groupoid $\check{C}(Y)$ and a $\SO_n^\delta$-valued cocycle for a principal bundle defined relative to the surjective submersion $Y\to M$. The 2-commutativity data in~\eqref{eq:orlift2} gives an isomorphism (possibly after refining the submersion) with the $\rmO_n^\delta$-valued cocycle specifying the given flat vector bundle. 
\end{ex}

\begin{defn}\label{defn:orientationline}
The \emph{orientation double cover} associated to a flat vector bundle is the $\{\pm 1\}$-principal bundle on~$X$, denoted $\Or_V\to X$, classified by the composition in $\LieGrpd[W^{-1}]$, 
\beq\label{eq:classifyor}
X \xrightarrow{V} \pt\sq \rmO_n^\delta\xrightarrow{\det} \pt\sq \{\pm 1\}
\eeq
where the rightmost arrow is induced by the determinant homomorphism. 
\end{defn}

\begin{rmk}
Explicitly, the composition~\eqref{eq:classifyor} gives an essential equivalence $Y\xrightarrow[]{\sim}X$ and a functor $Y\to \pt\sq \{\pm 1\}$ determined by a map $Y_1\to \{\pm 1\}$ that reads off whether the cocycle for $V$ preserves or reverses orientation. 
\end{rmk}

\begin{defn}
    A \emph{trivialization} of the orientation double cover is a section $X\to \Or_V$. Equivalently, a trivialization is a 2-commuting diagram in $\LieGrpd[W^{-1}]$,
\beq
\begin{tikzpicture}[baseline=(basepoint)];
\node (B) at (6,0) {$\pt$};
\node (D) at (4,-1) {$\pt\sq \rmO_n^\delta$};
\node (A) at (6,-1) {$\pt\sq\{\pm 1\}.$};
\node (C) at (0,-1) {$X$};
\draw[->,dashed] (C) to (B);
\draw[->] (B) to (A);
\draw[->] (D) to node [below] {$\det$} (A);
\draw[->] (C) to node [below] {$V$} (D);
\path (0,-.75) coordinate (basepoint);
\end{tikzpicture}\label{eq:ortriv}
\eeq
The \emph{category of vector bundles on $X$ with trivialized orientation double cover} has objects diagrams~\eqref{eq:ortriv} and morphisms isomorphisms between maps $X \to \pt\sq \rmO_n^\delta$ that are compatible with the trivialization of $\Or_V$.
\end{defn}

\begin{lem} \label{lem:orientations}
The category $\Bun_{\SO_n^\delta}(X)$ of oriented flat vector bundles is equivalent to the category of flat vector bundles with trivialized orientation double cover.  
\end{lem}

\bp
Since $\pt \sq \SO_n^\delta \to \pt \sq \rmO_n^\delta$ is the pullback of $\pt \to \pt \sq\{ \pm 1\}$ along the morphism $\det \colon \pt \sq \rmO_n^\delta \to \pt \sq \{\pm 1\}$, the orientation double cover fits into a pullback diagram as follows:
 \beq\begin{tikzpicture}[baseline=(basepoint)];
 \node (A) at (0,0) {$\Or_V$};
 \node (B) at (3,0) {$\pt \sq \SO_n^\delta$};
 \node (C) at (6,0) {$\pt$};
 \node (P) at (.7,-.4) {\scalebox{1.5}{$\lrcorner$}};
 \node (P) at (3.7,-.4) {\scalebox{1.5}{$\lrcorner$}};
 \node (AA) at (0,-1) {$X$};
 \node (BB) at (3,-1) {$\pt\sq \rmO_n^\delta$};
 \node (CC) at (6,-1) {$\pt\sq \{\pm 1\}.$};
 \draw[->] (A) to (B);
 \draw[->] (B) to (C);
 \draw[->] (AA) to (BB);
 \draw[->] (BB) to (CC);
 \draw[->] (A) to (AA);
 \draw[->] (B) to (BB);
 \draw[->] (C) to (CC);
 \path (0,-.5) coordinate (basepoint);
 \end{tikzpicture}\label{Eq:doublepullback}
 \eeq
The statement of the lemma then follows from the universal property of the pullback
\beq
\begin{tikzpicture}[baseline=(basepoint)];
\node (AA) at (-1,1) {$X$};
\node (A) at (0,0) {$\Or_V$};
\node (B) at (4,0) {$\pt\sq \SO_n^\delta$};
\node (D) at (4,-1) {$\pt\sq \rmO_n$};
\node (C) at (0,-1) {$X$};
\draw[->,dashed, bend left=10] (AA) to (B);
\draw[->,bend right] (AA) to node [left] {$\id$} (C);
\draw[->,dotted] (AA) to (A);
\draw[->] (A) to (B);
\draw[->] (A) to (C);
\draw[->,dashed] (C) to (B);
\draw[->] (B) to (D);
\draw[->] (C) to node [below] {$V$} (D);
\path (0,-.75) coordinate (basepoint);
\end{tikzpicture}\label{eq:univprop1}
\eeq
where the dashed arrows are the data of the oriented vector bundle $X\to \pt\sq \rmO_n^\delta$, and the dotted arrow is the data of a trivializing section of $\Or_V$. Hence, a lift determines a trivializing section $X\to \Or_V$ and conversely. Similarly, an isomorphism between lifts determines an isomorphism between trivializations of $\Or_V$. 
\ep

\begin{rmk}
We observe that the smooth functor $\pt\sq \SO_n^\delta\to \pt\sq \rmO_n^\delta$ is not a fibration of groupoids. However, it factors as a composition of smooth functors
\beq\label{eq:factorizationofor}
\pt\sq \SO_n^\delta  \xrightarrow{\sim} (\pt\coprod \pt)\sq \rmO_n^\delta \rightarrow \pt\sq \rmO_n^\delta
\eeq
where $(\pt\coprod \pt)\sq \rmO_n^\delta$ is the action groupoid for the $\rmO_n^\delta$-action on $\pt\coprod \pt\simeq \{\pm 1\}$ through the determinant homomorphism. In~\eqref{eq:factorizationofor}, the first arrow is an essential equivalence and the second arrow is a fibration of groupoids. Hence, replacing $\pt\sq \SO_n^\delta$ in~\eqref{Eq:doublepullback} with $(\pt\coprod \pt)\sq \rmO_n^\delta$, we find that the pullback exists in Lie groupoids and is the claimed orientation double cover. 
\end{rmk}

\begin{rmk}\label{rmk:orientationsareset}
In light of Lemma~\ref{lem:orientations}, an orientation of a vector bundle $V$ is a trivialization of a $\{\pm 1\}$-bundle on $X$. The category of trivializations of a $\{\pm 1\}$-bundle is equivalent to a set.
Hence orientations of a flat vector bundle $V$ naturally form a set. 
\end{rmk}

\begin{ex}\label{ex:Gequivorientation}
For a flat $G$-manifold $M$, a \emph{$G$-equivariant orientation} is a lift
\beq
\begin{tikzpicture}[baseline=(basepoint)];
\node (B) at (4,0) {$\pt\sq \SO_n^\delta$};
\node (D) at (4,-1) {$\pt\sq \rmO_n^\delta$};
\node (C) at (0,-1) {$M\sq G$};
\draw[->,dashed] (C) to (B);
\draw[->] (B) to (D);
\draw[->] (C) to node [below] {$TM$} (D);
\path (0,-.75) coordinate (basepoint);
\end{tikzpicture}\nonumber
\eeq
where $TM$ classifies the $G$-equivariant frame bundle of $M$. By Lemma~\ref{lem:orientations}, an equivariant orientation is equivalent to a $G$-invariant trivialization of the orientation double cover of~$TM$, i.e., a unit norm $G$-invariant section of $\Lambda^{\rm top}TM$ over $M$. 
\end{ex}

\subsection{Flat spin structures}\label{sec: flat spin}
The double cover $\Spin_n^\delta\to \SO_n^\delta$ postcomposed with the inclusion $\SO_n^\delta\hookrightarrow \rmO_n^\delta$ determines a functor
\beq\label{eq:spinlift}
\Bun_{\Spin_n^\delta}(X)\to \Bun_{O_n^\delta}(X)
\eeq
for any Lie groupoid $X$. 

\begin{defn}
Given a flat vector bundle $V\in \Bun_{\rmO_n^\delta}(X)$ the category of \emph{flat spin structures} on~$V$ is the homotopy fiber of~\eqref{eq:spinlift}. 
\end{defn}

Hence, a spin structure on $V$ is 2-commutative triangle in $\LieGrpd[W^{-1}]$
\beq
\begin{tikzpicture}[baseline=(basepoint)];
\node (B) at (4,0) {$\pt\sq \Spin_n^\delta$};
\node (D) at (4,-1) {$\pt\sq \rmO_n^\delta.$};
\node (C) at (0,-1) {$X$};
\draw[->,dashed] (C) to (B);
\draw[->] (B) to (D);
\draw[->] (C) to node [below] {$V$} (D);
\path (0,-.75) coordinate (basepoint);
\end{tikzpicture}\label{eq:spinlift2}
\eeq
The functor~\eqref{eq:spinlift} factors as
$$
\Bun_{\Spin_n^\delta}(X)\to \Bun_{\SO_n^\delta}(X)\to \Bun_{\rmO_n^\delta}(X).
$$
Hence, there is a natural functor from the category of vector bundles with spin structure to the category of oriented vector bundles.

\begin{ex}
In analogy to Example~\ref{ex:mfldorienation}, when $X=M$ is a manifold the dashed arrow in~\eqref{eq:spinlift2} specifies a \v Cech cocycle for a $\Spin_n^\delta$-principal bundle that (up to refinement) lifts the given cocycle for the flat vector bundle $V$. 
\end{ex}

\begin{rmk}
Flat spin structures have been studied in the context of index theory, e.g., see~\cite{Dekimpe,LUTOWSKI,PFAFFLE}.
\end{rmk}

\begin{defn}\label{defn:SW}
The \emph{Stiefel--Whitney gerbe of a flat vector bundle~$V$}, denoted $\SW_V$, is the pullback in $\LieGrpd[W^{-1}]$
\beq
\begin{tikzpicture}[baseline=(basepoint)];
\node (A) at (0,0) {$\SW_V$};
\node (B) at (3,0) {$\pt\sq \Spin_n^\delta$};
\node (C) at (0,-1.5) {$\Or_V$};
\node (D) at (3,-1.5) {$\pt\sq \SO_n^\delta$};
\node (P) at (.7,-.4) {\scalebox{1.5}{$\lrcorner$}};
\draw[->] (A) to  (B);
\draw[->] (B) to  (D);
\draw[->] (A) to  (C);
\draw[->] (C) to (D);
\path (0,-.75) coordinate (basepoint);
\end{tikzpicture}\label{eq:SWgerbe}
\eeq
of the fibration of groupoids $\pt\sq \Spin_n^\delta\to \pt\sq \SO_n^\delta$. This defines a Lie groupoid over~$X$ via the composition
\beq\label{eq:SWoverX}
\SW_V\to \Or_V\to X.
\eeq
A \emph{trivialization} of $\SW_V\to X$ is a section $X\to \SW_V$ (considered in $\LieGrpd[W^{-1}]$).  Trivializations are naturally the objects in a category whose morphisms are isomorphisms between sections. 
\end{defn}

\begin{rmk}
A trivializing section $X\to \SW_V$ induces a trivialization of $\Or_V$. We note that one can ask for an intermediate trivialization of $\SW_V\to \Or_V$ over the double cover $\Or_V$. Below, we will only consider trivializations of the composite~\eqref{eq:SWoverX}. \end{rmk}

\begin{rmk}
    The equivalence between trivializations of gerbes and sections in $\LieGrpd[W^{-1}]$ is standard; e.g., see~\cite[4.2]{BehrendXu}. 
\end{rmk}

\begin{lem}
The Stiefel--Whitney gerbe is a $\Z/2$-gerbe on the Lie groupoid $\Or_V$, i.e., a $(\pt\sq \Z/2)$-bundle determined by a 2-functor $\SW_V\colon \Or_V\to \pt\sq (\pt\sq (\Z/2))$. 
\end{lem}

\bp
We describe the pullback~\eqref{eq:SWgerbe} explicitly, using that the lax pullback along a fibration of (Lie) groupoids $\pt\sq \Spin_n^\delta\to \pt\sq \SO_n^\delta$ is equivalent to the strict pullback. Suppose that the map $\Or_V\to \pt\sq \SO_n^\delta$ in $\LieGrpd[W^{-1}]$ is determined by the zig-zag of smooth functors $\Or_V\xleftarrow[]{\sim} Y \to \pt\sq \SO_n^\delta$. Then $\SW_V$ has $Y_0$ as objects, and morphisms are the pullback (in smooth manifolds)
\beq
\begin{tikzpicture}[baseline=(basepoint)];
\node (A) at (0,0) {$\widehat{Y}_1$};
\node (B) at (3,0) {$\Spin_n^\delta$};
\node (C) at (0,-1.5) {$Y_1$};
\node (D) at (3,-1.5) {$\SO_n^\delta$};
\node (P) at (.7,-.4) {\scalebox{1.5}{$\lrcorner$}};
\draw[->] (A) to  (B);
\draw[->] (B) to  (D);
\draw[->] (A) to  (C);
\draw[->] (C) to (D);
\path (0,-.75) coordinate (basepoint);
\end{tikzpicture}\label{eq:SWgerbe2}
\eeq
where $\widehat{Y}_1\to Y_1$ is a principal $\Z/2$-bundle.

But since the topology on $\SO_n^\delta$ is discrete, the map $Y_1\to \SO_n^\delta$ is locally constant, and hence $\widehat{Y}_1\to Y_1$ is the trivial $\Z/2$-principal bundle. This implies that data of $\SW_V$ is the essential equivalence $Y\xrightarrow[]{\sim} \Or_V$ and a $\Z/2$-valued 2-cocycle $Y_1\times_{Y_0} Y_1\to \Z/2$. This proves the lemma. 
\ep

\begin{rmk}
When $X=M$ is a manifold, the 2-cocycle determining $\SW_V$ is a standard cocycle construction of the second Stiefel--Whitney class $w_2(V)$. More generally, the Stiefel--Whitney gerbe sits in a fiber sequence 
$$
\pt\sq \Z/2\to \SW_V\to \Or_V
$$
that is a \emph{$\Z/2$-central extension} in the sense of \cite[\S4-5]{BehrendXu}. Such central extensions are examples of twists for Thom classes in real K-theory, see~\cite[1.92 and 1.161]{vienna}. 
\end{rmk}

\begin{prop}\label{Prop:spinstrucequiv}
The category of flat spin structures on $V$ is equivalent to the category of trivializations of $SW_V\to X$. 
\end{prop}

\bp
Consider the diagram in $\LieGrpd[W^{-1}]$ of iterated pullbacks, 
\beq
\begin{tikzpicture}[baseline=(basepoint)];
\node (AA) at (-1,1) {$X$};
\node (A) at (0,0) {$\SW_V$};
\node (B) at (4,0) {$\pt\sq \Spin_n^\delta$};
\node (D) at (4,-1) {$\pt\sq \SO_n$};
\node (C) at (0,-1) {$\Or_V$};
\node (E) at (4,-2) {$\pt\sq \rmO_n$};
\node (F) at (0,-2) {$X$};
\draw[->] (C) to (F);
\draw[->] (D) to (E);
\draw[->,dashed, bend left=10] (AA) to (B);
\draw[->,bend right] (AA) to node [left] {$\id$} (F);
\draw[->,dotted] (AA) to (A);
\draw[->] (A) to (B);
\draw[->] (A) to (C);
\draw[->,dashed] (F) to (B);
\draw[->] (B) to (D);
\draw[->] (C) to (D);
\draw[->] (F) to node [below] {$V$} (E);
\path (0,-.75) coordinate (basepoint);
\end{tikzpicture}\label{eq:univprop2}
\eeq
where the dashed arrows are the data of a spin structure on $V$, and the dotted arrow is the data of a trivializing section. By the universal property of the pullback, these are equivalent data. \ep

\begin{ex}\label{ex:Gequivspin}
For a flat $G$-manifold $M$, a \emph{$G$-equivariant flat spin structure} is a lift
\beq
\begin{tikzpicture}[baseline=(basepoint)];
\node (B) at (4,0) {$\pt\sq \Spin_n^\delta$};
\node (D) at (4,-1) {$\pt\sq \rmO_n^\delta$};
\node (C) at (0,-1) {$M\sq G$};
\draw[->,dashed] (C) to (B);
\draw[->] (B) to (D);
\draw[->] (C) to node [below] {$TM$} (D);
\path (0,-.75) coordinate (basepoint);
\end{tikzpicture}\label{eq:Gequivspin}
\eeq
where $TM$ classifies the $G$-equivariant frame bundle of $M$. From Example~\ref{ex:Gequivorientation}, this lift determines a $G$-invariant orientation on~$M$. By inspection of~\eqref{eq:Gequivspin}, an equivariant spin structure is the data of a $\Spin_n^\delta$-principal bundle ~$P\to M$ covering the frame bundle of $TM$ where $P$ is furthermore equipped with a $G$-action and the map $P\to M$ is $G$-equivariant. An equivariant spin structure is the data of an orientation and a trivialization of the Stiefel--Whitney gerbe of $TM$ (e.g., see~\cite[Example 4.1.4]{WaldorfAGT}), together with equivariant descent data for these trivializations along the map $M\to M\sq G$ in $\LieGrpd[W^{-1}]$. 
\end{ex}

\subsection{Flat string structures}\label{sec: flat string}

The 1-homomorphism of 2-groups $\String_n^\delta\to \Spin_n^\delta\to \rmO_n^\delta$ determines a 2-functor
\beq\label{eq:stringlift}
\Bun_{\String_n^\delta}(X)\to \Bun_{\rmO_n^\delta}(X)
\eeq
for any Lie groupoid $X$. 

\begin{defn}\label{defn:stringstruc}
Given a flat vector bundle $V\in \Bun_{\rmO_n^\delta}(X)$ the bicategory of \emph{flat string structures} on~$V$ is the homotopy fiber of~\eqref{eq:stringlift}. 
\end{defn}

Hence, a string structure on $V$ is lifting 2-functor
\beq
\begin{tikzpicture}[baseline=(basepoint)];
\node (B) at (4,0) {$\pt\sq \String_n^\delta$};
\node (D) at (4,-1) {$\pt\sq \rmO_n^\delta.$};
\node (C) at (0,-1) {$X$};
\draw[->,dashed] (C) to (B);
\draw[->] (B) to (D);
\draw[->] (C) to node [below] {$V$} (D);
\path (0,-.75) coordinate (basepoint);
\end{tikzpicture}\nonumber
\eeq
Similarly to above, a string structure on $V$ determines a spin structure and an orientation. 

We recall that there is a 3-cocycle $\Spin_n^\delta\times \Spin_n^\delta\times \Spin_n^\delta\to U(1)$ that classifies the discrete string 2-group~\eqref{diag:discretestring}; we denote this 3-cocycle by $\CS^\delta$. By Proposition~\ref{prop:2grpconstr}, $\CS^\delta$ is unique up to coboundary. 

\begin{defn}\label{defn:CS2grb}
The (flat) \emph{Chern--Simons 2-gerbe} on $\pt\sq \Spin_n^\delta$ is specified by the 3-cocycle $\CS^\delta$ on $\Spin_n^\delta$. \end{defn}

\begin{rmk}
We observe that (by construction) the fractional Pontryagin class maps to the class underlying the 3-cocycle $\CS^\delta$, 
\begin{align*}
\Big[\frac{p_1}{2}\Big]\in \rmH^4(B\Spin_n;\Z)\simeq \rmH^3_{\rm SM}(\Spin_n;U(1))\xrightarrow{\delta^*} &\rmH^3_{\rm SM}(\Spin_n^\delta;U(1))\\
&\simeq \rmH^3(B\Spin_n^\delta;U(1))\ni [\CS^\delta].
\end{align*}
Hence, the Chern--Simons 2-gerbe in sense of Waldorff~\cite[\S2.1]{Waldorfstring} pulls back to the flat Chern--Simons 2-gerbe of Definition~\ref{defn:CS2grb}.
\end{rmk}

\begin{defn}\label{defn:CS}
The \emph{Chern--Simons 2-gerbe of a flat vector bundle~$V$}, denoted $\CS_V$, is the 2-gerbe over the Lie groupoid $\SW_V$ gotten by pulling back the flat Chern--Simons 2-gerbe along the top horizontal arrow $\SW_V\to \pi\sq\Spin_n^\delta$ in~\eqref{eq:SWgerbe} (see Example~\ref{ex:2gerbefromP} for a formula for this pullback). A \emph{trivialization} of $\CS_V\to X$ is a trivializing section $\sigma\colon X\to \SW_V$ and a trivialization of the 2-gerbe $\sigma^* \CS_V$ over $X$. We observe that trivializations of the Chern--Simons 2-gerbe form a (possibly empty) bicategory. 
\end{defn}

\begin{rmk} 
One can rephrase the above in the geometry of 2-stacks. The full 3-category of 2-stacks is understood, but technical, e.g., see~\cite{Nuiten}. Here we will only indicate the relevance for (geometric) string structures. 
Consider the iterated pullback diagram in 2-stacks 
\beq
\begin{tikzcd}
    \CS_V\arrow[r]\arrow[d]\arrow[dr, phantom, "\scalebox{1}{$\lrcorner$}" , pos=0.1] & \pt\sq\String_n^\delta \arrow[d]\\
    \SW_V\arrow[r]\arrow[d]\arrow[dr, phantom, "\scalebox{1}{$\lrcorner$}" , pos=0.1] & \pt\sq\Spin_n^\delta\arrow[d]\\
    \Or_V\arrow[r]\arrow[d]\arrow[dr, phantom, "\scalebox{1}{$\lrcorner$}" , pos=0.1]  & \pt\sq \SO_n^\delta\arrow[d]\\
    X \arrow[uuu,  bend left=50, dashed]\arrow[r, "V"']& \pt\sq\rmO_n^\delta
\end{tikzcd}
\eeq
where dashed arrow is a trivializing section of $\CS_V\to X$. By the universal property of the pullback, such a trivialization is equivalent to a map of 1-stacks $X\to \SW_V$ together with a trivialization of a 2-gerbe, as in Definition~\ref{defn:CS}. The advantage of Definition \ref{defn:CS} is that we need only work with the bicategory of trivializations of a 2-gerbe, avoiding 3-categories.
\end{rmk}

\begin{thm}\label{thm:flatstringdata}
Given a flat vector bundle $V\to X$, the bicategory of flat string structures on $V$ is equivalent to the bicategory of trivializations of $\CS_V\to X$. 
\end{thm}

\bp 
This follows from Proposition~\ref{Prop:spinstrucequiv} and Theorem~\ref{thm:2gerbetriv}. 
\ep

\begin{rmk} 
We note that there are two possible versions of string structure on a flat spin vector bundle. One involves (potentially non-flat) trivializations of the Chern--Simons 2-gerbe, whereas the second only considers flat trivializations. In terms of the extensions~\eqref{diag:discretestring} and the associated 2-group principal bundles, this distinction corresponds to whether one takes the usual topology on $U(1)$ as a Lie group, or its discrete topology. 
\end{rmk}

\begin{ex}
For a flat $G$-manifold $M$, a \emph{$G$-equivariant flat string structure} is a lift
\beq
\begin{tikzpicture}[baseline=(basepoint)];
\node (B) at (4,0) {$\pt\sq \String_n^\delta$};
\node (D) at (4,-1) {$\pt\sq \rmO_n^\delta$};
\node (C) at (0,-1) {$M\sq G$};
\draw[->,dashed] (C) to (B);
\draw[->] (B) to (D);
\draw[->] (C) to node [below] {$TM$} (D);
\path (0,-.75) coordinate (basepoint);
\end{tikzpicture}\label{eq:Gequivstring}
\eeq
where $TM$ classifies the $G$-equivariant frame bundle of $M$. 
By inspection of~\eqref{eq:Gequivstring}, an equivariant string structure is a $G$-equivariant $\String_n^\delta$-bundle on $M$ whose underlying $\Spin_n^\delta$-bundle is the spin bundle of $M$ for the underlying flat spin structure (see Example~\ref{ex:Gequivspin}). An equivariant string structure is therefore the data of an orientation, spin structure, and trivialization of the (nonequivariant) Chern--Simons 2-gerbe of $TM$ (see~\cite[\S5.1]{WaldorfAGT} or~\cite[\S2.1]{Waldorfstring}), together with equivariant descent data. An explicit (though tedious) exercise using Theorem~\ref{thm:2gerbetriv} gives an explicit description of this descent data; one recovers a flat version of Sharpe's heterotic $B$-field from~\cite[Equations~4 and~5]{sharpetorsionhet}. 
\end{ex}

\appendix 

\section{The discrete string 2-group}\label{sec:disstring}
The forgetful functor from geometric stacks to (discete) groupoids provides a discrete 2-group underlying any smooth 2-group~\cite[page 28]{SP11}. Below we give an explicit description of the discrete 2-group underlying the string 2-group, $\String_n$, i.e., the categorical central extension~\eqref{eq:stringgroupdef} 
 \cite[Theorem~2]{SP11}.

For a Lie group G and abelian Lie group $A$, smooth categorical central extensions
$$
1\to \pt\sq A\to \cG\to G\to 1
$$
are classified (up to 1-isomorphism) by $\rmH_{\rm SM}^3(G;A)$, the 3rd Segal--Mitchison cohomology group of $G$ with values in $A$~\cite[Theorem 1]{SP11}. This cohomology is computed by a double complex that depends on a simplicial cover of the nerve~$BG_\bullet$. 

\begin{lem} \label{lem:SMtoGRP}
When $G$ has the discrete topology, the Segal--Mitchison cohomology $\rmH^\bullet_{\rm SM}(G;A)$ is isomorphic to the usual group cohomology of $G$ valued in $A$. 
\end{lem}
\bp
In fact, the Segal--Mitchison bicomplex of $G$ with coefficients in $A$ is quasi-isomorphic to the complex computing the usual group cohomology of $G$. To construct this quasi-isomorphism, take the open cover of $G$ given by its points; this determines a simplicial cover of $B G_\bullet$. The resulting double complex is exact for the \v Cech differential, where the remaining differential is the one from group cohomology (e.g., see \cite[page~35]{SP11}). Hence we obtain a canonical quasi-isomorphism between the Segal--Mitchison complex and and the usual complex computing the group cohomology of~$G$. 
\ep

The categorical central extensions (and the Segal--Mitchison complex) are functorial in the groups $G$ and $A$~\cite[page~4]{SP11}. This provides a map from a discrete extension to the string extension, 
\beq\begin{tikzpicture}[baseline=(basepoint)];
\node (A) at (0,0) {$1$};
\node (B) at (2,0) {$\pt\sq U(1)$};
\node (C) at (4,0) {$\String_n$};
\node (D) at (6,0) {$\Spin_n$};
\node (E) at (8,0) {$1$};
\draw[->] (A) to (B);
\draw[->] (B) to (C);
\draw[->] (C) to (D);
\draw[->] (D) to (E);
\node (AA) at (0,1) {$1$};
\node (BB) at (2,1) {$\pt\sq U(1)^\delta$};
\node (CC) at (4,1) {$\String_n^\delta$};
\node (DD) at (6,1) {$\Spin_n^\delta$};
\node (EE) at (8,1) {$1$};
\draw[->] (AA) to (BB);
\draw[->] (BB) to (CC);
\draw[->] (CC) to (DD);
\draw[->] (DD) to (EE);
\draw[->] (BB) to (B);
\draw[->] (CC) to (C);
\draw[->] (DD) to (D);
\path (0,.5) coordinate (basepoint);
\end{tikzpicture}\label{diag:discretestring}
\eeq
where $\Spin_n^\delta$ and $U(1)^\delta$ are the respective groups with the discrete topology.
To spell out the map of extensions~\eqref{diag:discretestring},  by Lemma~\ref{lem:SMtoGRP} the Segal--Mitchison cocycle defining the lower extension determines (by precomposing with $\Spin_n^\delta\to \Spin_n$) a 3-cocycle $\CS^\delta\colon (\Spin_n^\delta)^{\times 3}\to U(1)$ classifying the upper extension. Here we use that any map $\CS^\delta\colon (\Spin_n^\delta)^{\times 3}\to U(1)$ is smooth; hence we may also take the discrete topology on $U(1)$. 

\begin{defn} 
The \emph{discrete string group} is the 2-group $\String_n^\delta$ in \eqref{diag:discretestring}.
\end{defn}

\section{Overview of 2-gerbes}\label{sec:2gerbe}

Bundle 2-gerbes are the natural categorification of the 2-stack of gerbes, which in turn are a categorification of the 1-stack of line bundles. Early work on 2-gerbes includes~\cite{Breen,CMW,Gajer,Stevenson}. Modern approaches using $\infty$-stacks have streamlined the theory considerably; a non-exhaustive list of references are~\cite[\S7.2.2]{LurieTopos}, \cite{NikolausSchreiberStevenson1,NikolausSchreiberStevenson2}, and \cite{Nuiten}. The full apparatus of $\infty$-stacks is both unnecessary and beyond the scope of this paper. For our purposes it suffices to define 2-gerbes in terms of \v{C}ech data analogous to that appearing in the definition of $A$-gerbes in Examples \ref{ex: gerbes} and \ref{ex: gerbes localized}. The abelian group $A$ will be omitted from the notation when it is clear from context.

\begin{defn}\label{def: 2-gerbe fixed base}
    Let $X$ be a Lie groupoid. A flat \emph{$A$-2-gerbe on $X$ that trivializes on $X_0$} is given by a locally constant \v Cech 3-cocycle $\lambda \colon X_3 \to A$. In particular, the \emph{trivial} 2-gerbe is given by the constant map $1_A \colon X_3 \to A$.

    A \emph{1-isomorphism} $\lambda_1 \to \lambda_2$ between a pair of such 2-gerbes is a locally constant 2-cochain $\gamma \colon X_2 \to A$ such that $d\gamma = \frac{\lambda_2}{\lambda_1}$: that is, for $(\phi, \psi, \theta) \in X_3$, 
    \begin{align*}
        \frac{\gamma(\psi, \theta), \gamma(\phi, \psi \circ \theta)}{\gamma(\phi \circ \psi, \theta)\gamma(\phi, \psi)} = \frac{\lambda_2(\phi, \psi, \theta)}{\lambda_1(\phi, \psi, \theta)}.
    \end{align*}
    In particular, a \emph{trivialization} of a 2-gerbe $\lambda$ that trivializes on $X_0$ is a 1-isomorphism $1_A \to \lambda$, i.e. a locally constant 2-cochain $\gamma \colon X_2 \to A$ such that $d\gamma = \lambda$. 

    A \emph{2-isomorphism} $\gamma_1 \to \gamma_2$ is a locally constant 1-cochain $\eta \colon X_1 \to A$ with $d\eta = \frac{\gamma_1}{\gamma_2}$:
    \begin{align*}
        \frac{\eta(\psi)\eta(\phi)}{\eta(\phi \circ \psi)} = \frac{\gamma_1(\phi, \psi)}{\gamma_2(\phi, \psi)} \; \text{for } (\phi, \psi) \in X_2.
    \end{align*} 

    A \emph{3-isomorphism} $\eta_1 \to \eta_2$ is a locally constant 0-cochain $\omega \colon X_0 \to A$ with $d\omega = \frac{\eta_2}{\eta_1}$:
    \begin{align*}
        \frac{\omega(x')}{\omega(x)} = \frac{\eta_2(\phi)}{\eta_1(\phi)} \; \text{ for } \phi \colon x \to x' \text{ in } X.
    \end{align*}

    This forms a strict 3-category, with composition induced from multiplication in $A$. We denote this 3-category by $2\Gerbe_A^\pre(X)$. 
\end{defn}

\begin{prop}\label{prop: pullback magic for 2-gerbes}
    The assignment $X \mapsto 2\Gerbe_A^\pre(X)$ forms a prestack.
\end{prop}

\begin{proof}
It is clear that a smooth functor $f \colon X \to Y$ induces a pullback $f^* \colon 2\Gerbe_A^\pre(Y) \to 2\Gerbe_A^\pre(X)$ by precomposition of cochains, analogously to Definition~\ref{def: pullback}. 

It remains to show that a natural transformation $\tau \colon f_1 \Rightarrow f_2$ of smooth functors induces a natural transformation of pullbacks $\tau^* \colon f_1^* \Rightarrow f_2^*$. Let $\lambda$ be an object in $2\Gerbe_A^\pre(Y)$. We define a 1-morphism $ \tau^*_\lambda \colon f_1^*(\lambda) \to f_2^*(\lambda)$ by 
    \begin{align*}
        \tau^*_\lambda (\phi, \psi) = \frac{\lambda(f_2(\phi), \tau(x'), f_1(\psi)}{\lambda(f_2(\phi), f_2(\psi), \tau(x)) \lambda(\tau(x''), f_1(\phi), f_1(\psi))} \; \text{ for } x \xrightarrow{\psi} x' \xrightarrow{\phi} x'' \text{ in } X. 
    \end{align*}
    The 3-cocycle condition on $\lambda$ implies that $d\lambda (\phi, \psi, \theta)= \frac{\lambda(f_2(\phi), f_2(\psi), f_2(\theta))}{\lambda(f_1(\phi), f_1(\psi), f_1(\theta))}$ for all $(\phi, \psi, \theta)$ in $X_3$, as required. 
    It is natural in $\lambda$ and hence defines the required natural transformation. 
\end{proof}

With this established, we adopt the following definition of 2-gerbe, e.g., see \cite[Definition 4.48]{NikolausSchreiberStevenson1}. 

\begin{defn} 
Let $X$ be a Lie groupoid. A \emph{flat 2-gerbe on $X$} is an essential equivalence $ f \colon Y \xrightarrow[]{\sim} X$ and a locally constant 3-cocycle $\lambda \colon Y_3 \to A$. A \emph{1-isomorphism} $(f_1 \colon Y_1 \to X, \lambda_1) \to (f_2 \colon Y_2 \to X, \lambda_2)$ of 2-gerbes is given by a tuple $(Z,g_1,g_2,\tau, \gamma)$, where $Z$ is a Lie groupoid, $g_1 \colon Z \to Y_1$ and $g_2 \colon Z \to Y_2$ are smooth functors, $\tau$ is a smooth natural transformation $f _1\circ g_1 \to f_2 \circ g_2$, and $\gamma \colon Z_2 \to A$ gives a 1-isomorphism $g_1^*(\lambda_1) \to g_2^* (\lambda_2)$. That is, for $(\phi, \psi, \theta) \in Z_3$, we have 
    \begin{align*}
        \frac{\gamma(\psi, \theta)\gamma(\phi, \psi \circ \theta)}{\gamma(\phi \circ \psi, \theta) \gamma(\phi, \psi)} = \frac{\lambda_2(g_2(\phi), g_2(\psi), g_2(\theta))}{\lambda_1(g_1(\phi), g_1(\psi), g_1(\theta))}.
    \end{align*}
Similarly, we have 2- and 3-morphisms between 2-gerbes defined analogously using Definition~\ref{def: 2-gerbe fixed base} together with refinements of essential equivalences $Y\xrightarrow[]{\sim}X$. The collection of 2-gerbes on $X$ forms a $(3,1)$-category $2\Gerbe_A(X)$.
\end{defn}

\begin{rmk}
More formally, the 3-category $2\Gerbe_A(X)$ is the stackification of the 3-category $2\Gerbe_A^\pre(X)$. This stackification arises from a higher-categorical localization of $2\Gerbe_A^\pre(X)$, e.g., see \cite{Nuiten} for an exposition that generalizes~\cite{Pronk}. 
\end{rmk}
In light of the previous remark, morphisms between a pair of 2-gerbes form a (possibly empty) bicategory. In particular, trivializations of a given flat 2-gerbe $(f \colon Y \to X, \lambda)$ are the objects of the (possibly empty) bicategory ${\rm Triv}(f, \lambda)$. 
This bicategory itself can be identified with the 2-stackification of the bicategory of trivializations of $\lambda$ in $2\Gerbe_A^\pre(Y)$. Up to isomorphism, the objects of ${\rm Triv}(f, \lambda)$  are of the form $(g \colon Z \to Y, \gamma \colon 1_A \to g^*\lambda)$. 

\begin{lem}\label{lem: trivializations of 2-gerbe torsor for gerbes}
    Let $(f \colon Y \to X, \lambda)$ be a flat 2-gerbe on $X$. If it is non-empty, the bicategory ${\rm Triv}(f, \lambda)$ of trivializations of $(f, \lambda)$ is non-canonically equivalent to the bicategory $\Gerbe_A(X)$, where a choice of equivalence  ${\rm Triv}(f, \lambda) \simeq \Gerbe_A(X)$ is determined by a choice of trivialization $(g, \gamma)$ of~$(f, \lambda)$. 
\end{lem}

\begin{proof}
    It suffices to show that the choice of trivialization $\gamma$ of $g^*\lambda$ determines an equivalence (in fact an isomorphism) of the bicategory of trivializations of $g^*\lambda$ in $2\Gerbe_A^\pre(Z)$ with the bicategory $\Gerbe_A^\pre(Z)$; upon 2-stackification this yields the desired equivalence. 

    Let $\gamma_1$ be any other trivialization of $g^*\lambda$; then as in definition \ref{def: 2-gerbe fixed base} we have $d\gamma = d\gamma_1 = g^* \lambda$, so $d\frac{\gamma_1}{\gamma}= 1$, and $\frac{\gamma_1}{\gamma}$ is an object of $\Gerbe_A^\pre(Z)$. This gives the bijection on the level of objects. 

    If $\gamma_1, \gamma_2$ are two trivializations of $g^*\lambda$ and $\eta$ is a 2-morphism between them, then we have $\frac{\eta(\psi)\eta(\phi)}{\eta(\phi \circ \psi)} = \frac{\gamma_1}{\gamma_2} = \frac{\gamma_1/\gamma}{\gamma_2/\gamma}$. So $\eta$ is also a 1-morphism between the gerbes $\gamma_1/\gamma$ and $\gamma_2/\gamma$.

    Finally, if $\eta_1, \eta_2$ are two 2-morphisms between trivializations of $g^*\lambda$, and $\omega$ is a 3-morphism between them, we have $d\omega = \frac{\eta_2}{\eta_1}$, which implies that $\omega$ is also a 2-morphism between the $\eta_i$'s viewed as 1-morphisms in $\Gerbe_A^\pre(Z)$. 

    This completes the proof. 
\end{proof}

\section{Proof of Theorem~\ref{thm:pronk}}\label{sec: localization app}

Recall that Theorem~\ref{thm:pronk} states that the bicategory $\Bun_\cG^\pre$ admits a calculus of right fractions with respect to the class of 1-morphisms $\calw$ of the form $(f,h,\eta)$ where $f$ is an essential equivalence. Proving this requires verification of conditions BF1-BF5 from \cite[\S 2.1]{Pronk}; we take the equivalent (diagrammatic) version of BF4 from \cite{vazquez2021fs}.
Performing this localization yields the bicategory $\Bun_\cG:=\Bun_\cG^\pre[\calw^{-1}]$, the stack of $\cG$-bundles on Lie groupoids.

\begin{lem}\label{lem: pronk1}
    Let $(f, h, \eta): (\rho_1, \gamma_1)_X\to (\rho_2, \gamma_2)_Y$ be a 1-morphism in $\Bun_\cG^\pre$. If $(f, h, \eta)$ is an equivalence, then $f$ is an essential equivalence. This is \cite[BF1]{Pronk}.
\end{lem}

\bp
If $(f, h, \eta)$ is an equivalence, there exists $(f_2, h_2, \eta_2)$ such that $(f, h,\eta)\circ(f_2, h_2,\eta_2)\cong id_{(\rho_2,\gamma_2)_Y}$ and $(f_2,h_2,\eta_2)_Y \circ(f,h,\eta)\cong id_{(\rho_1 ,\gamma_1 )_X}$; in particular $f \circ f_2 \cong \id_Y$ and $f_2 \circ f \cong \id_X$. Hence $f$ is an equivalence of Lie groupoids and 
(in particular) an essential equivalence. 
\ep

\begin{lem}\label{lem: pronk2}
   Given $(f_4, h_4, \eta_4): (\rho_1,\gamma_1)_X\to (\rho_3,\gamma_3)_Z, (f_3, h_3, \eta_3): (\rho_2,\gamma_2)_Y\to (\rho_3,\gamma_3)_Z$ in $\Bun_\cG^\pre$, where $(f_4, h_4, \eta_4)$ is in $\calw$, there exists the data filling the diagram below:
\beq\label{eq:BF3square}
       \begin{tikzcd}[column sep=huge]
           (\rho_0,\gamma_0)_V \arrow{r}{(f_1,h_1,\eta_1)} \arrow[swap, name=U]{d}{(f_2, h_2, \eta_2)} & (\rho_2, \gamma_2)_Y \arrow{d}{(f_3,h_3,\eta_3)}\\
           (\rho_1, \gamma_1)_X \arrow[swap]{r}{(f_4,h_4,\eta_4)} \ar[ru,shorten <>=20pt,Rightarrow, "\cong"']{}{(\tau,\omega)}& (\rho_3,\gamma_3)_Z
\end{tikzcd}
\eeq
    where $(f_1,h_1,\eta_1)$ is in $\calw$ and $(\tau, \omega)$ is a 2-isomorphism. This is BF3 of \cite{Pronk}.
\end{lem}

\bp
The desired data in terms of the Lie groupoid structure is given by the lax pullback $V$ of $X$ and $Y$ over $Z$, see \cite{Moerdijk88}. 
This pullback has objects $V_0:= X_0\times_{Z_0} Z_1\times_{Z_0} Y_0$ and morphisms $V_1:= eq\Big( X_1\times(V_0\times V_0)\times Y_1 \rightrightarrows Z_1\Big)$. Unpacking this, objects are of the form $(x,\phi, y)$ where $\phi:f_4(x)\to f_3(y)$ in $Z_1$. A morphism $\Phi$ from $(x,\phi, y)$ to $(x',\phi', y')$ consists of $\Phi_X\in X_1, \Phi_Y\in Y_1$ such that $f_3(\Phi_Y)\circ \phi =\phi'\circ f_4(\Phi_X)$ in $Z_1$.

We denote the natural projections from $V$ onto $Y$ and $X$ by $f_1, f_2$ respectively. The natural transformation $\tau \colon f_4 \circ f_2 \to f_3 \circ f_1$ in the diagram~\eqref{eq:BF3square} is given by $V_0 \to Z_1$, $(x, \phi, y) \mapsto \phi$. The equality $f_3(\phi_Y)\circ \phi =\phi'\circ f_4(\phi_X)$ verifies the naturality condition for $\tau$.

By Lemma~\ref{lem: factor 1-morph}, the lower horizontal arrow in~\eqref{eq:BF3square} factors as
\beq\label{eq:factorBF3}
    \begin{tikzcd}[column sep=huge]
            &  & (\rho_2, \gamma_2)_Y \arrow{d}{(f_3,h_3,\eta_3)}\\
           (\rho_1, \gamma_1)_X \arrow[swap]{r}{(1,h_4,\eta_4)} & (f_4^*(\rho_3,\gamma_3))_X \arrow[swap]{r}{(f_4,1,1)}& (\rho_3,\gamma_3)_Z .
    \end{tikzcd}
\eeq
We will construct the square~\eqref{eq:BF3square} in two steps, via squares over each of the horizontal arrows in~\eqref{eq:factorBF3}. First take the square
\beq\label{eq:BF3square1}
    \begin{tikzcd}[column sep=huge]
           (f_1^*(\rho_2,\gamma_2))_V \arrow{r}{(f_1,1,1)} \arrow[swap, name=U]{d}{(f_2, h, \eta)} & (\rho_2, \gamma_2)_Y \arrow{d}{(f_3,h_3,\eta_3)}\\
           (f_4^*(\rho_3,\gamma_3))_X \arrow[swap]{r}{(f_4,1,1)}& (\rho_3,\gamma_3)_Z
       \end{tikzcd}
\eeq
where $(h,\eta)$ is the composition of the 1-morphisms  
$$f_1^*(\rho_2,\gamma_2)\xrightarrow{f_1^*(h_3,\eta_3)}f_1^*f_3^*(\rho_3,\gamma_3)\xrightarrow{(h_\tau,\eta_\tau)}f_2^*f_4^*(\rho_3,\gamma_3)$$
for $(h_\tau,\eta_\tau)$ the value of the natural transformation $\tilde{\tau}\colon (f_3\circ f_1)^*\Rightarrow (f_4\circ f_2)^*$ on $(\rho_3,\gamma_3)$, see Proposition~\ref{prop:pullbackmagic}. The 2-isomorphism filling~\eqref{eq:BF3square1} square is $(\tau, \omega:=1_A)$. To construct the remaining square whose lower horizontal arrow is the leftmost arrow in~\eqref{eq:factorBF3}, consider
\begin{center}
    \begin{tikzcd}[column sep=huge]
             & (f_1^*(\rho_2,\gamma_2))_V \arrow{d}{(f_2,h,\eta)}\\
            (\rho_1,\gamma_1)_X \arrow[swap]{r}{(1, h_4,\eta_4)}& (f_4^*(\rho_3,\gamma_3))_X.
       \end{tikzcd}
\end{center}
As $(1,h_4,\eta_4)$ admits an inverse, consider the composition
$$f_1^*(\rho_2,\gamma_2)\xrightarrow{(f_2,h,\eta)} f_4^*(\rho_3,\gamma_3)\xrightarrow{(1,h_4,\eta_4)^{-1}} (\rho_1,\gamma_1),$$
which admits a factoring as in Lemma~\ref{lem: factor 1-morph}
\begin{center}
    \begin{tikzcd}
        f_1^*(\rho_2,\gamma_2) \arrow{rr} \arrow[swap]{dr}{(1,h_1,\eta_1)^{-1}} && (\rho_1,\gamma_1).\\
        & f_2^*(\rho_1,\gamma_1)\arrow[swap]{ur}{(f_2, 1, 1)} &
    \end{tikzcd}
\end{center}
Taking $(\rho_0,\gamma_0)_V:=(f_2^*(\rho_1,\gamma_1))_V$, we obtain 
\begin{center}
    \begin{tikzcd}[column sep=huge]
            (\rho_0,\gamma_0)_V \arrow{r}{(1, h_1, \eta_1)} \arrow[swap, name=U]{d}{(f_2, 1, 1)} & (f_1^*(\rho_2,\gamma_2))_V \arrow{d}{(f_2,h,\eta)}\\
            (\rho_1,\gamma_1)_X \arrow[swap]{r}{(1, h_4,\eta_4)}& (f_4^*(\rho_3,\gamma_3))_X
       \end{tikzcd}
\end{center}
where the 2-commutativity data comes from the uniquely specified 2-morphism $(1,h_1,\eta_1)\simeq ((1,h_1,\eta_1)^{-1})^{-1}$ determined by~\eqref{eq:unitcounitnormalized}.
Composing these two squares gives the desired data for the square in BF3 of \cite{Pronk}.
\ep 

\begin{lem} \label{lem: pronk3}
    Suppose the data in the diagram below (in $\Bun_\cG^\pre$)
    \begin{center}
        \begin{tikzcd}[column sep=huge]
            (\rho_1,\gamma_1)_X \arrow[bend left=40]{r}{(f_1, h_1,\eta_1)}\arrow[swap, bend right=40]{r}{(f_2,h_2,\eta_2)} & (\rho_2,\gamma_2)_Y\arrow{r}{(f_3,h_3,\eta_3)} &(\rho_3,\gamma_3)_Z,\\
        \end{tikzcd}  
    \end{center}
     where $(f_3,h_3,\eta_3)$ is in $\calw$, and where there exists a 2-morphism $(\tau,\omega): (f_3,h_3,\eta_3)\circ(f_1,h_1,\eta_1)\Rightarrow (f_3,h_3,\eta_3)\circ(f_2,h_2,\eta_2)$. Then there exists a 2-morphism $(\tilde{\tau},\tilde{\omega}): (f_1,h_1,\eta_1)\Rightarrow(f_2,h_2,\eta_2)$. This is the first part of BF4 of \cite{Pronk}; alternatively it is condition 1-Frc of \cite{vazquez2021fs}
\end{lem}

\bp
Note that $\tau$ gives a map $\tau:X_0\to Z_1$. Because $f_3$ is an essential equivalence, see \eqref{eq:essentialequiv}, we have the following pullback diagram
\begin{center}
    \begin{tikzcd}
        X_0 \arrow[bend left=20]{rrd}{\tau}\arrow[swap, bend right=20]{rdd}{f_1\times f_2} \arrow[dashed]{dr}{\tilde\tau}& &\\ 
        & Y_1 \arrow[r]\arrow[d] & Z_1 \arrow[d]\\
        & Y_0\times Y_0 \arrow[r] & Z_0\times Z_0.\\
    \end{tikzcd}
\end{center}
 For the naturality conditions $\tilde\tau$ satisfies, see~\cite{Pronk} Section 4. 

Define $\tilde\omega:=\omega: X_0\to A$. To check that this satisfies the compatibility conditions needed to define the desired 2-morphism $(\tilde\tau,\tilde\omega)$ (from \eqref{compat for 2-morph}) note that the conditions to have the original 2-morphism $(\tau,\omega)$ give that 
\begin{align*}
    h_3(f_2(x))h_2(x) &= \rho_3(\tau(x))h_3(f_1(x))h_1(x) \\
    &= h_3(f_2(x))\rho_2(\tau(x))h_1(x)\\
\end{align*}
So $\rho_2(\tau(x))h_1(x)=h_2(x)$. 
Additionally,
\begin{align*} 
  \frac{\tilde{\omega}(x)}{\tilde{\omega}(x')}&= \frac{\omega(x)}{\omega(x')} 
  = \frac{\eta_{31}(\phi)}{\eta_{32}(\phi)}\frac{\gamma_3(\tau(x'), f_{31}(\phi))}{\gamma_3(f_{32}(\phi), \tau(x)} \\[5pt]
  & \phantom{={}}\frac{\alpha(\rho_3(\tau(x')), \rho_3(f_{31}(\phi)), h_3(f_1(x))h_1(x))}{\alpha(\rho_3(\tau(x')), h_3(f_1(x'))h_1(x'), \rho_1(\phi))\alpha(\rho_3(f_{32}(\phi)), \rho_3(\tau(x)), h_3(f_1(x))h_1(x)}\\[5pt]
  &= \frac{\eta_1(\phi)}{\eta_2(\phi)}\frac{\gamma_2(\tilde\tau(x'),f_1(\phi))}{\gamma_2(f_2(\phi),\tilde\tau(x))}\\[5pt]
  &\phantom{={}} \frac{\alpha (\rho_2(\tilde{\tau}(x')),\rho_2(f_1(\phi)), h_1(x))}{\alpha(\rho_2(\tilde{\tau}(x')), h_1(x'), \rho_1(\phi))\alpha(\rho_2(f_2(\phi)), \rho_2(\tilde{\tau}(x)), h_1(x))}
\end{align*}
for $\phi \colon x \to x'$ in $X$. 
This calculation involves repeated application of the cocycle condition for $\alpha$ and the conditions for $\rho_i, h_i, i=1,2,3$. 
This shows that $\tilde{\omega}:=\omega: X_0\to A$ satisfies the necessarily compatibility conditions to be a 2-morphism between the desired 1-morphisms.
\ep 

\begin{lem}\label{lem: pronk4}
    Given the data in the diagram below (in $\Bun_\cG^\pre$),
    \begin{center}
        \begin{tikzcd}[column sep=10em]
            (\rho_1,\gamma_1)_X \arrow[bend left=40, ""{name=U, below}]{r}{(f_1, h_1,\eta_1)}\arrow[swap, bend right=40, ""{name=D}]{r}{(f_2,h_2,\eta_2)} & (\rho_2,\gamma_2)_Y\arrow{r}{(f_3,h_3,\eta_3)} &(\rho_3,\gamma_3)_Z\\
            \arrow[Rightarrow, bend right=40, from=U, to=D, shorten <>=10pt,swap]{}{(\tau_1,\omega_1)}
            \arrow[Rightarrow, bend left=40, from=U, to=D, shorten <>=10pt]{}{(\tau_2,\omega_2)}
        \end{tikzcd}  
    \end{center}
    where $f_3$ is an essential equivalence and $(f_3,h_3,\eta_3)\star (\tau_1,\omega_1)=(f_3,h_3,\eta_3)\star (\tau_2, \omega_2)$, then $(\tau_1,\omega_1)=(\tau_2,\omega_2)$. This is the second part of BF4 of \cite{Pronk}; alternatively it is condition 2-Frc of \cite{vazquez2021fs}.
\end{lem}

\bp
Because $f_3$ is an essential equivalence, we obtain the pullback square
\begin{center}
   \begin{tikzcd}
        X_0 \arrow[bend left=20]{rrd}{f_3\star\tau_1=f_3\star\tau_2}\arrow[swap, bend right=20]{rdd}{f_1\times f_2} \arrow[dashed]{dr}{\tau_1=\tau_2}& &\\ 
        & Y_1 \arrow{r}{f_3}\arrow{d}{s\times t} & Z_1 \arrow{d}{s\times t}\\
        & Y_0\times Y_0 \arrow{r}{f_3\times f_3} & Z_0\times Z_0\\
    \end{tikzcd}
\end{center}
By the universal property, the dashed arrow from $X_0\to Y_1$ is unique; both $\tau_1, \tau_2$ fill that arrow, so $\tau_1=\tau_2: X_0\to Y_1$.

Let us denote $\id_{(f_3, h_3, \eta_3)}$ by $(\tau_3, \omega_3)$, so $\tau_3:Y_0\to Z_1$ sends $y\mapsto \id: f_3(y)\to f_3(y)$, and $\omega_3:Y_0\to A$ is the constant map to $1_A$. From the equation for horizontal composition given in~\eqref{eq: horiz comp 2-morph} and using the fact that $\alpha$ is normalized, 
\begin{align*}
    \omega_3 \star \omega_1(x) = & \omega_3(f_1(x))\omega_1(x)\eta_3(\tau_1(x)) \gamma_3(f_3(\tau_1(x)), \tau_2(f_1(x))) \\ 
 & \phantom{={}} \frac{\alpha(\rho_3(f_3(\tau_1(x))), \rho_3(\tau_3(f_1(x))), h_2(f_1(x))h_1(x))}{\alpha(\rho_3(\tau_3(f_1(x))), h_2(f_1(x)), h_1(x))}\\[5pt]
 &\phantom{={}} \frac{\alpha(h_3(f_2(x)), \rho_2(\tau_1(x)), h_1(x))}{\alpha(\rho_3(f_3(\tau_1(x))), h_3(f_1(x)), h_1(x))}\\[5pt]
 &= \omega_1 (x) \eta_3(\tau_1(x)) \gamma_3(f_3(\tau_1(x)), \tau_1(f_1(x))) \\[5pt]
 & \phantom{={}} \frac{\alpha(\rho_3(f_3(\tau_1(x))), 1_G, h_2(f_1(x))h_1(x))\alpha(h_3(f_2(x)), \rho'(\tau_1(x)), h_1(x))}{\alpha(1_G, h_2(f_1(x)), h_1(x))\alpha(\rho_3(f_3(\tau_1(x))), h_3(f_1(x)), h_1(x))}\\[5pt]
 &= \omega_1 (x) \eta_3(\tau_1(x)) \gamma_3(f_3(\tau_1(x)), \tau_3(f_1(x)))\\[5pt]
 & \phantom{={}}   \frac{\alpha(h_3(f_2(x)), \rho_2(\tau_1(x)), h_1(x))}{\alpha(\rho_3(f_3(\tau_1(x))), h_3(f_1(x)), h_1(x))}\\[2em]
\omega_3\star \omega_2(x) = & \omega_3(f_1(x))\omega_2(x)\eta_3(\tau_2(x)) \gamma_3(f_3(\tau_2(x)), \tau_3(f_1(x))) \\[5pt]
 & \phantom{={}} \frac{\alpha(\rho_3(f_3(\tau_2(x))), \rho_3(\tau_3(f_1(x))), h_2(f_1(x))h_1(x))}{\alpha(\rho_3(\tau_3(f_1(x))), h_2(f_1(x)), h_1(x))}\\[5pt]
 &\phantom{={}} \frac{\alpha(h_3(f_2(x)), \rho_2(\tau_2(x)), h_1(x))}{\alpha(\rho_3(f_3(\tau_2(x))), h_3(f_1(x)), h_1(x))}\\[5pt]
 &= \omega_2 (x) \eta_3(\tau_2(x)) \gamma_3(f_3(\tau_2(x)), \tau_3(f_1(x)))\\[5pt]
 &\phantom{={}} \frac{\alpha(h_3(f_2(x)), \rho_2(\tau_2(x)), h_1(x))}{\alpha(\rho_3(f_3(\tau_2(x))), h_3(f_1(x)), h_1(x))}\\[5pt]
\end{align*}
Since $\omega_3\star \omega_1=\omega_3\star \omega_2$ and $\tau_1=\tau_2$, the remaining terms in the above equations cancel out, giving $\omega_1=\omega_2: X_0\to A$.
\ep

 \begin{proof}[Proof of Theorem~\ref{thm:pronk}]
Note that the composition of 1-morphisms in $\calw$ is again in $\calw$; this is condition BF2 of \cite{Pronk}. If there is a 2-isomorphism between two 1-morphisms where one of the 1-morphisms in is $\calw$, then the other 1-morphism is also in $\calw$; this is condition BF5 of \cite{Pronk}. Both these conditions are satisfied trivially. Together with lemmas \ref{lem: pronk1}, \ref{lem: pronk2}, \ref{lem: pronk3},  and \ref{lem: pronk4}, this completes the proof. 
\ep

\bibliographystyle{amsplain}

\begin{thebibliography}{1}
\bibitem{baas_dundas_rognes_2004}
N.~A. Baas, B.~I. Dundas, and J.~Rognes, \emph{Two-vector bundles and forms of
  elliptic cohomology}, London Mathematical Society Lecture Note Series,
  p.~18–45, Cambridge University Press, 2004.

\bibitem{BL04}
J.~Baez and A.~Lauda, \emph{Higher-dimensional algebra. {V}. 2-groups}, Theory
  Appl. Categ. \textbf{12} (2004).

\bibitem{BaezSchreiber}
J.~Baez and U.~Schreiber, \emph{Higher gauge theory}, Categories in Algebra,
  Geometry and Mathematical Physics, eds. A. Davydov et al, Contemp. Math.
  \textbf{431} (2007).

\bibitem{BaezStevenson}
J.~Baez and D.~Stevenson, \emph{The classifying space of a topological
  2-group}, Algebraic Topology Abel Symposia \textbf{4} (2009).

\bibitem{BCSS}
J.~Baez, A.~Crans, D.~Stevenson, and U.~Schreiber, \emph{{From
  loop groups to 2-groups}}, Homology, Homotopy and Applications \textbf{9}
  (2007), no.~2, 101 -- 135.

\bibitem{Bartels}
T.~Bartels, \emph{Higher gauge theory {I}: 2-bundles}, preprint (2004).

\bibitem{BehrendXu}
K.~Behrend and P.~Xu, \emph{Differentiable stacks and gerbes}, Journal of
  Symplectic Geometry \textbf{9} (2011).

\bibitem{Breen}
L.~Breen, \emph{On the classification of 2-gerbes and 2-stacks}, Ast\`erisque
  \textbf{225} (1994).

\bibitem{Bunk}
S.~Bunk, \emph{Principal $\infty$-bundles and smooth string group models},
  Hamburger Beitraege Nr. \textbf{858} (2020).

\bibitem{CMW}
A.~Carey, M.~Murray, and B.L. Wang, \emph{Higher bundle gerbes and cohomology
  classes in gauge theories}, Journal of Geometry and Physics \textbf{21}
  (1997).

\bibitem{Costello1}
K.~Costello, \emph{A geometric construction of the {Witten} genus {I}},
  Proceedings of the International Congress of Mathematicians 2010 (2011).

\bibitem{Costello2}
\bysame, \emph{A geometric construction of the {Witten} genus {II}}, ArXiv
  preprint, \url{http://arxiv.org/abs/1112.0816} (2011).

\bibitem{Dekimpe}
K.~Dekimpe, M.~Sadowski, and A.~Szczepański, \emph{Spin structures on flat
  manifolds}, Mh Math \textbf{148} (2006), 283--296.

\bibitem{DouglasHenriques}
C.~Douglas and A.~Henriques, \emph{Topological modular forms and conformal
  nets}, Mathematical Foundations of QFT and Perturbative String Theory, Proc.
  Sympos. Pure Math. \textbf{83} (2011).

\bibitem{vienna}
D.~Freed, \emph{Lectures on twisted {K}-theory and orientifolds}, lectures at
  {ESI} {Vienna} (2012).

\bibitem{Gajer}
P.~Gajer, \emph{Geometry of {Deligne} cohomology}, Invent. Math. \textbf{127}
  (1997).

\bibitem{Andre}
A.~Henriques, \emph{Integrating $l_\infty$-algebras}, Compositio Mathematica
  \textbf{144} (2008).

\bibitem{HU2004325}
P.~Hu and I.~Kriz, \emph{Conformal field theory and elliptic cohomology},
  Advances in Mathematics \textbf{189} (2004), no.~2, 325--412.

\bibitem{JohnsonYau}
N.~Johnson and D.~Yau, \emph{2-dimensional categories}, Oxford University
  Press, 2021.

\bibitem{Killingback}
T.~P. Killingback, \emph{World-sheet anomalies and loop geometry}, Nuclear
  Phys. B \textbf{288} (1987).

\bibitem{Lerman}
E.~Lerman, \emph{Orbifolds as stacks?}, Enseign. Math. \textbf{2} (2010).

\bibitem{LurieTopos}
J.~Lurie, \emph{Higher topos theory}, Annals of Mathematics Studies 170,
  Princeton University Press, 2009.

\bibitem{Lurie}
J.~Lurie, \emph{A survey of elliptic cohomology}, Algebraic Topology (N.~Baas,
  E.~Friedlander, J.~Bj\"orn, and P.~{O}st\ae{r}, eds.), vol.~4, Springer
  Berlin Heidelberg, 2009.

\bibitem{LUTOWSKI}
R.~Lutowski and B.~Putrycz, \emph{Spin structures on flat manifolds}, Journal
  of Algebra \textbf{436} (2015), 277--291.

\bibitem{Mackenzie}
K. Mackenzie, \emph{General theory of Lie groupoids and Lie algebroids},
  Cambridge University Press
  (2005).


\bibitem{Moerdijk88}
I. Moerdijk, \emph{The classifying topos of a continuous groupoid. i},
  Transactions of the American Mathematical Society \textbf{310} (1988), no.~2,
  629--668.

\bibitem{Moerdijk2002}
\bysame, \emph{Orbifolds as groupoids: an introduction}, ArXiv
  preprint, \url{https://arxiv.org/pdf/math/0203100}
  (2002).

\bibitem{Murray}
M.~K. Murray, \emph{Bundle gerbes}, Journal of the London Mathematical Society
  \textbf{54} (1996), no.~2, 403--416.

\bibitem{NSW}
T.~Nikolaus, C.~Sachse, and C.~Wockel, \emph{A smooth model for the string
  group}, International Mathematics Research Notices \textbf{2013} (2013).

\bibitem{NikolausSchreiberStevenson1}
T.~Nikolaus, U.~Schreiber, and D.~Stevenson, \emph{Principal $\infty$-bundles:
  {General} theory}, Journal of Homotopy and Related Structures \textbf{10}
  (2015).

\bibitem{NikolausSchreiberStevenson2}
\bysame, \emph{Principal $\infty$-bundles: {Presentations}}, Journal of
  Homotopy and Related Structures \textbf{10} (2015).

\bibitem{NikolausSchweigert}
T.~Nikolaus and C.~Schweigert, \emph{Equivariance in higher geometry}, Advances
  in Mathematics \textbf{226} (2011).

\bibitem{NikolausWaldorf}
T.~Nikolaus and K.~Waldorf, \emph{Four equivalent versions of non-abelian
  gerbes}, Pacific J. Math. \textbf{264} (2013).

\bibitem{Nuiten}
J.~Nuiten, \emph{Higher stacks as a category of fractions}, preprint (2016).

\bibitem{PFAFFLE}
F.~Pfaffle, \emph{{The Dirac spectrum of Bieberbach manifolds}}, Journal of
  Geometry and Physics \textbf{35} (2000), no.~4, 367--385.

\bibitem{Pronk}
D.~Pronk, \emph{Etendues and stacks as bicategories of fractions}, Compositio
  Math. \textbf{102} (1996).

\bibitem{SP11}
C.~J. Schommer-Pries, \emph{Central extensions of smooth 2-groups and a
  finite-dimensional string 2-group}, Geometry \& Topology \textbf{15} (2011).

\bibitem{SchreiberWaldorf}
U.~Schreiber and K.~Waldorf, \emph{Smooth functors vs. differential forms},
  Homology Homotopy Appl. \textbf{13} (2011).

\bibitem{Segal_Elliptic}
G.~Segal, \emph{Elliptic cohomology}, S\'eminaire N. Bourbaki \textbf{695}
  (1988).

\bibitem{SegalCFT}
\bysame, \emph{The definition of conformal field theory}, Topology, geometry
  and quantum field theory, London Math. Soc. Lecture Note Ser., vol. 308,
  Cambridge Univ. Press, Cambridge, 2004, pp.~421--577.

\bibitem{sharpetorsionhet}
E.~Sharpe, \emph{Discrete torsion in perturbative heterotic string theory},
  Phys. Rev. D \textbf{68} (2003).

\bibitem{sinh}
H.~Sinh, \emph{Gr-categories}, Universit\'e {Paris VII} \textbf{doctoral
  thesis} (1975).

\bibitem{Stevenson}
D.~Stevenson, \emph{Bundle 2-gerbes}, Proc. Lond. Math. Soc. (2004).

\bibitem{Stolz}
S.~Stolz, \emph{A conjecture concerning positive ricci curvature and the witten
  genus}, Mathematische Annalen \textbf{304} (1996).

\bibitem{ST04}
S.~Stolz and P.~Teichner, \emph{What is an elliptic object?}, Topology,
  geometry and quantum field theory, London Math. Soc. LNS 308, Cambridge Univ.
  Press (2004), 247--343.

\bibitem{ST11}
S.~Stolz and P.~Teichner, \emph{Supersymmetric field theories and generalized
  cohomology}, Mathematical foundations of quantum field theory and
  perturbative string theory \textbf{83} (2011), 279--340.

\bibitem{vazquez2021fs}
P.~Bustillo Vazquez, D.~Pronk, and M.~Szyld, \emph{The three f's for
  bicategories i: Localization by fractions is exact}, 2021.

\bibitem{Waldorf}
K.~Waldorf, \emph{A construction of string 2-group models using a
  transgression-regression technique}, Analysis, geometry and quantum field
  theory \textbf{584} (2012).

\bibitem{Waldorfstring}
\bysame, \emph{String connections and {Chern--Simons} theory}, Transactions of
  the American Mathematical Society \textbf{365} (2013).

\bibitem{WaldorfAGT}
\bysame, \emph{Spin structures on loop spaces that characterize string
  manifolds}, Algebr. Geom. Topol. \textbf{16} (2016).

\bibitem{Whitehead1941}
J.H.C.~Whitehead,
  \emph{On adding relations to homotopy groups},
Annals of Mathematics (1941).

\bibitem{Whitehead1949}
\bysame,
  \emph{Combinatorial homotopy. II}, Bulletin of the American Mathematical Society (1949).

\bibitem{WittenDirac}
E.~Witten, \emph{The index of the {Dirac} operator in loop space}, In Elliptic
  curves and modular forms in algebraic topology (Princeton, NJ, 1986), volume
  1326 of Lecture Notes in Math. (1988).

\bibitem{Wockel}
C.~Wockel, \emph{A global perspective to gerbes and their gauge stacks}, Forum
  Math. \textbf{23} (2011).

\bibitem{Killingback}
T.P. Killingback, \emph{World-sheet anomalies and loop geometry}, Nuclear Physics B Volume 288, Pages 578-588 (1987).

\bibitem{Bunke}
U. Bunke, \emph{String structures and trivialisations of a Pfaffian line bundle}, Commun. Math. Phys. 307 (2011).

\end{thebibliography}

\end{document}